\documentclass{elsarticle}
\usepackage{algorithm}
\usepackage{algorithmic}
\usepackage{amsfonts, amsmath, amsopn, amssymb, amsthm}
\usepackage{url}
\usepackage{bm}
\usepackage{datetime2}
\usepackage[inline, shortlabels]{enumitem}
\usepackage{graphicx}
\usepackage[section]{placeins}
\usepackage{xspace}

\newtheorem{theorem}{Theorem}

\newdefinition{rmk}{Remark}
\newproof{pf}{Proof}
\newproof{pot}{Proof of Theorem \ref{thm2}}

\usepackage{upquote,textcomp}

\usepackage{bm}
\usepackage{mathrsfs}
\usepackage{dsfont}
\usepackage{color}

\newcounter{todocounter}

\usepackage{textcomp}

\newcommand{\diff}[1]{\textcolor{red}{#1}}

\newcommand{\chol}{\texttt{chol}\xspace}
\newcommand{\triu}{\texttt{triu}\xspace}

\newcommand{\nan}{\texttt{NaN}\xspace}

\newcommand{\spR}{\mathbb{R}}				

\newcommand{\s}{$s$\xspace}

\newcommand{\vq}{\bm{q}}

\newcommand{\vr}{\bm{r}}

\newcommand{\vt}{\bm{t}}
\newcommand{\vu}{\bm{u}}
\newcommand{\vv}{\bm{v}}

\newcommand{\vw}{\bm{w}}
\newcommand{\vx}{\bm{x}}

\newcommand{\vz}{\bm{z}}

\newcommand{\vB}{\bm{B}}

\newcommand{\vD}{\bm{D}}

\newcommand{\vQ}{\bm{Q}}
\newcommand{\vQhat}{\widehat{\vQ}}
\newcommand{\vQbar}{\bar{\vQ}}
\newcommand{\vR}{\bm{R}}

\newcommand{\vS}{\bm{S}}
\newcommand{\vT}{\bm{T}}
\newcommand{\vU}{\bm{U}}

\newcommand{\vV}{\bm{V}}

\newcommand{\vW}{\bm{W}}

\newcommand{\vX}{\bm{X}}

\newcommand{\vY}{\bm{Y}}

\newcommand{\vZ}{\bm{Z}}
\newcommand{\vZero}{\bm{0}}

\newcommand{\Rbar}{\bar{R}}

\newcommand{\Tbar}{\bar{T}}

\newcommand{\RR}{\mathcal{R}}
\newcommand{\RRbar}{\bar{\RR}}
\renewcommand{\SS}{\mathcal{S}}

\newcommand{\TT}{\mathcal{T}}
\newcommand{\TTbar}{\bar{\TT}}

\newcommand{\VV}{\mathcal{V}}

\newcommand{\bQQ}{\bm{\mathcal{Q}}}
\newcommand{\bQQbar}{\bar{\bQQ}}

\newcommand{\bUU}{\bm{\mathcal{U}}}

\newcommand{\bVV}{\bm{\mathcal{V}}}

\newcommand{\bXX}{\bm{\mathcal{X}}}

\newcommand{\inv}{{-1}}
\newcommand{\tinv}{{-T}}

\newcommand{\bigO}[1]{\mathcal{O}\left(#1\right)}

\newcommand{\norm}[1]{\left\lVert#1\right\rVert}

\newcommand{\normF}[1]{\norm{#1}_{\text{F}}}

\makeatletter
\newsavebox{\@brx}
\newcommand{\llangle}[1][]{\savebox{\@brx}{\(\m@th{#1\langle}\)}%
	\mathopen{\copy\@brx\kern-0.5\wd\@brx\usebox{\@brx}}}
\newcommand{\rrangle}[1][]{\savebox{\@brx}{\(\m@th{#1\rangle}\)}%
	\mathclose{\copy\@brx\kern-0.5\wd\@brx\usebox{\@brx}}}
\makeatother

\newcommand{\MATLAB}{\textsc{Matlab}\xspace}

\newcommand{\rpltol}{\delta}
\newcommand{\eps}{\varepsilon}

\newcommand{\resTS}{\Delta_{TS}}
\newcommand{\resTTSS}{\Delta_{\TT\SS}}
\newcommand{\resQR}{\Delta_{\vQbar \Rbar}}
\newcommand{\resQQRR}{\Delta_{\bQQbar \RRbar}}
\newcommand{\resTR}{\Gamma_{TR}}
\newcommand{\resTTRR}{\Gamma_{\TT\RR}}

\newcommand{\randuniform}{\texttt{rand\_uniform}\xspace}
\newcommand{\randnormal}{\texttt{rand\_normal}\xspace}
\newcommand{\rankdef}{\texttt{rank\_def}\xspace}
\newcommand{\laeuchli}{\texttt{laeuchli}\xspace}
\newcommand{\monomial}{\texttt{monomial}\xspace}
\newcommand{\stewart}{\texttt{stewart}\xspace}
\newcommand{\stewartextreme}{\texttt{stewart\_extreme}\xspace}

\newcommand{\glued}{\texttt{glued}\xspace}
\newcommand{\sstep}{\texttt{s-step}\xspace}
\newcommand{\newton}{\texttt{newton}\xspace}

\newcommand{\BCGS}{\texttt{BCGS}\xspace}	
\newcommand{\BCGSP}{\texttt{BCGS-P}\xspace}	
\newcommand{\BCGSPIP}{\texttt{BCGS-PIP}\xspace}	
\newcommand{\BCGSPIO}{\texttt{BCGS-PIO}\xspace}	
\newcommand{\BCGSRO}{\texttt{BCGS+}\xspace}	
\newcommand{\BCGSIRO}{\texttt{BCGSI+}\xspace}	
\newcommand{\BCGSIROT}{\texttt{BCGSI+T}\xspace}	
\newcommand{\BCGSIROLS}{\texttt{BCGSI+LS}\xspace}	
\newcommand{\BCGSSROR}{\texttt{BCGSS+rpl}\xspace}	
\newcommand{\BMGS}{\texttt{BMGS}\xspace}	
\newcommand{\BMGST}{\texttt{BMGS-T}\xspace}	
\newcommand{\BMGSSVL}{\texttt{BMGS-SVL}\xspace}	
\newcommand{\BMGSLTS}{\texttt{BMGS-LTS}\xspace}	
\newcommand{\BMGSCWY}{\texttt{BMGS-CWY}\xspace}	
\newcommand{\BMGSICWY}{\texttt{BMGS-ICWY}\xspace}	
\newcommand{\DBMGS}{\texttt{DBMGS}\xspace}	

\newcommand{\IO}[1]{\texttt{IntraOrtho}\left(#1\right)}	
\newcommand{\IOnoarg}{\texttt{IntraOrtho}\xspace}	
\newcommand{\IOnoargs}{\texttt{IntraOrtho}s\xspace}	
\newcommand{\CGS}{\texttt{CGS}\xspace}	
\newcommand{\CGSP}{\texttt{CGS-P}\xspace}	
\newcommand{\CGSRO}{\texttt{CGS+}\xspace}	
\newcommand{\CGSIRO}{\texttt{CGSI+}\xspace}	
\newcommand{\CGSIROLS}{\texttt{CGSI+LS}\xspace}	
\newcommand{\CGSSRO}{\texttt{CGSS+}\xspace}	
\newcommand{\CGSSROR}{\texttt{CGSS+rpl}\xspace}	
\newcommand{\MGS}{\texttt{MGS}\xspace}	
\newcommand{\MGSRO}{\texttt{MGS+}\xspace}	
\newcommand{\MGSIRO}{\texttt{MGSI+}\xspace}	
\newcommand{\MGSSVL}{\texttt{MGS-SVL}\xspace}	
\newcommand{\MGSLTS}{\texttt{MGS-LTS}\xspace}	
\newcommand{\MGSICWY}{\texttt{MGS-ICWY}\xspace}	
\newcommand{\MGSCWY}{\texttt{MGS-CWY}\xspace}	
\newcommand{\HouseQR}{\texttt{HouseQR}\xspace}	
\newcommand{\GivensQR}{\texttt{GivensQR}\xspace}	
\newcommand{\TSQR}{\texttt{TSQR}\xspace}	
\newcommand{\CholQR}{\texttt{CholQR}\xspace}	
\newcommand{\mCholQR}{\texttt{mCholQR}\xspace}	
\newcommand{\CholQRRO}{\texttt{CholQR+}\xspace}	
\newcommand{\ShCholQRRORO}{\texttt{ShCholQR++}\xspace}	

\begin{document}
\begin{frontmatter}
\title{Block Gram-Schmidt algorithms and their stability properties}

\author[1]{Erin Carson\corref{cor1}%
\fnref{fn1,fn2}}
\ead{carson@karlin.mff.cuni.cz}

\author[1]{Kathryn Lund\fnref{fn1}}
\ead{kathryn.d.lund@gmail.com}

\author[2]{Miroslav Rozlo\v{z}n\'{i}k\fnref{fn3}}
\ead{miro@math.cas.cz}

\author[3]{Stephen Thomas\fnref{fn2,fn4}}
\ead{stephethomas@gmail.com}

\cortext[cor1]{Corresponding author}
\fntext[fn1]{Supported by Charles University PRIMUS project no.\ PRIMUS/19/SCI/11 and Charles University Research program no.\ UNCE/SCI/023.}
\fntext[fn2]{Supported by the Exascale Computing Project (17-SC-20-SC), a collaborative effort of the U.S. Department of Energy Office of Science and the National Nuclear Security Administration.}
\fntext[fn3]{Supported by the Academy of Sciences of the Czech Republic (RVO 67985840) and by the Grant Agency of the Czech Republic, Grant No.\ 20-01074S.}
\fntext[fn4]{National Renewable Energy Laboratory is operated by Alliance for Sustainable
Energy, LLC, for the U.S. Department of Energy (DOE) under Contract No.
DE-AC36-08GO28308.}

\address[1]{Faculty of Mathematics and Physics, Charles University, Prague, Czech Republic}
\address[2]{Institute of Mathematics, Czech Academy of Sciences, Prague, Czech Republic}
\address[3]{National Renewable Energy Laboratory, Boulder, Colorado, USA}

\begin{abstract}
	Block Gram-Schmidt algorithms serve as essential kernels in many scientific computing applications, but for many commonly used variants, a rigorous treatment of their stability properties remains open.  This work provides a comprehensive categorization of block Gram-Schmidt algorithms, particularly those used in Krylov subspace methods to build orthonormal bases one block vector at a time.  Known stability results are assembled, and new results are summarized or conjectured for important communication-reducing variants.  Additionally, new block versions of low-synchronization variants are derived, and their efficacy and stability are demonstrated for a wide range of challenging examples.  
	Numerical examples are computed with a versatile \MATLAB package hosted at \url{https://github.com/katlund/BlockStab}, and scripts for reproducing all results in the paper are provided.  Block Gram-Schmidt implementations in popular software packages are discussed, along with a number of open problems. An appendix containing all algorithms type-set in a uniform fashion is provided.
\end{abstract}


\end{frontmatter}


\section{Introduction} \label{sec:intro}
Given a matrix $\bXX \in \spR^{m \times n}$, $m \gg n$, we treat block Gram-Schmidt (BGS) algorithms that return an ``economic" QR factorization
\[
\bXX = \bQQ \RR,
\]
where $\bQQ \in \spR^{m \times n}$ has orthonormal columns which span the same space as the columns of $\bXX$ and $\RR \in \spR^{n \times n}$ is upper triangular with positive diagonal entries. 

There are many applications in which block Gram-Schmidt algorithms are more suitable than their single-vector counterparts. Block Krylov subspace methods (KSMs), based on underlying block Arnoldi/Lanczos procedures, are widely used for solving clustered eigenvalue problems (in which the block size is chosen based on the number of eigenvalues in a cluster) as well as linear systems with multiple right-hand sides. For seminal work in this area, see, e.g., Golub and Underwood~\cite{GolubUnderwood1977}, Stewart~\cite{Stewart2001}, and O'Leary~\cite{OLeary1980}. Block approaches can also have an advantage from a performance standpoint, as they replace BLAS-2 (i.e., vector-wise) operations with the more cache-friendly BLAS-3 (i.e., matrix-wise) operations; see, e.g., the work of Baker, Dennis, and Jessup~\cite{BakerDennisJessup2006}.

Block-partitioning strategies are also a key component in new algorithms designed to reduce communication in high-performance computing (HPC). The $s$-step (or, \emph{communication-avoiding}) Arnoldi/GMRES algorithms are based on a block orthogonalization strategy in which \emph{inter-block orthogonalization} is accomplished via a BGS method; for details see~\cite[Section 8]{BallardCarsonDemmel2014} and the references therein. There is also recent work by Grigori, Moufawad, and Nataf on ``enlarged'' KSMs, which are a special case of block KSMs~\cite{GrigoriMoufawadNataf2016}; here the block size is based on a partitioning of the problem and is designed to reduce communication cost.

Implicit in all formulations of BGS is the choice of an \emph{intra-block orthogonalization} scheme, to which we refer throughout as the \IOnoarg or ``muscle" of a BGS ``skeleton."  An \IOnoarg is any routine that takes a tall and skinny matrix (i.e., block vector) $\vX \in \spR^{m \times s}$ and returns an economic QR factorization $\vQ R = \vX$, with $\vQ \in \spR^{m \times s}$ and $R \in \spR^{s \times s}$
. Popular choices include classical Gram-Schmidt (with reorthogonalization), modified Gram-Schmidt, Householder-based QR (\HouseQR), Tall-Skinny QR (\TSQR)~\cite{DemmelGrigoriHoemmen2008}, and Cholesky-based QR (\CholQR).  Backward stability analyses for all these algorithms are well established and understood, and they are foundational for the analysis of BGS methods; see, e.g., \cite{Bjorck1967, BjorckPaige1992, Higham2002, PaigeRozloznikStrakos2006, YamamotoNakatsukasaYanagisawa2015, FukayaKannanNakatsukasa2020}, and especially the survey by Leon, Bj\"orck, and Gander \cite{LeonBjorckGander2013}.

On the other hand, stability analysis for BGS, especially in the form of block Arnoldi, is either non-existent or applies only to a small set of special cases. We focus on the \emph{loss of orthogonality} among the computed basis vectors as the primary measure of stability.  To be precise, we let $\bQQbar \in \spR^{m \times n}$ and $\RRbar \in \spR^{n \times n}$ denote the computed versions of $\bQQ \in \spR^{m \times n}$, $\RR \in \spR^{n \times n}$ for a given algorithm, where $\bQQ \RR = \bXX$ is the QR decomposition of $\bXX$ in exact arithmetic.  We define \emph{loss of orthogonality} as the quantity
\begin{equation} \label{eq:LOO}
	\norm{I_n - \bQQbar^T \bQQbar},
\end{equation}
for some norm $\norm{\cdot}$. We take $\norm{\cdot}$ to be the Euclidean norm, unless otherwise noted. This quantity often plays an important role in down-stream applications, such as the analysis of (block) GMRES/Arnoldi methods.

Letting $\eps$ denote the unit roundoff, we seek bounds on \eqref{eq:LOO} in terms of $\bigO{\eps}$ and $\kappa(\bXX)$, which we define as the 2-norm condition number of $\bXX$ or ratio between its largest and smallest singular values.  We say that an algorithm is \emph{unconditionally stable} if \eqref{eq:LOO} is bounded by $\bigO{\eps}$ for any matrix $\bXX$.  In contrast, we say an algorithm is \emph{conditionally stable} if the bound on \eqref{eq:LOO} is in terms of $\kappa(\bXX)$ or only holds with some constraint on $\kappa(\bXX)$. Note that this use differs from the notion of ``conditional (backward) stability'' often used in the literature. 

Another quantity that plays an important role in stability studies is the \emph{relative residual}:
\begin{equation} \label{eq:res}
	\frac{\norm{\bQQbar \RRbar - \bXX}}{\norm{\bXX}}.
\end{equation}
We take for granted that all algorithms considered here (with the exception of, perhaps, \CGSSROR, \BCGSSROR, \BMGSCWY and \BMGSICWY) produce a residual bounded by $\bigO{\eps}$.

A quantity that features strongly in the stability analysis of \CGS-based algorithms is what we call the \emph{relative Cholesky residual},
\begin{equation} \label{eq:res_chol}
	\frac{\norm{\bXX^T \bXX - \RRbar^T \RRbar}}{\norm{\bXX}^2},
\end{equation}
which measures how closely an algorithm approximates the Cholesky decomposition of $\bXX^T \bXX$.

Recent work by Barlow \cite{Barlow2019}, as well as important results by Barlow and Smoktunowicz \cite{BarlowSmoktunowicz2013}, Vanderstraeten \cite{Vanderstraeten2000}, and Jalby and Philippe \cite{JalbyPhilippe1991}, appear to be the only rigorous stability treatments of BGS. Some of the most influential papers on block Krylov subspace methods devote little or no attention to the stability analysis of the underlying block Arnoldi algorithm, and implicitly, BGS, routine.  Consider, for example, the following highly cited books and papers, most of which implicitly or explicitly use an algorithm derived from the block modified Gram-Schmidt skeleton, referred to here as \BMGS and given as Algorithm~\ref{alg:BMGS}:
\begin{itemize}
	\item 1980: O'Leary \cite{OLeary1980} proposes block Conjugate Gradients (BCG) for solving linear systems with multiple right-hand sides.  Although BCG is not necessarily derived from \BMGS, O'Leary makes a recommendation relevant to this conversation, namely, that either a ``QR algorithm or a modified Gram-Schmidt algorithm" be used as an \IOnoarg.
	\item 1993: Sadkane \cite{Sadkane1993} uses a block Arnoldi (Algorithm 1 in that text) that is based on \BMGS to compute the leading eigenpairs of large sparse unsymmetric matrices.  The ``QR method on CRAY2" is specified as the \IOnoarg, presumably a Householder-based QR (\HouseQR) routine.
	\item 1996: Simoncini and Gallopoulos \cite{SimonciniGallopoulos1996a} use a block Arnoldi algorithm as a part of block GMRES but they do not provide pseudocode or specify what is used as the \IOnoarg.
	\item 2003: Saad \cite{Saad2003} provides block Arnoldi routines based on both block Classical 	(Algorithm~6.22) and block Modified (Algorithm~6.23) Gram-Schmidt, with the \IOnoarg specified only as a ``QR factorization," again presumably \HouseQR.
	\item 2005: Morgan \cite{Morgan2005} uses Ruhe's variant of block Arnoldi (see, e.g., \cite[Algorithm~6.24]{Saad2003}) to implement block GMRES and block QMR and recommends full reorthogonalization (line (7) of Block-GMRES-DR) to increase stability, because good eigenvalue approximations are needed for deflation.  Ruhe's block Arnoldi variant is not attuned for BLAS-3 operations.
	\item 2006: Baker, Dennis, and Jessup \cite{BakerDennisJessup2006} propose a variant of block GMRES tuned to minimize memory movement. They use no reorthogonalization, and their block Arnoldi (lines 4-11 in Figure 1: B-LGMRES$(m,k)$) is based on \BMGS with ``QR factorization" specified as the \IOnoarg (probably a memory-sensitive version of \HouseQR).
	\item 2007: Gutknecht \cite{Gutknecht2007} discusses a number of important theoretical properties of block KSMs, particularly the issue of rank deficiency among columns of computed basis vectors. Block Arnoldi (Algorithm 9.1) is based on \BMGS and again, only a ``QR factorization" is specified as the \IOnoarg, presumably \HouseQR.
\end{itemize}
The intention of this list is merely to highlight a hidden danger: often, block Gram-Schmidt algorithms are used in applications without a full specification of the algorithm and without an accompanying stability analysis. To our knowledge, no one has formally proven stability bounds for \BMGS with \HouseQR, and yet the numerical linear algebra community (and beyond) has been using this algorithm ubiquitously for decades. Even worse, it \textit{is} known that \BMGS with, for example, \MGS, can suffer from steep orthogonality loss (see \cite{JalbyPhilippe1991}, our Section~\ref{sec:bgs_stab_1}, as well as Section~\ref{sec:BMGS_stab}). 

Our goal with the present work is to unify the literature on BGS algorithms in a common language to clarify what is known and what is unknown in terms of their numerical behavior. 
Our contributions are as follows:
\begin{itemize}
\item We provide a classification of existing variants of BGS algorithms within a unifying framework, based on a ``skeleton-muscle'' analogy first proposed by Hoemmen \cite{Hoemmen2010};
\item We assemble existing stability bounds and note commonly used skeleton-muscle combinations for which analysis remains open;
\item We derive a new block generalization of a low-synchronization algorithm of \'{S}wirydowicz, et al.\ \cite{SwirydowiczLangouAnanthan2021};
\item We give new insights into the mechanism by which stability is improved in block low-synchronization variants and also demonstrate the incompatibility between different classes of low-synchronization skeletons and muscles;
\item We show how the results of Jalby and Philippe can be adapted to give bounds on the loss of orthogonality for \BMGS combined with an unconditionally stable \IOnoarg (such as \HouseQR or \TSQR); and 
\item Using our versatile, publicly-available \MATLAB package developed to accompany this work, we perform extensive numerical experiments on various skeleton-muscle combinations for a variety of well-known test problems.
\end{itemize}


The rest of the paper proceeds as follows. In Section~\ref{sec:bgs_skel_musc} we discuss the history and communication properties of all the BGS skeletons we have identified.  
Some of them are new, in
particular, we derive a new block generalization of one of the low-synchronization algorithms by \'{S}wirydowicz, et al.\ \cite{SwirydowiczLangouAnanthan2021}.  Sections~\ref{sec:bgs_stab_1} and \ref{sec:bgs_stab_2} treat stability properties in more detail.  Section~\ref{sec:bgs_stab_1} contains heatmaps that allow for a quick comparison of many BGS variants all at once for a fixed matrix, while Section~\ref{sec:bgs_stab_2} focuses on one skeleton at a time and how its stability is affected across different \IOnoargs and condition numbers.  We give an overview of known mixed-precision BGS algorithms in Section~\ref{sec:mixed_prec}.  In Section~\ref{sec:software}, we examine which variants of BGS are implemented in well-known software.  We conclude the survey in Section~\ref{sec:conclusions} and identify open problems and future directions. \ref{sec:musc-alg} and \ref{sec:sror} contain pseudocode and \MATLAB scripts, respectively, for algorithms not included in the main text, and \ref{sec:scripts} contains \MATLAB command line calls for reproducing the plots in Sections~\ref{sec:bgs_stab_1} and \ref{sec:bgs_stab_2}.

\section{Block Gram-Schmidt variants and a skeleton-muscle analogy} \label{sec:bgs_skel_musc}
We assume that the matrix $\bXX \in \spR^{m \times n}$, $m \gg n$ is partitioned into a set of $p$ {\em block vectors}, each of size $m \times s$ ($n = ps$), i.e., 
\[
\bXX = [\vX_1 \,|\, \vX_2 \,|\, \cdots \,|\, \vX_p],
\]
and that a block Gram-Schmidt (BGS) algorithm returns an ``economic" QR factorization $\bXX = \bQQ \RR$, 
where $\bQQ = [\vQ_1 \,|\, \vQ_2 \,|\, \cdots \,|\, \vQ_p] \in \spR^{m \times n}$ has the same block-partitioned structure as $\bXX$, and $\RR \in \spR^{n \times n}$ can also be thought of as a $p \times p$ matrix with matrix-valued entries of size $s \times s$.\footnote{To help remember what each index represents throughout the text, note that usually $m > n > p > s$, which is also an alphabetical ordering.}

To simplify complex indexing, we first note that matrices like $\bQQ$ and $\RR$ have implicit tensor structure: $\bQQ$ can be viewed as a third-order tensor, being a vector of block vectors, and $\RR$ as a fourth-order tensor, being a matrix of matrices.  None of the algorithms we consider requires further partitioning explicitly, so when we index $\bQQ$ or $\RR$, we are always referring to the corresponding contiguous block component, be it an $m \times s$ block vector or $s \times s$ block entry, respectively.  This is best demonstrated by some examples:
\[
\bQQ_{1:j} = [\vQ_1 \,|\, \cdots \,|\, \vQ_j] \in \spR^{m \times sj}.
\]
where $1:j = \{1, 2, \dots, j\}$ is an indexing vector.  Similarly,
\[
\RR_{1:j,k} =
\begin{bmatrix} R_{1,k}	\\ R_{2,k}	\\ \vdots 	\\ R_{j,k}	\\ \end{bmatrix}
\in \spR^{sj \times s}.
\]
For all Cholesky-based algorithms, we let $R = \chol(A)$ denote an algorithm that takes a symmetric, numerically positive-definite matrix $A \in \spR^{s \times s}$ and returns an upper triangular matrix $R \in \spR^{s \times s}$.

We focus on scenarios where $\bXX$ is a set of block vectors that may not all be available at once, as is the case in block Arnoldi or block KSMs, wherein the set of block vectors is built by successively applying an operator $A$ to previously generated basis vectors.  As such, we only consider BGS variants that do not require access to all of $\bXX$ at once; in particular, we do not examine methods like CAQR from \cite{DemmelGrigoriHoemmen2012}.  We further require that BGS algorithms work left-to-right-- i.e., only information from the first $k$ block vectors of $\bQQ$ and $\bXX$ is necessary for generating $\vQ_{k+1}$-- and that they feature BLAS-3 operations. We also note that $p$ (the number of block vectors) may not be known a priori, as is the case in block KSMs, which will continue iterating until the specified convergence criterion is met.

Although we do not examine the performance of BGS methods in detail here, we do touch on their asymptotic communication properties, because they are crucial in a number of high-performance applications. In a simplified setting, the cost of an algorithm can be modeled in terms of the required computation, or number of floating point operations performed, and the required communication. By \emph{communication}, we mean the movement of data, both between levels of the memory hierarchy in sequential implementations and between parallel processors in parallel implementations. It is well established that communication and, in particular, synchronization between parallel processors, is the dominant cost (in terms of both time and energy) in large-scale settings; see, e.g., Bienz et al.~\cite{BienzGroppOlson2019}. It is therefore of interest to understand the potential trade-offs between the numerical properties of loss of orthogonality and stability in finite precision and the cost of communication in terms of number of messages and number of words moved.

There are a multitude of viable BGS algorithms which can all be succinctly described via the skeleton-muscle analogy from Hoemmen's thesis \cite{Hoemmen2010}.  We thoroughly wear out this metaphor and attempt to identify all viable skeleton and muscle options that have been considered before in the literature, as well as propose a few new ones.

\begin{table}[htbp!]
	\caption{Acronyms for algorithms.  For the block version of a Gram-Schmidt method, just add a ``B" prefix. \label{tab:acronyms}}
	\begin{tabular}{l|l}
		\CGS		& classical Gram-Schmidt	\\
		\CGSP		& \CGS, Pythagorean variant	\\
		\CGSRO		& \CGS, run twice	\\
		\CGSIRO		& \CGS with inner reorthogonalization	\\
		\CGSIROLS	& \CGSIRO, low-synchronization variant	\\
		\CGSSRO		& \CGS with selective reorthogonalization	\\
		\CGSSROR	& \CGS with selective reorthogonalization and replacement	\\
		\MGS		& modified Gram-Schmidt	\\
		\MGSRO		& \MGS run twice	\\
		\MGSIRO		& \MGS with inner reorthogonalization	\\
		\MGSSVL		& \MGS with Schreiber-Van-Loan reformulation	\\
		\MGSLTS		& \MGS with lower triangular solve	\\
		\MGSICWY	& \MGS with inverse compact WY	\\
		\MGSCWY		& \MGS with compact WY	\\
		\HouseQR	& QR via Householder reflections	\\
		\GivensQR	& QR via Givens rotations	\\
		\TSQR		& Tall-Skinny QR	\\
		\CholQR		& Cholesky QR	\\
		\mCholQR	& mixed-precision Cholesky QR	\\
		\CholQRRO	& \CholQR run twice	\\
		\ShCholQRRORO & shifted Cholesky QR with two stages of reorthogonalization
	\end{tabular}
\end{table}

The analogy is best understood via the literal meanings of skeleton and muscle in mammals.  Without muscles\footnote{...and tendons, ligaments, fasciae, bursae, etc.  See, e.g., \url{https://en.wikipedia.org/wiki/Human_musculoskeletal_system}}, we would be motionless bags of organs and bones.  Muscles provide the tension needed to stand, propel forward, grasp objects, and push things.  In much the same way, a BGS routine without an \IOnoarg specified is like a skeleton without muscles. However, just as we can still learn much about the human body by looking at bones, we can learn much from BGS skeletons without being distracted by the details of a specific \IOnoarg.  In the following sections, we examine the history, development, and communication properties of two classes of BGS: classical and modified.  Unless otherwise noted, single-column versions of block algorithms can be obtained by setting $s = 1$ and replacing $\IO{\cdot}$
with $\norm{\cdot}$.

Before proceeding, some comments about notation are warranted: throughout the text, we will describe a skeleton with a specific muscle choice as $\texttt{BGS} \circ \IOnoarg$, where $\texttt{BGS}$ is a general skeleton, $\IOnoarg$ is a general muscle, and $\circ$ stands for composition.  We also use our own naming system for algorithms, rather than what is proposed in their paper of origin, due to suffixes like ``2" being used inconsistently to denote reorthogonalization, BLAS2-featuring, or simply a second version of an
algorithm.  In general, we use the suffix $+$ to  denote a reorthogonalized variant.  A summary of acronyms used throughout the paper is provided in Table~\ref{tab:acronyms}.

\subsection{Block classical Gram-Schmidt skeletons} \label{sec:skel-BCGS} 
\subsubsection{Block Classical Gram-Schmidt (\BCGS)} \label{sec:BCGS_intro}
The \BCGS algorithm is a straightforward generalization of \CGS obtained by replacing vectors with block vectors and norms with \IOnoargs; see Algorithm~\ref{alg:BCGS}.  \BCGS is especially attractive in massively parallel settings because it passes $O(p)$ fewer messages asymptotically than block Modified Gram-Schmidt (\BMGS). See \cite[Table 2.4]{Hoemmen2010} and also discussions in \cite{BallardCarsonDemmel2014, DemmelGrigoriHoemmen2008, DemmelGrigoriHoemmen2012} for more details regarding the communication properties of \BCGS.

\begin{algorithm}[htbp!]
	\caption{$[\bQQ, \RR] = \BCGS(\bXX)$ \label{alg:BCGS}}
	\begin{algorithmic}[1]
		\STATE{Allocate memory for $\bQQ$ and $\RR$}
		\STATE{$[\vQ_1, R_{11}] = \IO{\vX_1}$}
		\FOR{$k = 1, \ldots, p-1$}
		\STATE{$\RR_{1:k,k+1} = \bQQ_{1:k}^T \vX_{k+1}$}
		\STATE{$\vW = \vX_{k+1} - \bQQ_{1:k} \RR_{1:k,k+1}$}
		\STATE{$[\vQ_{k+1}, R_{k+1,k+1}] = \IO{\vW}$}
		\ENDFOR
		\RETURN $\bQQ = [\vQ_1, \ldots, \vQ_p]$, $\RR = (R_{jk})$
	\end{algorithmic}
\end{algorithm}

By itself, \BCGS is not often considered a viable skeleton, especially not prior to its use in communication-avoiding Krylov subspace methods \cite{DemmelGrigoriHoemmen2008}.  This is likely due to the fact that \BCGS inherits many of the bad stability properties of \CGS.  In particular, like \CGS, the loss of orthogonality for \BCGS can be worse than quadratic in terms of $\kappa(\bXX)$; see Section~\ref{sec:BCGS_stab}.  To overcome this issue for \CGS, Smoktunowicz et al. \cite{SmoktunowiczBarlowLangou2006} introduced a variant of \CGS in which the diagonal
entries of the $R$-factor are computed via the Pythagorean theorem.  The resulting algorithm (\CGSP) has $\bigO{\eps} \kappa^2(\bXX)$, as long as $\bigO{\eps} \kappa^2(\bXX) < 1$, which can be a steep limitation for many practical applications.

A similar correction is possible for \BCGS, but a block generalization of the Pythagorean theorem, stated in the following theorem (see also \cite{CarsonLundRozloznik2021}), is needed instead.
\begin{theorem} \label{thm:block_pythag}
	Let full-rank block vectors $\vX, \vY, \vZ \in \spR^{n \times s}$ be such that $\vX = \vY + \vZ$ and $\vY \perp \vZ$ in the block sense, i.e., the spaces spanned by the columns of $\vY$ and $\vZ$ are orthogonal to each other. Then
	\begin{equation}
	\label{eq:pip}
		\vX^T \vX = \vY^T \vY + \vZ^T \vZ.
	\end{equation}
	In particular, if $\vX = \vQ R$, $\vY = \vU S$, and $\vZ = \vV T$ are economic QR decompositions (i.e., $\vQ, \vU, \vV \in \spR^{n \times s}$ are orthonormal, and $R, S, T \in \spR^{s \times s}$ are upper triangular), then
	\begin{equation}
	\label{eq:pio}
		R^T R = S^T S + T^T T.
	\end{equation}
\end{theorem}
We call the resulting algorithm BCGSP (i.e., \BCGS with a block Pythagorean correction). Using Theorem~\ref{thm:block_pythag}, we can derive two \BCGSP variants, one based on \eqref{eq:pip}, which we call \BCGSPIP (for ``Pythagorean inner product''), and one based on \eqref{eq:pio}, which we call \BCGSPIO (for ``Pythagorean intra-orthogonalization''). These variants are shown in Algorithms~\ref{alg:BCGSPIP} and \ref{alg:BCGSPIO}, respectively. 
Stability analysis for both \BCGSP variants can be found in \cite{CarsonLundRozloznik2021}; see also Section~\ref{sec:BCGS_stab}.  We note that the communication requirements of \BCGSP can be made asymptotically comparable to those of \BCGS, as long as the steps computing inner products or intra-orthogonalizations are coupled appropriately.  We note that an algorithm equivalent to \BCGSPIP was developed in \cite[Figure~2]{Yamazakietal2020}.

\begin{algorithm}[htbp!]
	\caption{$[\bQQ, \RR] = \BCGSPIP(\bXX)$ \label{alg:BCGSPIP}}
	\begin{algorithmic}[1]
		\STATE{Allocate memory for $\bQQ$ and $\RR$}
		\STATE{$[\vQ_1, R_{11}] = \IO{\vX_1}$}
		\FOR{$k = 1, \ldots, p-1$}
		\STATE{$\begin{bmatrix}	\RR_{1:k,k+1} \\ \mathcal{Z}_{k+1} \end{bmatrix} = [\bQQ_{1:k} \,\, \vX_{k+1}]^T \vX_{k+1}$}
		\STATE{$R_{k+1,k+1} = \chol(\mathcal{Z}_{k+1} - \RR_{1:k,k+1}^T \RR_{1:k,k+1})$}
		\STATE{$\vW = \vX_{k+1} - \bQQ_{1:k} \RR_{1:k,k+1}$}
		\STATE{$\vQ_{k+1} = \vW R_{k+1,k+1}^\inv$}
		\ENDFOR
		\RETURN $\bQQ = [\vQ_1, \ldots, \vQ_p]$, $\RR = (R_{jk})$
	\end{algorithmic}
\end{algorithm}

\begin{algorithm}[htbp!]
	\caption{$[\bQQ, \RR] = \BCGSPIO(\bXX)$ \label{alg:BCGSPIO}}
	\begin{algorithmic}[1]
		\STATE{Allocate memory for $\bQQ$ and $\RR$}
		\STATE{$[\vQ_1, R_{11}] = \IO{\vX_1}$}
		\FOR{$k = 1, \ldots, p-1$}
		\STATE{$\RR_{1:k,k+1} = \bQQ_{1:k}^T \vX_{k+1}$}
		\STATE{$\left[\sim, \begin{bmatrix} T_{k+1} & \\ & P_{k+1} \end{bmatrix} \right] = \IO{\begin{bmatrix} \vX_{k+1} & \\ & \RR_{1:k,k+1}\end{bmatrix}}$}
		\STATE{$R_{k+1,k+1} = \chol(T_{k+1}^T T_{k+1} - P_{k+1}^T P_{k+1})$}
		\STATE{$\vW = \vX_{k+1} - \bQQ_{1:k} \RR_{1:k,k+1}$}
		\STATE{$\vQ_{k+1} = \vW R_{k+1,k+1}^\inv$}
		\ENDFOR
		\RETURN $\bQQ = [\vQ_1, \ldots, \vQ_p]$, $\RR = (R_{jk})$
	\end{algorithmic}
\end{algorithm}

\subsubsection{\BCGS with reorthogonalization} \label{sec:BCGS_reortho_intro}
It is well known that reorthogonalization can stabilize \CGS; see, e.g., \cite{Abdelmalek1971}.  We note that there are essentially two reorthogonalization approaches: simply running the algorithm twice (\CGSRO) and performing each inner loop twice in a row (\CGSIRO, where the \texttt{I} stands for ``inner loop").  While the two approaches differ in floating-point arithmetic, they satisfy the same error bounds and have the same asymptotic communication cost.  In particular, they achieve $\bigO{\eps}$ loss of orthogonality, formulated as the ``twice is enough" principle in \cite{Parlett1998}, and double the communication cost of \CGS.  A variant of \CGSIRO that selectively chooses whether or not to run an inner loop twice is called \CGSSRO (the extra `S' for ``selective''), which can reduce the total number of floating-point operations, usually at the expense of increased communication due to the norm computations used to determine which vectors to reorthogonalize.  Many different selection criteria have been proposed; see, e.g., \cite{DanielGraggKaufman1976, Hoffmann1989, LeonBjorckGander2013, Stewart1998, Vanderstraeten2000}.

Reorthogonalization techniques can be easily generalized to block algorithms. \BCGSIRO (Algorithm~\ref{alg:BCGSIRO}) has received increasing attention in recent years and was analyzed in detail by Barlow and Smoktunowicz in 2013
\cite{BarlowSmoktunowicz2013}.  We explore some additional stability properties in Section~\ref{sec:BCGSIRO_stab}.

As an aside, note that running an \IOnoarg twice on each vector in \BCGS would generally not salvage orthogonality the way \BCGSRO and \BCGSIRO do.  The problem is that even with an incredibly stable muscle, \BCGS itself is unstable as a skeleton; see Section~\ref{sec:BCGS_stab} for more details, in particular,
the behavior of $\BCGS \circ \HouseQR$.

\begin{algorithm}[htbp!]
	\caption{$[\bQQ, \RR] = \BCGSIRO(\bXX)$ \label{alg:BCGSIRO}}	
	\begin{algorithmic}[1]
		\STATE{Allocate memory for $\bQQ$ and $\RR$}
		\STATE{$[\vQ_1, R_{11}] = \IO{\vX_1}$}
		\FOR{$k = 1, \ldots, p-1$}
		\STATE{$\RR_{1:k,k+1}^{(1)} = \bQQ_{1:k}^T \vX_{k+1}$ \% first BCGS step}
		\STATE{$\vW = \vX_{k+1} - \bQQ_{1:k} \RR_{1:k,k+1}^{(1)}$}
		\STATE{$[\vQhat, R_{k+1,k+1}^{(1)}] = \IO{\vW}$}
		\STATE{$\RR_{1:k,k+1}^{(2)} = \bQQ_{1:k}^T \vQhat$ \% second BCGS step}
		\STATE{$\vW = \vQhat - \bQQ_{1:k} \RR_{1:k,k+1}^{(2)}$}
		\STATE{$[\vQ_{k+1}, R_{k+1,k+1}^{(2)}] = \IO{\vW}$}
		\STATE{$\RR_{1:k,k+1} =  \RR_{1:k,k+1}^{(1)} + \RR_{1:k,k+1}^{(2)} R_{k+1,k+1}^{(1)}$ \% combine both steps}
		\STATE{$R_{k+1,k+1} = R_{k+1,k+1}^{(2)} R_{k+1,k+1}^{(1)}$}
		\ENDFOR
		\RETURN $\bQQ = [\vQ_1, \ldots, \vQ_p]$, $\RR = (R_{jk})$
	\end{algorithmic}
\end{algorithm}

\subsubsection{\BCGS with selective replacement and reorthogonalization} \label{sec:BCGSSROR_intro}
Reorthogonalization is not enough to recover stability for some (pathologically bad) matrices.  In response to such situations, Stewart \cite{Stewart2008} developed what we call \CGS with Selective Reorthogonalization and Replacement (\CGSSROR, Algorithm~\ref{alg:CGSSROR}), which replaces low-quality vectors with random ones of the same magnitude.  He also developed a block variant (\BCGSSROR, Algorithm~\ref{alg:BCGSSROR}), designed to work specifically with \CGSSROR as its muscle.

One of the primary reasons Stewart developed this variant was to account for \emph{orthogonalization faults}, which may occur when applying an \IOnoarg twice to the same block vector, possibly causing the norm of some columns to drop drastically due to extreme cancellation error and therefore lose orthogonality with respect to the columns of previously computed orthogonal block vectors.  Hoemmen describes this situation as ``naive reorthogonalization" in Section~2.4.7 of his thesis \cite{Hoemmen2010} and provides a theoretical example in Appendix~C.2 therein.

Both algorithms are highly technical and are perhaps best understood via the body of \MATLAB code accompanying this paper, especially the functions \texttt{cgs\_step\_sror} and \texttt{bcgs\_step\_sror}, which are reproduced in \ref{sec:sror}\footnote{These routines are precisely Stewart's original implementations, accessed in February 2020 at a cached version of \url{ftp://ftp.umiacs.umd.edu/pub/stewart/reports/Contents.html}. It is unclear if this page will remain accessible in the future.  Because the full algorithms are not printed in the paper \cite{Stewart2008}, we have therefore incorporated the original implementations into our package, to facilitate and preserve accessibility.}.  These routines are in fact the core of the algorithms and perform a careful reorthogonalization step that checks the quality of the reorthogonalization via a series of before-and-after norm comparisons and replaces reorthogonalized vectors with random ones if the quality is not satisfactory.  The subroutine \texttt{cgs\_step\_sror} additionally replaces identically zero vectors with random ones of small norm, an essential feature that is not explicitly mentioned in Stewart's paper \cite{Stewart2008} but is present in his \MATLAB implementations.  This feature is key in ensuring that \CGSSROR and \BCGSSROR can handle severely rank-deficient matrices, i.e., that they return a full-rank, orthogonal $\bQQ$ while passing rank deficiency onto $\RR$.

A major downside to \CGSSROR is that it can be as expensive as \HouseQR, due to the many norm computations and the potential for more than two orthogonalization steps per vector.  \BCGSSROR suffers the same fate and is generally ill-suited to high-performance contexts, because the frequent column-wise norm computations create communication bottlenecks.

It is also possible to formulate such an algorithm based on \MGS; see, e.g., \cite[Algorithm~13]{Hoemmen2010}).  We do not explore this further here, however, because there is little practical benefit from the communication point-of-view.

\begin{algorithm}[htbp!]
	\caption{$[\bQQ, \RR] = \BCGSSROR(\bXX, \rpltol)$ \label{alg:BCGSSROR}}
	\begin{algorithmic}[1]
		\STATE{Allocate memory for $\bQQ$ and $\RR$}
		\STATE{$[\vQ_1, \sim, R_{11}] = \texttt{bcgs\_step\_sror}(\vZero, \vX_1, \rpltol)$}
		\FOR{$k = 1, \ldots, p-1$}
		\STATE{$[\vQ_{k+1}, \RR_{1:k,k+1}, R_{k+1,k+1}] = \texttt{bcgs\_step\_sror}(\bQQ_{1:k}, \vX_{k+1}, \rpltol)$}
		\ENDFOR
		\RETURN{$\bQQ = [\vQ_1, \ldots, \vQ_s]$, $\RR = (R_{jk})$}
	\end{algorithmic}
\end{algorithm}

\subsubsection{One-sync \BCGSIRO (\BCGSIROLS)} \label{sec:BCGSIROLS_intro}
Recently, a number of low-synchronization variants of Gram-Schmidt have been proposed \cite{SwirydowiczLangouAnanthan2021}. Here, a \emph{synchronization point} refers to a global reduction which requires $\log P$ time to complete for $P$ processors or MPI ranks.  This occurs with operations such as inner products and norms because large vectors are often stored and operated on in a distributed fashion.

\CGSIRO (if implemented like Algorithm~\ref{alg:BCGSIRO} with $s=1$) has up to four synchronization points per column, while Algorithm~3 from \cite{SwirydowiczLangouAnanthan2021}-- which we denote here as \CGSIROLS, where ``LS" stands for ``low-synchronization"-- has just one per column. See Algorithm~\ref{alg:CGSIROLS}.  Note that it is formulated a bit differently from the presentation in \cite{SwirydowiczLangouAnanthan2021}; in particular, lines 21-23 introduce a single final synchronization point necessary for orthogonalizing the last column.

\begin{algorithm}[htbp!]
	\caption{$[\vQ, R] = \CGSIROLS(\vX)$ \label{alg:CGSIROLS}}
	\begin{algorithmic}[1]
		\STATE{Allocate memory for $\vQ$ and $R$}
		\STATE{$\vu = \vx_1$}
		\FOR{$k = 2, \ldots s$}
		\IF{$k = 2$}
		\STATE{$[r_{k-1,k-1}^2 \,\, \rho] = \vu^T [\vu \,\, \vx_k]$}
		\ELSIF{$k > 2$}
		\STATE{$\begin{bmatrix} \vw & \vz \\ \omega & \zeta \end{bmatrix}
			= [\vQ_{1:k-2} \,\, \vu]^T [\vu \,\, \vx_k ]$}
		\STATE{$[r_{k-1,k-1}^2 \,\, \rho] = [\omega \,\, \zeta] - \vw^T [\vw \,\, \vz]$}
		\ENDIF
		\STATE{$r_{k-1,k} = \rho / r_{k-1,k-1}$}
		\IF{$k = 2$}
		\STATE{$\vq_{k-1} = \vu / r_{k-1,k-1}$}
		\ELSIF{$k > 2$}
		\STATE{$R_{1:k-2,k-1} = R_{1:k-2,k-1} + \vw$}
		\STATE{$R_{1:k-2,k} = \vz$}
		\STATE{$\vq_{k-1} = (\vu - \vQ_{1:k-2} \vw) / r_{k-1,k-1}$}
		\ENDIF
		\STATE{$\vu = \vx_k - \vQ_{1:k-1} R_{1:k-1,k}$}
		\ENDFOR
		\STATE{$\begin{bmatrix} \vw \\ \omega \end{bmatrix} = [\vQ_{1:s-1} \,\, \vu]^T \vu$}
		\STATE{$r_{s,s}^2 = \omega - \vw^T \vw$}
		\STATE{$R_{1:s-1,s} = R_{1:s-1,s} + \vw$}
		\STATE{$\vq_s = (\vu - \vQ_{1:s-1} \vw) / r_{s,s}$}
		\RETURN{$\vQ = [\vq_1, \ldots, \vq_s]$, $R = (r_{jk})$}
	\end{algorithmic}
\end{algorithm}

Note that \CGSIROLS differs from Cholesky-based QR approaches such as \CholQR, which has a single synchronization point for the entire algorithm. Reorthogonalized \CholQR (\CholQRRO) \cite{YamamotoNakatsukasaYanagisawa2015} and shifted and reorthogonalized \CholQR (\ShCholQRRORO) \cite{FukayaKannanNakatsukasa2020} require two and three total synchronization points, respectively.  See Algorithms~\ref{alg:CholQR}-\ref{alg:ShCholQRRORO} to compare details.

The low-synchronization algorithms from \cite{SwirydowiczLangouAnanthan2021} are based on two ideas. The first is to compute a strictly lower triangular matrix $L$ (i.e., with zeros on the diagonal) one row or block of rows at a time in a single global reduction to account for all the inner products needed in the current iteration.  Each row $L_{k-1,1:k-2} = (\vQ_{1:k-2}^T \vq_{k-1})^T$ is obtained one at a time within the current step.  The second idea is to lag the normalization step and merge it into this single reduction.

Ruhe \cite{Ruhe83} observed that \MGS and \CGS could be interpreted as Gauss-Seidel and Gauss-Jacobi iterations, respectively, for solving the normal equations where the associated orthogonal projector is given as
\[
	I - \vQ_{1:k-1} T_{1:k-1,1:k-1} \vQ_{1:k-1}^T, \mbox{ for } T_{1:k-1,1:k-1} \approx (\vQ_{1:k-1}^T \vQ_{1:k-1})^{-1}
\]
For \CGSIROLS, the matrix $T$ would be given iteratively by the following (originally deduced in \cite{SwirydowiczLangouAnanthan2021}):
\[
	T_{1:k-1,1:k-1} = I - L_{1:k-1,1:k-1} - L_{1:k-1,1:k-1}^T.
\]
Another way to think of what \CGSIROLS does per step is the following: one can split the auxiliary matrix $T_{1:k-1,1:k-1}$ into two parts and apply them across two iterations. The lower triangular matrix $I - L_{1:k-1,1:k-1}$ is applied first, followed by a lagged correction
\[
	R_{1:k-2,k-1} = R_{1:k-2,k-1} + \vw, \quad \vw = -L_{k-1,1:k-2}^T,
\]
which is a delayed reorthogonalization step in the next iteration; see line~14 in Algorithm~\ref{alg:CGSIROLS}, but note that normalization there is delayed. Thus, the notion of reorthogonalization is modified and occurs ``on-the-fly," as opposed to requiring a complete second pass of the algorithm.  We emphasize, however, that $T$ and $L$ are not computed explicitly in Algorithm~\ref{alg:CGSIROLS}.

A block generalization of Algorithm~\ref{alg:CGSIROLS} is rather straightforward; see Algorithm~\ref{alg:BCGSIROLS}.  One only has to be careful with the diagonal entries of $\RR$, which are now $s \times s$ matrices.  In particular, line~5 from Algorithm~\ref{alg:CGSIROLS} must be replaced with a Cholesky factorization, and instead of dividing by $r_{k-1}$ in lines~10, 12, and 16, it is necessary to invert either $R_{k-1,k-1}$ or its transpose; compare with lines~10, 12, and 16 of Algorithm~\ref{alg:BCGSIROLS}.  Note that line~10 in particular is tricky: by combining previous quantities, it holds that
\[
P = \vU^T \vX_k \mbox{ or } \vU^T \big(I - \bQQ_{1:k-2} \bQQ_{1:k-2}^T \big) \vX_k,
\]
so we must apply $R_{k-1,k-1}^\tinv$ from the left in order to properly scale the $\vU$ ``hidden" in $P$.

Our block generalization in Algorithm~\ref{alg:BCGSIROLS} is essentially the same as \cite[Figure~3]{Yamazakietal2020}. There are no existing theoretical studies on the stability of \BCGSIROLS. However, it is clear that unlike \BCGSIRO, \BCGSIROLS is not guaranteed to exhibit $\bigO{\eps}$ loss of orthogonality even with an unconditionally stable \IOnoarg, likely due to the delayed operations; see Sections~\ref{sec:bgs_stab_1} and \ref{sec:BCGSIROLS_stab}.

\begin{algorithm}[htbp!]
	\caption{$[\bQQ, \RR] = \BCGSIROLS(\bXX)$ \label{alg:BCGSIROLS}}
	\begin{algorithmic}[1]
		\STATE{Allocate memory for $\bQQ$ and $\RR$}
		\STATE{$\vU = \vX_1$}
		\FOR{$k = 2, \ldots p$}
		\IF{$k = 2$}
		\STATE{$[R_{k-1,k-1}^T R_{k-1,k-1} \,\, P] = \vU^T [\vU \,\, \vX_k]$}
		\ELSIF{$k > 2$}
		\STATE{$\begin{bmatrix} \vW & \vZ \\ \Omega & Z \end{bmatrix}
			= [\bQQ_{1:k-2} \,\, \vU]^T [\vU \,\, \vX_k ]$}
		\STATE{$[R_{k-1,k-1}^T R_{k-1,k-1} \,\, P] = [\Omega \,\, Z] - \vW^T [\vW \,\, \vZ]$}
		\ENDIF
		\STATE{$R_{k-1,k} = R_{k-1,k-1}^\tinv P$}
		\IF{$k = 2$}
		\STATE{$\vQ_{k-1} = \vU R_{k-1,k-1}^\inv$}
		\ELSIF{$k > 2$}
		\STATE{$\RR_{1:k-2,k-1} = \RR_{1:k-2,k-1} + \vW$}
		\STATE{$\RR_{1:k-2,k} = \vZ$}
		\STATE{$\vQ_{k-1} = (\vU - \bQQ_{1:k-2} \vW) R_{k-1,k-1}^\inv$}
		\ENDIF
		\STATE{$\vU = \vX_k - \bQQ_{1:k-1} \RR_{1:k-1,k}$}
		\ENDFOR
		\STATE{$\begin{bmatrix} \vW \\ \Omega \end{bmatrix} = [\bQQ_{1:s-1} \,\, \vU]^T \vU$}
		\STATE{$R_{s,s}^T R_{s,s} = \Omega - \vW^T \vW$}
		\STATE{$\RR_{1:s-1,s} = \RR_{1:s-1,s} + \vW$}
		\STATE{$\vQ_s = (\vU - \bQQ_{1:s-1} \vW) R_{s,s}^\inv$}
		\RETURN{$\bQQ = [\vQ_1, \ldots, \vQ_s]$, $\RR = (R_{jk})$}
	\end{algorithmic}
\end{algorithm}

\subsection{Block modified Gram-Schmidt skeletons} \label{sec:skel-BMGS}
\subsubsection{Block Modified Gram-Schmidt (\BMGS)} \label{sec:BMGS_intro}
Much like \BCGS, \BMGS is obtained through a straightforward replacement of the column vectors of \MGS with block vectors; see Algorithm~\ref{alg:BMGS}.  An initial stability study can be found in \cite{JalbyPhilippe1991}, and the algorithm was perhaps first considered in a high-performance setting by Vital \cite{Vital1990}.

A quick comparison between Algorithms~\ref{alg:BCGS} and \ref{alg:BMGS} reveals that \BMGS has many more synchronization points than \BCGS and therefore suffers from higher communication demands.  Specifically, \BCGS has just one (block) inner product per block vector, whereas \BMGS splits this up into $k$ smaller (still block) inner products for the $(k+1)$st block vector.  At the same time, it is precisely this strategy that generally makes \BMGS more stable than \BCGS, which we discuss in more detail in Section~\ref{sec:BMGS_stab}.

\begin{algorithm}[htbp!]
	\caption{$[\bQQ, \RR] = \BMGS(\bXX)$ \label{alg:BMGS}}
	\begin{algorithmic}[1]
		\STATE{Allocate memory for $\bQQ$ and $\RR$}
		\STATE{$[\vQ_1, R_{11}] = \IO{\vX_1}$}
		\FOR{$k = 1, \ldots, p-1$}
		\STATE{$\vW = \vX_{k+1}$}
		\FOR{$j = 1, \ldots, k$}
		\STATE{$R_{j,k+1} = \vQ_j^T\vW$}
		\STATE{$\vW = \vW - \vQ_j R_{j,k+1}$}
		\ENDFOR
		\STATE{$[\vQ_{k+1}, R_{k+1,k+1}] = \IO{\vW}$}  \label{line:BMGS}
		\ENDFOR
		\RETURN{$\bQQ = [\vQ_1, \ldots, \vQ_p]$, $\RR = (R_{jk})$}
	\end{algorithmic}
\end{algorithm}

\subsubsection{Low-sync \BMGS variants} \label{sec:BMGS_low_sync_intro}
Although increases in floating-point operations or communication often correlate with better stability properties, several new low-synchronization variants of \MGS have been developed recently that appear to challenge this maxim.  We describe all four such variants here and show how they can be turned into block algorithms.  A full stability analysis for most of these algorithms remains open, but we outline a path forward in Section~\ref{sec:BMGS_low_sync_stab}.

\paragraph{Low-sync \MGS variants}
Recall that \CGS (Algorithm~\ref{alg:BCGS} with $s=1$) has two synchronization points per vector, namely the inner product in line 4 and the norm (\IOnoarg) in line 6, while \MGS (Algorithm~\ref{alg:BMGS} with $s=1$) has a linearly increasing number of $k+1$ synchronization points per vector, where $k$ increases linearly with the vector index.  The core steps of these algorithms can be expressed as the application of a projector or a sequence of projectors on the ``next" vector:
\begin{align}
	\CGS: & \quad (I - \vQ_{1:k} \vQ_{1:k}^T) \vx_{k+1} \label{eq:CGS_proj} \\
	\MGS: & \quad (I - \vq_k \vq_k^T) \cdots (I - \vq_1 \vq_1^T) \vx_{k+1} \label{eq:MGS_proj}
\end{align}
It is well known that the application of the single projector \eqref{eq:CGS_proj} is unstable, but breaking it up like in \eqref{eq:MGS_proj} improves stability.  A third option would be to compute \CGS's projector slightly differently, for example, by ``cushioning" it with a correction matrix $C_k \in \spR^{s \times s}$:
\[
(I - \vQ_{1:k} C_k \vQ_{1:k}^T) \vx_{k+1}.
\]
An algorithm built from such projectors would then communicate like \CGS.  The challenge, then, is to find $C_k$ so that orthogonality is lost like $\bigO{\eps} \kappa(\bXX)$ or better.  The matrix $T$ from Section~\ref{sec:BCGSIROLS_intro} would be one such candidate.

Recently, four such algorithms have been developed.  The first two (\MGSSVL and \MGSLTS) only seek to eliminate the increasing number of inner products and have two synchronization points per vector; the latter two (\MGSICWY and \MGSCWY) use a technique called ``normalization lagging," introduced by Kim and Chronopoulos (1992) \cite{Kim92}, in order to bring the number of synchronization points per vector down to one:
\begin{itemize}
	\item \MGSSVL \cite[Func.~3.1]{Barlow2019}: Barlow calls this variant ``MGS2." We prefer the suffix ``SVL" for ``Schreiber and Van Loan," on whose work \cite{SchreiberVanLoan1989} the development of this algorithm is largely based. See Algorithm~\ref{alg:MGSSVL} in \ref{sec:musc-alg}.
	\item \MGSLTS \cite[Alg.~4]{SwirydowiczLangouAnanthan2021}: This variant is not assigned a name in \cite{SwirydowiczLangouAnanthan2021}; it is designated only as ``Algorithm 4."  We refer to it with the suffix ``LTS" in this note, which stands for ``Lower Triangular Solve," due to how the correction matrix is handled.  See Algorithm~\ref{alg:MGSLTS} in \ref{sec:musc-alg}.
	\item \MGSCWY \cite[Alg.~6]{SwirydowiczLangouAnanthan2021}: The suffix here stands for ``Compact WY," and the algorithm turns out to be the one-sync version of \MGSSVL. See Algorithm~\ref{alg:MGSCWY} in \ref{sec:musc-alg}.
	\item \MGSICWY \cite[Alg.~5]{SwirydowiczLangouAnanthan2021}: The suffix here stands for ``Inverse Compact WY," referring to an alternative formulation of the \MGS projector \eqref{eq:MGS_proj}.  It turns out that \MGSICWY is the one-sync version of \MGSLTS.  See Algorithm~\ref{alg:MGSICWY} in \ref{sec:musc-alg}.
\end{itemize}
The differences between the four algorithms are highlighted in \diff{red}.

The forms of the projectors for these variants are quite similar, although the correction matrix itself is different for all algorithms, especially in floating-point arithmetic:
\begin{align}
	\MGSSVL, \MGSCWY:	& \quad (I - \vQ_{1:k} T_{1:k,1:k}^T \vQ_{1:k}^T) \vx_{k+1} \label{eq:MGSSVL_proj}\\
	\MGSLTS, \MGSICWY:	& \quad (I - \vQ_{1:k} T_{1:k,1:k}^\tinv \vQ_{1:k}^T) \vx_{k+1} \label{eq:MGSLTS_proj}
\end{align}
The projector for \MGSCWY has the same form as that of \MGSSVL, because both represent the product of rank-one elementary projectors as a recursively constructed triangular matrix $T$. The same is true for \MGSICWY and \MGSLTS; see, e.g., Walker \cite{Walker88} and the technical report by Sun \cite{Sun96}. To see how \MGSCWY and \MGSSVL share the same form of projector directly from the algorithms, start with \MGSCWY, and assume $k > 1$.  Then from line 14 of Algorithm~\ref{alg:MGSCWY}, we have that
\begin{align*}
	\vu
	& = \vw - \vQ_{:,1:k} R_{1:k,k+1} \\
	& = \vw - \vQ_{:,1:k} T_{1:k,1:k}^T \begin{bmatrix} \vr \\ \rho/ r_{k,k} \end{bmatrix} \mbox{, from line 12}\\
	& = \vw - \vQ_{:,1:k} T_{1:k,1:k}^T \begin{bmatrix} \vQ_{1:k-1}^T \vw \\ \vu^T\vw / r_{k,k} \end{bmatrix} \mbox{, from line 9} \\
	& = \vx_{k+1} - \vQ_{:,1:k} T_{1:k,1:k}^T \vQ_{1:k}^T \vx_{k+1} \mbox{, from lines 13 and 5}.
\end{align*}
By similar means, \MGSLTS and \MGSICWY can be related.

We also note that \MGSSVL and \MGSCWY use only matrix-matrix multiplication to apply their correction matrices, while \MGSLTS and \MGSICWY use lower triangular solves.  Depending on how these kernels are implemented, one variant may be preferred over another in practice. The $T^T$ matrix comes from the $(2,2)$ block of the transposed Householder matrix $U^T$ described by Barlow; see \cite{Barlow2019}, equation (2.17) on page 1262.

\paragraph{Block generalizations}
Block generalizations of both \MGSSVL (Algorithm~\ref{alg:BMGSSVL}) and \MGSLTS (Algorithm~\ref{alg:BMGSLTS}) are quite straightforward.  In \cite{Barlow2019}, \BMGSSVL appears as \texttt{MGS3} and \texttt{BMGS\_H} with \IOnoargs \MGSSVL (Algorithm~\ref{alg:MGSSVL}) and \HouseQR (i.e., Householder-based QR, as in, e.g., \cite{GolubVanLoan2013}), respectively.  \BMGSLTS as a BGS variant is new, as far as we are aware.

One quirk about these low-sync variants is that, unlike \BCGS or \BMGS, they require an \IOnoarg that produces a $T$-factor, in addition to $\vQ$ and $R$. For algorithms that do not explicitly produce such a factor, one must assume that $T = I$.  This correction matrix is necessary to ensuring the stability of the overall algorithm.  Furthermore, there is a kind of compatibility requirement between $T$-producing \IOnoargs and low-sync skeletons: combinations like $\BMGS \circ \MGSSVL$ or $\BMGSSVL \circ \MGSLTS$, for example, may not produce the expected stability behavior.  We explore this issue further in Section~\ref{sec:bgs_stab_1}.

It is also possible to derive block generalizations of \MGSCWY and \MGSICWY that are still one-sync and mostly stable; see \cite[Figure~3]{Yamazakietal2020}.  Instead of the implied square root in line~7 of Algorithms~\ref{alg:MGSCWY} and \ref{alg:MGSICWY}, a Cholesky factorization is needed to recover the block diagonal entry $R_{k,k}$.  The computation of $\vQ_k$ is then completed by inverting the $s \times s$ upper triangular matrix $R_{k,k}$.  As with \BCGSIROLS (cf.\ Section~\ref{sec:BCGSIROLS_intro}), one has to take care about when to apply $R_{k,k}^\inv$ or $R_{k,k}^\tinv$: line~12 of both Algorithms~\ref{alg:BMGSCWY} and \ref{alg:BMGSICWY} require applying $R_{k,k}^\tinv$ to properly scale the ``hidden" $\vU$ in $P$.  Note that both the Cholesky factorization and matrix inverse can be computed locally, because $s$ is small.  The final orthogonalization in lines~16-17 of Algorithms~\ref{alg:MGSCWY} and \ref{alg:MGSICWY} can be replaced by an \IOnoarg or another Cholesky factorization; we opt for an \IOnoarg in line~16 of Algorithms~\ref{alg:BMGSCWY} and \ref{alg:BMGSICWY}.

For further discussion of the stability properties of these low-sync block methods, see Section~\ref{sec:BMGS_low_sync_stab}.

\begin{algorithm}[htbp!]
	\caption{$[\bQQ, \RR, \TT] = \BMGSSVL(\bXX)$ \label{alg:BMGSSVL}}
	\begin{algorithmic}[1]
		\STATE{Allocate memory for $\bQQ$, $\RR$, and $\TT$}
		\STATE{$[\vQ_1, R_{11}, T_{11}] = \IO{\vX_1}$}
		\FOR{$k = 1, \ldots, p-1$}
		\STATE{$\RR_{1:k,k+1} = \diff{\TT_{1:k,1:k}^T} \big( \bQQ_{1:k}^T \vX_{k+1}\big )$}
		\STATE{$\vW = \vX_{k+1} - \bQQ_{1:k} \RR_{1:k,k+1}$}
		\STATE{$[\vQ_{k+1}, R_{k+1,k+1}, T_{k+1,k+1}] = \IO{\vW}$}
		\STATE{\diff{$\TT_{1:k,k+1} = -\TT_{1:k,1:k} \big( \bQQ_{1:k}^T\vQ_{k+1} \big) T_{k+1,k+1}$}}
		\ENDFOR
		\RETURN{$\bQQ = [\vQ_1, \ldots, \vQ_p]$, $\RR = (R_{jk})$, $\TT = (T_{jk})$}
	\end{algorithmic}
\end{algorithm}

\begin{algorithm}[htbp!]
	\caption{$[\bQQ, \RR, \TT] = \BMGSLTS(\vX)$ \label{alg:BMGSLTS}}
	\begin{algorithmic}[1]
		\STATE{Allocate memory for $\bQQ$, $\RR$, and $\TT$}
		\STATE{$[\vQ_1, R_{11}, T_{11}] = \IO{\vX_1}$}
		\FOR{$k = 1, \ldots, p-1$}
		\STATE{$\RR_{1:k,k+1} = \diff{\TT_{1:k,1:k}^\tinv} \big( \bQQ_{1:k}^T \vX_{k+1} \big)$}
		\STATE{$\vW = \vX_{k+1} - \bQQ_{1:k} \RR_{1:k,k+1}$}
		\STATE{$[\vQ_{k+1}, R_{k+1,k+1}, T_{k+1,k+1}] = \IO{\vW}$}
		\STATE{\diff{$\TT_{1:k,k+1} = \big( \bQQ_{1:k}^T \vQ_{k+1} \big) T_{k+1,k+1}$}}
		\ENDFOR
		\RETURN{$\bQQ = [\vQ_1, \ldots, \vQ_p]$, $\RR = (R_{jk})$, $\TT = (T_{jk})$}
	\end{algorithmic}
\end{algorithm}

\begin{algorithm}[htbp!]
	\caption{$[\bQQ, \RR, \TT] = \BMGSCWY(\bXX)$ \label{alg:BMGSCWY}}
	\begin{algorithmic}[1]
		\STATE{Allocate memory for $\bQQ$ and $\RR$}
		\STATE{$\TT = I_n$}
		\STATE{$\vU = \vX_1$}
		\FOR{$k = 1, \ldots, p-1$}
		\STATE{$\vW = \vX_{k+1}$}
		\IF{$k = 1$}
		\STATE{$\begin{bmatrix} R_{k,k}^T R_{k,k} & P \end{bmatrix} = \vU^T [\vU \,\, \vW]$ \% recover $R_{k,k}$ with \chol}
		\ELSIF{$k > 1$}
		\STATE{$\begin{bmatrix} \vT & \vR \\ R_{k,k}^T R_{k,k} & P \end{bmatrix} = [\bQQ_{1:k-1} \,\, \vU]^T [\vU \,\, \vW]$ \% recover $R_{k,k}$ with \chol}
		\STATE{$\TT_{1:k-1,k} = \diff{-\TT_{1:k-1,1:k-1} (\vT R_{k,k}^\inv)}$}
		\ENDIF
		\STATE{$\RR_{1:k,k+1} = \diff{\TT_{1:k,1:k}^T} \begin{bmatrix} \vR \\ R_{k,k}^\tinv P \end{bmatrix}$}
		\STATE{$\vQ_k = \vU R_{k,k}^\inv$}
		\STATE{$\vU = \vW - \bQQ_{:,1:k} \RR_{1:k,k+1}$}
		\ENDFOR
		\STATE{$[\vQ_p, R_{p,p}] = \IO{\vU}$}
		\RETURN{$\bQQ = [\vQ_1, \ldots, \vQ_p]$, $\RR = (R_{jk})$, $\TT = (T_{jk})$}
	\end{algorithmic}
\end{algorithm}

\begin{algorithm}[htbp!]
	\caption{$[\bQQ, \RR, \TT] = \BMGSICWY(\bXX)$ \label{alg:BMGSICWY}}
	\begin{algorithmic}[1]
		\STATE{Allocate memory for $\bQQ$ and $\RR$}
		\STATE{$\TT = I_n$}
		\STATE{$\vU = \vX_1$}
		\FOR{$k = 1, \ldots, p-1$}
		\STATE{$\vW = \vX_{k+1}$}
		\IF{$k = 1$}
		\STATE{$\begin{bmatrix} R_{k,k}^T R_{k,k} & P \end{bmatrix} = \vU^T [\vU \,\, \vW]$ \% recover $R_{k,k}$ with \chol} 
		\ELSIF{$k > 1$}
		\STATE{$\begin{bmatrix} \vT & \vR \\ R_{k,k}^T R_{k,k} & P \end{bmatrix} = [\bQQ_{1:k-1} \,\, \vU]^T [\vU \,\, \vW]$ \% recover $R_{k,k}$ with \chol}
		\STATE{$\TT_{1:k-1,k} = \diff{\vT R_{k,k}^\inv}$}
		\ENDIF
		\STATE{$\RR_{1:k,k+1} = \diff{\TT_{1:k,1:k}^\tinv} \begin{bmatrix} \vR \\ R_{k,k}^\tinv P \end{bmatrix}$}
		\STATE{$\vQ_k = \vU R_{k,k}^\inv$}
		\STATE{$\vU = \vW - \bQQ_{:,1:k} \RR_{1:k,k+1}$}
		\ENDFOR
		\STATE{$[\vQ_p, R_{p,p}] = \IO{\vU}$}
		\RETURN{$\bQQ = [\vQ_1, \ldots, \vQ_p]$, $\RR = (R_{jk})$, $\TT = (T_{jk})$}
	\end{algorithmic}
\end{algorithm}
	
\subsubsection{Dynamic \BMGS (\DBMGS)} \label{sec:DBMGS_intro}
Dynamic \BMGS is essentially \BMGS with variable block sizes determined adaptively in order to reduce the condition numbers of block vectors to be orthogonalized.  Vanderstraeten first proposed \DBMGS in 2000, particularly with \MGS in mind as the \IOnoarg \cite{Vanderstraeten2000}.  The communication costs associated with this approach can be expected to be somewhere between that of \MGS (i.e., with a block size of 1) and that of \BMGS using the specified maximum block size, depending on the numerical properties of the input data.  We will not devote
further discussion to \DBMGS here, because we assume that the block partitioning of $\bXX$ is fixed a priori.  However, a stability analysis of \DBMGS would greatly inform that of adaptive \s-step algorithms \cite{Carson2020}.

\section{Skeleton stability in terms of muscle stability: overview} \label{sec:bgs_stab_1}

Viewing BGS algorithms through the skeleton-muscle framework begs the following questions:
\begin{enumerate}[(Q.1)]
	\item If we use an unconditionally stable muscle, what is the best a skeleton can do? \label{q1}
	\item What are the minimum requirements on the muscle such that a given skeleton is stable? \label{q2}
\end{enumerate}

To answer either question, it will be helpful to keep in mind the stability properties of different muscles, summarized in Table~\ref{tab:musc_upper_bounds}.  Additionally, we summarize known and conjectured results for block algorithms in Table~\ref{tab:skel_upper_bounds}.  For the specific assumptions on the muscles that lead to these stability bounds, see the exposition in Section~\ref{sec:bgs_stab_2}.

\begin{table}[htbp!]
	\caption{Upper bounds on the loss of orthogonality in the $\vQbar$ factor for various \IOnoargs, along with proof references.  Note that all conditions have hidden constants in terms of polynomials of the dimensions $m$ and $s$.  A superscript dagger $^\dagger$ indicates that the result is conjectured based on numerical observations, but not yet proven. \label{tab:musc_upper_bounds}}
	\begin{center}
		\begin{tabular}{l|c|c|c}
		\IOnoarg	& $\norm{I - \vQbar^T \vQbar}_2$	& Assumption on $\kappa(\vX)$	& Reference(s)	\\ \hline
		\CGS		& $\bigO{\eps} \kappa^{n-1}(\vX)$ & $\bigO{\eps} \kappa(\vX) <1$		& \cite{Kielbasinski1971} \\
		\CGSP		& $\bigO{\eps} \kappa^2(\vX)$	& $\bigO{\eps} \kappa^2(\vX) <1$	& \cite{SmoktunowiczBarlowLangou2006} \\
		\CholQR		& $\bigO{\eps} \kappa^2(\vX)$	& $\bigO{\eps} \kappa^2(\vX) <1$	& \cite{YamamotoNakatsukasaYanagisawa2015}	\\
		\MGS		& $\bigO{\eps} \kappa(\vX)$		& $\bigO{\eps} \kappa(\vX) <1$		& \cite{Bjorck1967}	\\
		\MGSSVL		& $\bigO{\eps} \kappa(\vX)$		& $\bigO{\eps} \kappa(\vX) <1$		& \cite{Barlow2019}	\\
		\MGSLTS		& $\bigO{\eps} \kappa(\vX)^\dagger$		& $\bigO{\eps} \kappa(\vX) <1^\dagger$			& conjecture \cite{SwirydowiczLangouAnanthan2021}	\\
		\MGSCWY		& $\bigO{\eps} \kappa(\vX)^\dagger$		& $\bigO{\eps} \kappa(\vX) <1^\dagger$			& conjecture \cite{SwirydowiczLangouAnanthan2021}	\\
		\MGSICWY	& $\bigO{\eps} \kappa(\vX)^\dagger$		& $\bigO{\eps} \kappa(\vX) <1^\dagger$			& conjecture \cite{SwirydowiczLangouAnanthan2021}	\\
		\CholQRRO	& $\bigO{\eps}$					& $\bigO{\eps} \kappa^2(\vX) <1$ 	& \cite{YamamotoNakatsukasaYanagisawa2015}	\\
		\CGSIRO		& $\bigO{\eps}$					& $\bigO{\eps} \kappa(\vX) <1$		& \cite{Abdelmalek1971, BarlowSmoktunowicz2013, GiraudLangouRozloznikEshof2005}	\\
		\CGSSRO		& $\bigO{\eps}$					& $\bigO{\eps} \kappa(\vX) <1$		& \cite{DanielGraggKaufman1976, Hoffmann1989}	\\
		\CGSIROLS	& $\bigO{\eps}^\dagger$			& $\bigO{\eps} \kappa(\vX) <1^\dagger$					& conjecture \cite{SwirydowiczLangouAnanthan2021}	\\
		\MGSRO		& $\bigO{\eps}$					& $\bigO{\eps} \kappa(\vX) <1$		& \cite{JalbyPhilippe1991, GiraudLangou2002}	\\
		\MGSIRO		& $\bigO{\eps}$					& $\bigO{\eps} \kappa(\vX) <1$		& \cite{Hoffmann1989, Gander1980}	\\
		\ShCholQRRORO& $\bigO{\eps}$				& $\bigO{\eps} \kappa(\vX) <1$		& \cite{FukayaKannanNakatsukasa2020}	\\
		\CGSSROR	& $\bigO{\eps}^\dagger$			& none$^\dagger$					& conjecture \cite{Stewart2008}\\
		\HouseQR	& $\bigO{\eps}$					& none								& \cite[Sec.~19.3]{Higham2002}	\\
		\GivensQR	& $\bigO{\eps}$					& none								& \cite[Sec.~19.6]{Higham2002}	\\
		\TSQR		& $\bigO{\eps}$					& none								& \cite{DemmelGrigoriHoemmen2012,MoriYamamotoZhang2012}	\\
		\end{tabular}
	\end{center}
\end{table}

\begin{table}[htbp!]
	\caption{Upper bounds on the loss of orthogonality in the $\bQQbar$ factor for BGS skeletons \textbf{composed with unconditionally stable \IOnoargs}, along with proof references.  Note that all conditions have hidden constants in terms of the dimensional parameters $m$, $n$, and $s$.  A superscript dagger $^\dagger$ indicates that the result is conjectured but lacks a rigorous proof.  \label{tab:skel_upper_bounds}}
\begin{center}
	\begin{tabular}{l|c|c|c}
		\texttt{BGS}& $\norm{I - \bQQbar^T \bQQbar}_2$ 	& Assumption on $\kappa(\bXX)$				& Reference(s)	\\ \hline
		\BCGS			& $\bigO{\eps} \kappa^{n-1}(\bXX)^\dagger$	& $\bigO{\eps} \kappa(\bXX) < 1^\dagger$	& conjecture \\
		\BCGSP			& $\bigO{\eps} \kappa^2(\bXX)$		& $\bigO{\eps} \kappa^2(\bXX) < 1$			& \cite{CarsonLundRozloznik2021} \\
		\BMGS			& $\bigO{\eps}\kappa(\bXX)$					& $\bigO{\eps} \kappa(\bXX) < 1$			& \cite{JalbyPhilippe1991}, here \\
		\BMGSSVL		& $\bigO{\eps}\kappa(\bXX)$					& $\bigO{\eps} \kappa(\bXX) < 1$			& \cite{Barlow2019} \\
		\BMGSLTS		& $\bigO{\eps}\kappa(\bXX)^\dagger$			& $\bigO{\eps} \kappa(\bXX) < 1^\dagger$	& conjecture, here\\
		\BMGSCWY		& $\bigO{\eps}\kappa^2(\bXX)^\dagger$			& $\bigO{\eps} \kappa^2(\bXX) < 1^\dagger$	& conjecture, here \\
		\BMGSICWY		& $\bigO{\eps}\kappa^2(\bXX)^\dagger$			& $\bigO{\eps} \kappa^2(\bXX) < 1^\dagger$	& conjecture, here \\
		\BCGSIRO		& $\bigO{\eps}$								& $\bigO{\eps} \kappa(\bXX) < 1$			& \cite{BarlowSmoktunowicz2013}	\\
		\BCGSSROR		& $\bigO{\eps}^\dagger$						& none$^\dagger$							& conjecture \cite{Stewart2008}\\
		\BCGSIROLS  &  $\bigO{\eps}\kappa^2(\bXX)^\dagger$	& $\bigO{\eps} \kappa^2(\bXX) < 1^\dagger$	& conjecture, here
	\end{tabular}
\end{center}
\end{table}

To demonstrate stability properties for multiple skeleton-muscle combinations at once, we build heat-maps for a number of test matrices $\bXX \in \spR^{m \times ps}$, where $m = 10000$ (number of rows), $p = 50$ (number of block vectors), and $s = 10$ (columns per block vector).  Additional matrix properties are given in Table~\ref{tab:matrix_props}.  All heat-map plots are run in \MATLAB 2019b on the r3d3 login node of the Sn{\v e}hurka cluster, which consists of an Intel Xeon Processor E5-2620 with 15MB cache \@ 2.00 GHz and operating system 18.04.4 LTS.\footnote{\url{http://cluster.karlin.mff.cuni.cz/}.} Our code is publicly available at \url{https://github.com/katlund/BlockStab}. Our test problems include:

\begin{itemize}
	\item \randuniform: The entries of $\bXX$ are drawn randomly from the uniform distribution via the \texttt{rand} command in \MATLAB.
	\item \randnormal: The entries of $\bXX$ are drawn randomly from the normal distribution via the \texttt{randn} command in \MATLAB.
	\item \rankdef: $\bXX$ is first generated like \randnormal.  Then the first block vector is set to 100 times the last block vector, i.e., $\vX_1 = 100 \vX_{p}$, to ensure that the matrix is numerically rank-deficient and badly scaled.	
	\item \laeuchli: $\bXX$ is a L\"auchli matrix of the form
	\[
	\bXX = \begin{bmatrix}
	1 		& 1 	& \cdots	& 1 	\\
	\eta	& 		& 		 	&	  	\\
	& \eta	&		 	&		\\
	&		& \ddots 	&		\\
	&		&			& \eta
	\end{bmatrix},
	\quad \eta \in (\eps, \sqrt{\eps}),
	\]
	where $\eta$ is drawn randomly from a scaled uniform distribution.  This matrix is interesting, because columns are only barely linearly independent.
	\item \monomial: A diagonal $m \times m$ operator $A$ with evenly distributed eigenvalues in $(\frac{1}{10},10)$ is defined, and $p$ vectors $\vv_k$, $k = 1, \ldots, p$, are randomly generated from the uniform distribution and normalized.  The matrix $\bXX$ is then defined as the concatenation of $p$ block vectors
	\[
	\vX_k = [\vv_k \,|\, A \vv_k \,|\, \cdots \,|\, A^{s-1} \vv_k].
	\]
	\item \sstep: $\bXX$ is built similarly to \monomial, but now $\vv_k$ is the normalized version of $A^{s-1} \vv_{k-1}$.
	\item \newton: $\bXX$ is built like \sstep but with Newton polynomials instead of monomials.
	\item \stewart: $\bXX = \bUU \Sigma \bVV^T$, where $\bUU$ and $\bVV$ are random unitary matrices and the diagonal of $\Sigma$ is the geometric sequence from $1$ to $10^{-20}$.  The first and 25th columns of $\bXX$ are the same, and the 35th column is identically zero.
	\item \stewartextreme: $\bXX = \bUU \Sigma \bVV^T$, where $\bUU$ and $\bVV$ are random unitary matrices and the first half of the diagonal of $\Sigma$ is the geometric sequence from $1$ to $10^{-10}$, while the second half is identically zero.
\end{itemize}

\begin{table}[htbp!]
	\caption{Matrix properties. Note that singular values are computed via \texttt{svd} after $\bXX$ is explicitly formed.  When exact values are known, which is the case with the smallest singular value of \rankdef, and all the singular values of \stewart and \stewartextreme, those values are listed instead. \label{tab:matrix_props}}
	\begin{center}
		\begin{tabular}{c|c|c|c}
			Matrix ID				& $\sigma_{1}$	& $\sigma_{n}$	& $\kappa(\bXX)$ \\ \hline
			\texttt{\randuniform}	& 1.12e+03		& 2.25e+01		& 4.96e+01		\\ \hline 
			\texttt{\randnormal}	& 1.22e+02		& 7.77e+01		& 1.46e+00		\\ \hline 
			\texttt{\rankdef}		& 1.02e+04		& 0				& \texttt{Inf}	\\ \hline 
			\texttt{\laeuchli}		& 2.24e+01		& 2.02e-11		& 1.11e+12		\\ \hline 
			\texttt{\monomial}		& 2.32e+08		& 3.04e-04		& 7.63e+11		\\ \hline
		    \texttt{\sstep}			& 2.08e+01 		& 3.21e-17 		& 6.50e+17		\\ \hline 
			\texttt{\newton}		& 3.15e+00 		& 7.77e-04 		& 4.06e+03		\\ \hline			\texttt{\stewart}		& 1				& 1e-20			& 1e+20			\\ \hline 
			\texttt{\stewartextreme}& 1				& 0				& \texttt{Inf}
		\end{tabular}
	\end{center}
\end{table}

Note that while \stewart and \stewartextreme are taken from \cite{Stewart2008}, matrices like these have been used to study stability properties in Gram-Schmidt algorithms perhaps as early as Hoffman \cite{Hoffmann1989}.  However, as we shall see, only Stewart's algorithm \BCGSSROR and a couple variations of \BCGSIRO are consistently robust enough to reliably factor these matrices.  We therefore believe the names \stewart and \stewartextreme are well earned.

For each matrix, two heat-maps are produced -- one for loss of orthogonality and one for relative residual -- allowing for a quick and straightforward comparison of all methods at once.  The colorbars indicate the value range.  Note that \BCGSSROR is only compatible with \CGSSRO and \CGSSROR, so a \texttt{NaN} is returned for all other combinations. Note also that \BCGSSROR and \CGSSROR are denoted as `BCGSS+R' and `CGSS+R', respectively, in the axis labels.  To demonstrate the effect of replacement, \CGSSRO has a replacement tolerance $\rpltol = 0$, while \CGSSROR has $\rpltol = 100$, which Stewart \cite{Stewart2008} notes is relatively ``aggressive."  \BMGSSVL and \BMGSLTS have similar behavior, so we only include the former; likewise for \BMGSCWY and \BMGSICWY.  Recall that \BCGSIROLS does not require an \IOnoarg, and \BMGSCWY only at the last step, so the reported values are identical for their columns.  A \texttt{NaN} value is also returned when algorithms using Cholesky (\CholQR, \CholQRRO, \ShCholQRRORO, \BCGSIROLS, and \BMGSCWY) encounter a matrix that is not numerically positive definite, because \MATLAB's built-in \chol throws an error.

\subsection{\randuniform and \randnormal} \label{sec:rand}
\begin{figure}[htbp!]
	\begin{center}
		\caption{Measurements for \randuniform \label{fig:rand_uniform}}
		\begin{tabular}{cc}
			\includegraphics[width=.45\textwidth]{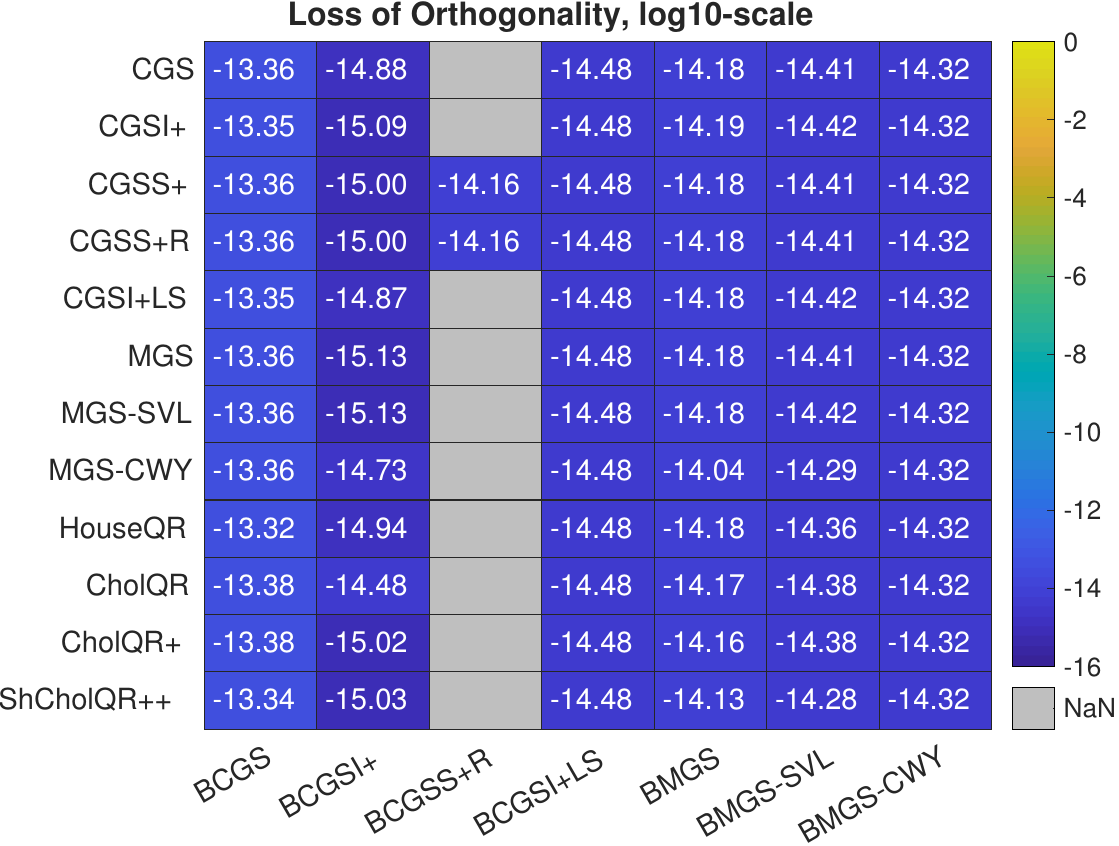} &
			\includegraphics[width=.45\textwidth]{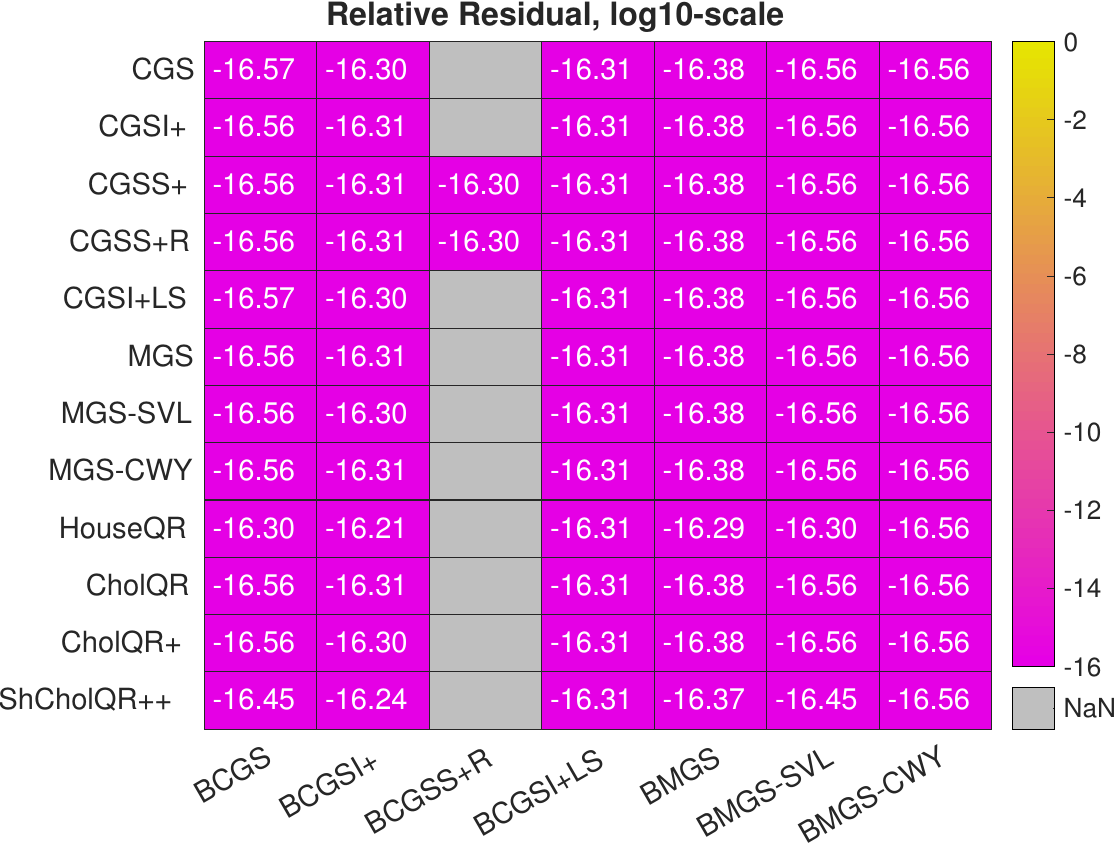}
		\end{tabular}
	\end{center}
\end{figure}

\begin{figure}[htbp!]
	\begin{center}
		\caption{Measurements for \randnormal \label{fig:rand_normal}}
		\begin{tabular}{cc}
			\includegraphics[width=.45\textwidth]{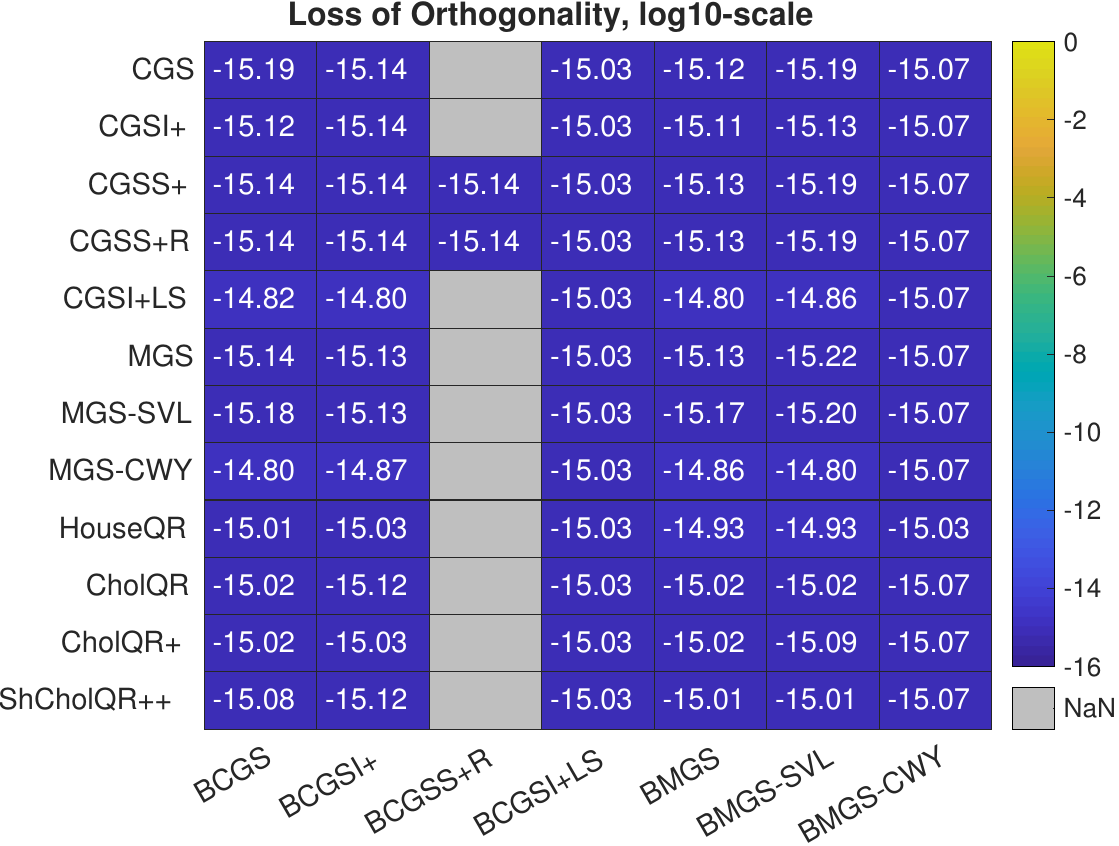} &
			\includegraphics[width=.45\textwidth]{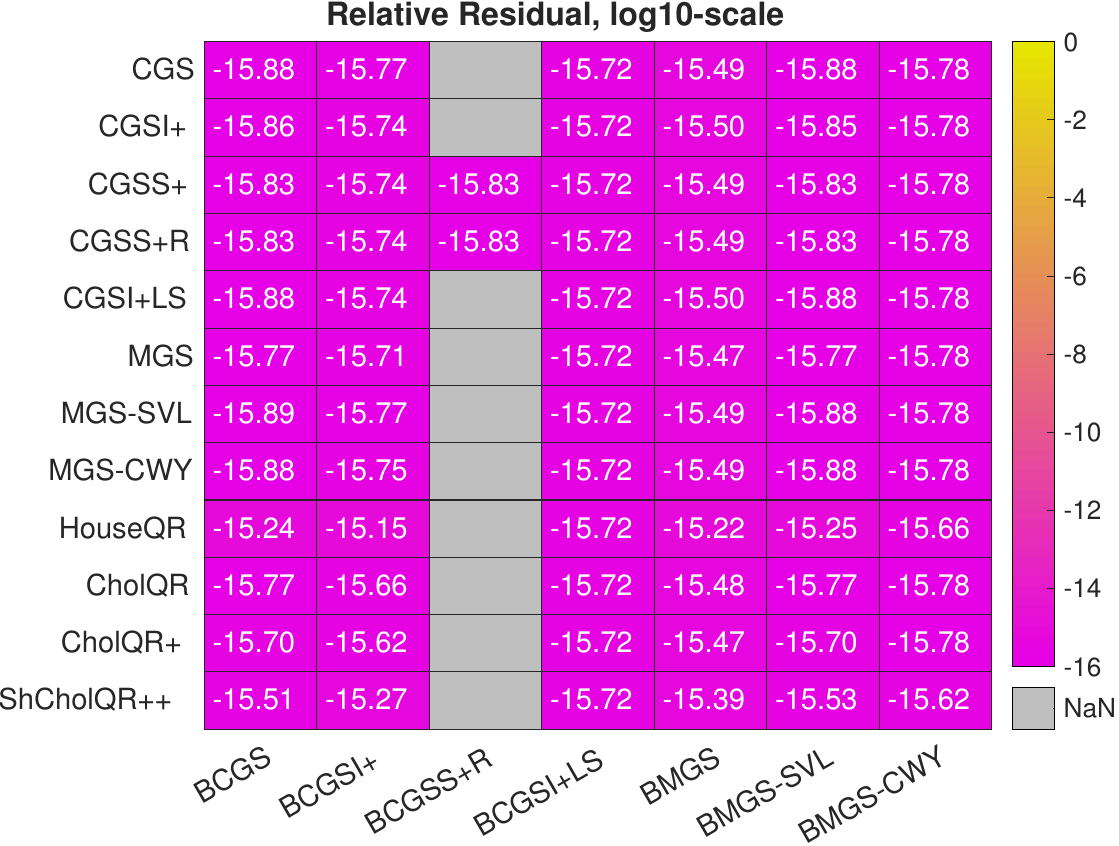}
		\end{tabular}
	\end{center}
\end{figure}

These tests serve primarily as sanity checks.  In comparing Figures~\ref{fig:rand_uniform} and \ref{fig:rand_normal}, there is little difference overall among the algorithms for either matrix.  The only noticeable differences arise for \randuniform (Figure~\ref{fig:rand_uniform}) where we see that \BCGS algorithms lose the most orthogonality, while \BCGSIRO the least; indeed, there is an order of magnitude difference between the two columns.  Other small differences can be seen for the muscles \CGSIROLS and \MGSCWY in Figure~\ref{fig:rand_normal}; for all methods they perform slightly worse.

\subsection{\rankdef} \label{sec:rank_def}
\begin{figure}[htbp!]
	\begin{center}
		\caption{Measurements for \rankdef \label{fig:rank_def}}
		\begin{tabular}{cc}
			\includegraphics[width=.45\textwidth]{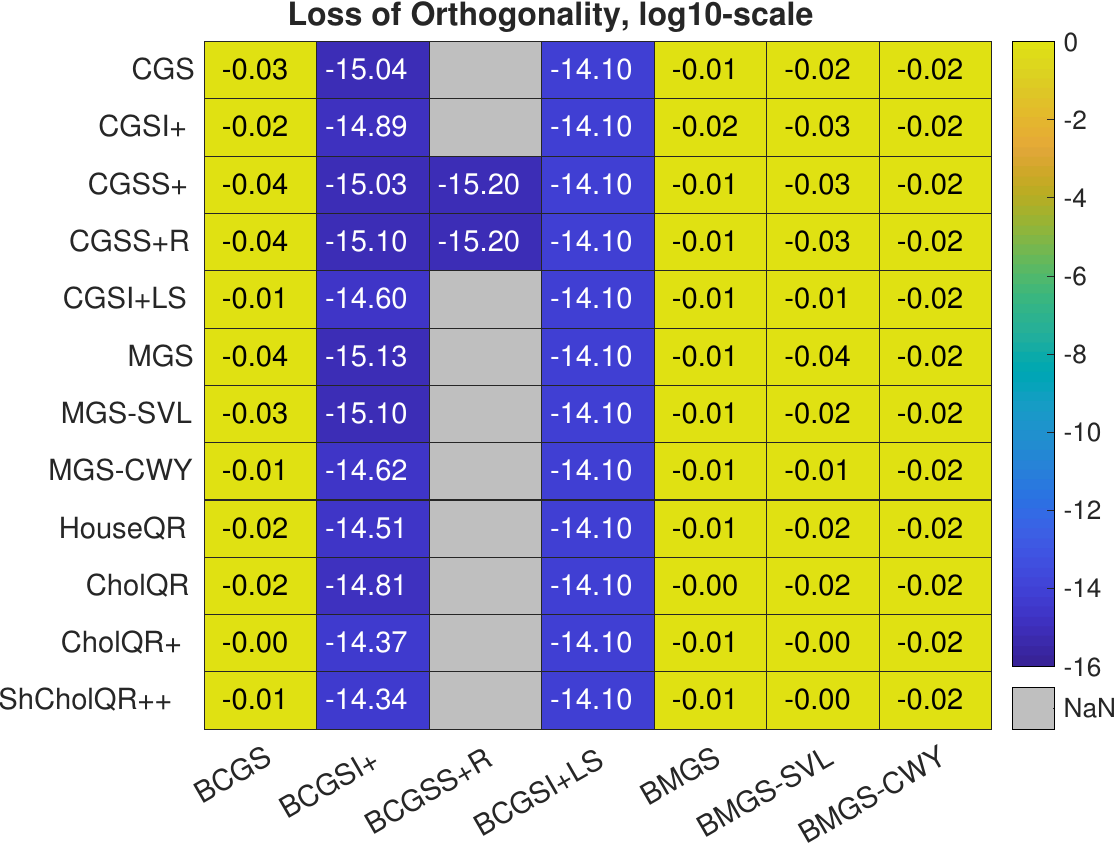} &
			\includegraphics[width=.45\textwidth]{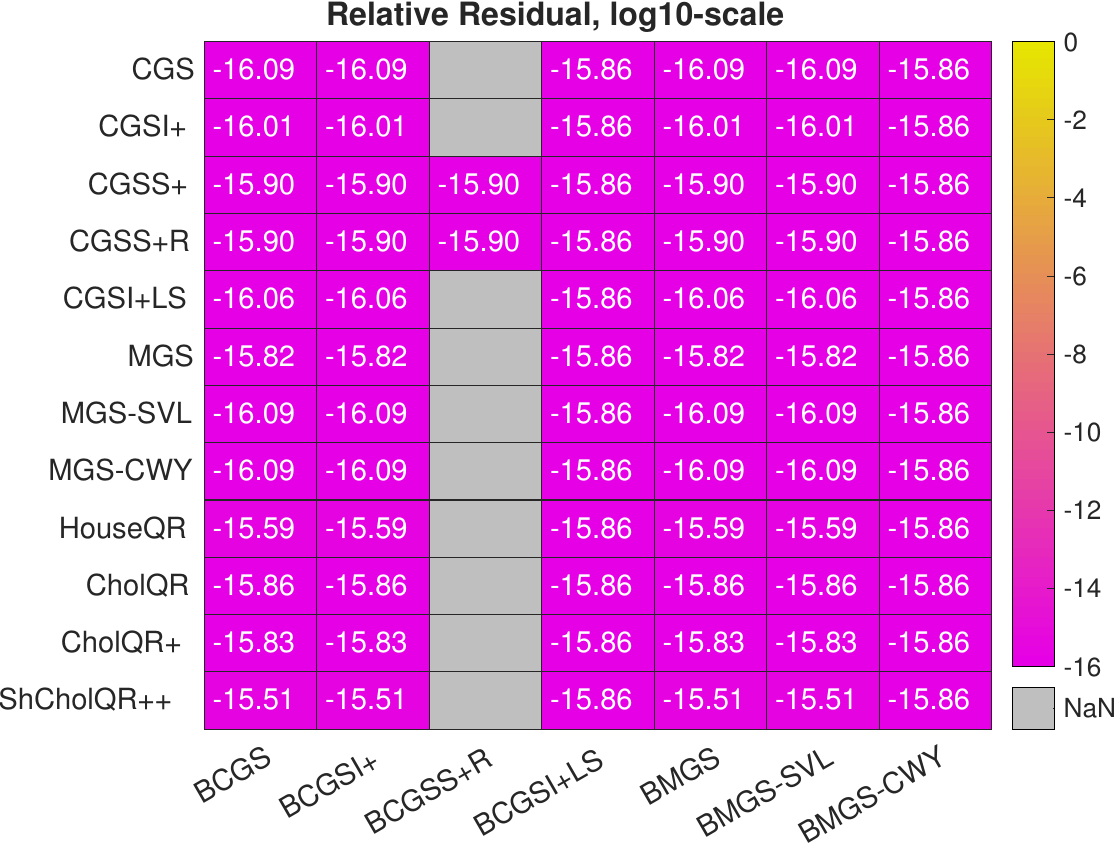}
		\end{tabular}
	\end{center}
\end{figure}

In Figure~\ref{fig:rank_def}, we see that only the skeletons with reorthogonalization are robust enough to handle the matrix.  It is further interesting to note that \BCGSIRO performs nearly the same regardless of \IOnoarg. 

\subsection{\laeuchli} \label{sec:laeuchli}
\begin{figure}[htbp!]
	\begin{center}
		\caption{Measurements for \laeuchli \label{fig:laeuchli}}
		\begin{tabular}{cc}
			\includegraphics[width=.45\textwidth]{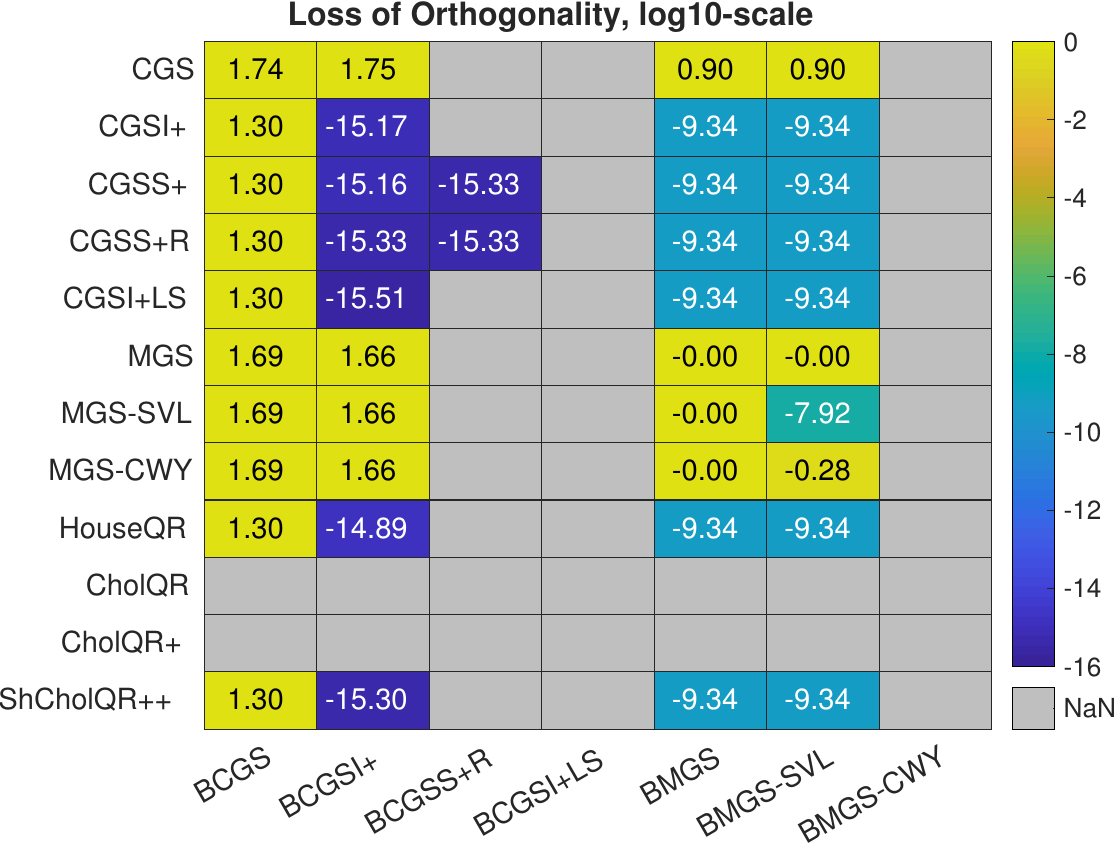} &
			\includegraphics[width=.45\textwidth]{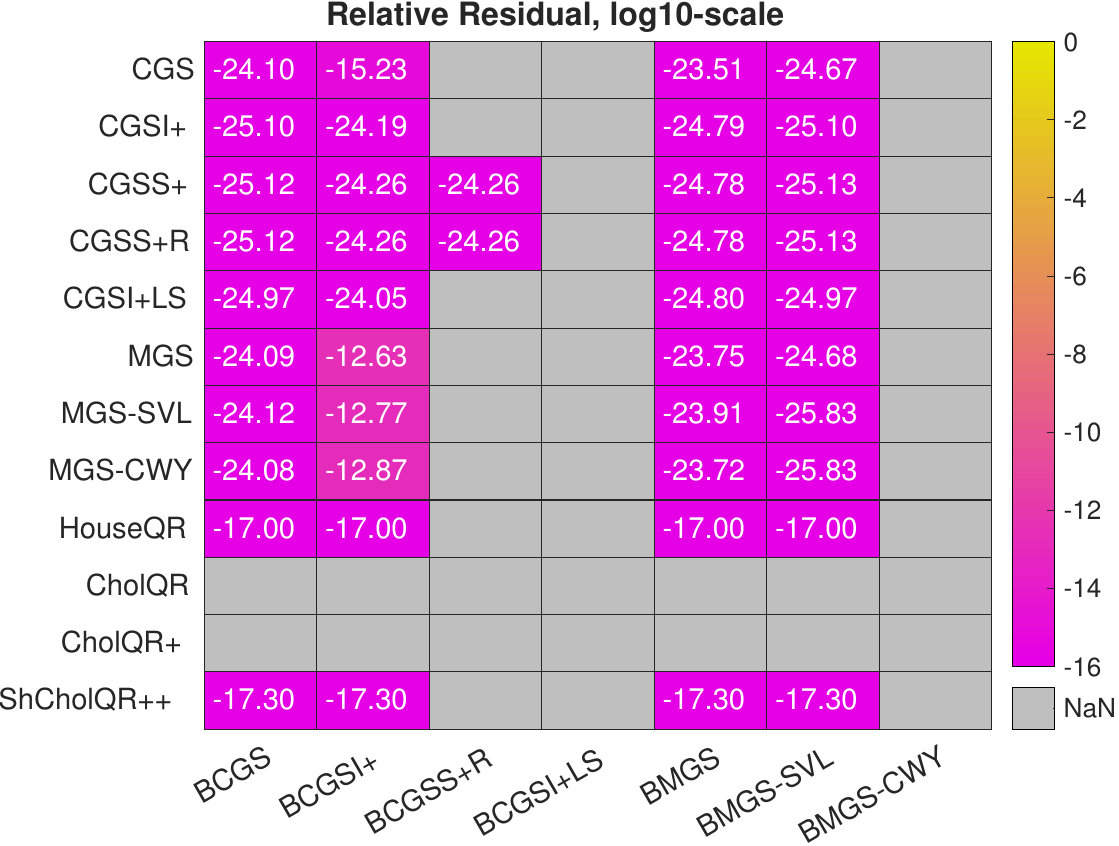}
		\end{tabular}
		\\
		\includegraphics[width=.25\textwidth]{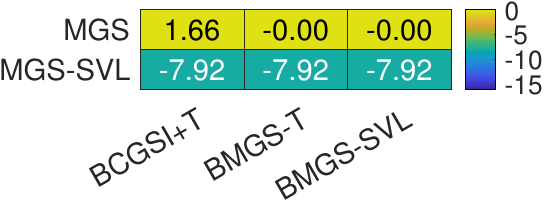}
	\end{center}
\end{figure}

The L\"auchli matrix reveals many interesting nuances. We first note in Figure~\ref{fig:laeuchli} that for all algorithms that run to completion, except $\BCGSIRO$ with $\CGS$, $\MGS$, $\MGSSVL$, and $\MGSCWY$, the relative residual exceeds double precision.  This is due to the nature of the L\"auchli matrix itself: because entries are either zero or close to $\eps$, there is a high rate of cancellation.  Incidentally, this same cancellation is the reason that neither \CholQR nor \CholQRRO are viable \IOnoargs, and \BCGSIROLS and \BMGSCWY are not viable skeletons: $\vX_1^T \vX_1$ is within $\bigO{\eps}$ of a matrix of all ones, which is not strictly positive definite and therefore violates the requirements in \MATLAB's \chol.

The other interesting point relates to the correction matrix $T$ of \MGSSVL.  Both \BCGSIRO and \BMGS suffer a total loss of orthogonality with \MGSSVL, but \BMGSSVL does not, since the \BMGSSVL skeleton was designed to incorporate this $T$ matrix from the \MGSSVL muscle. There is an easy fix for \BCGSIRO and \BMGS: we can incorporate the $T$ output of \MGSSVL into both algorithms by replacing every action of $\vQ_k$ with $\vQ_k T_{kk}$.  The second loss of orthogonality table in Figure~\ref{fig:laeuchli} demonstrates the outcome, with a ``T" suffix denoting the altered algorithms.  With the T-fix, both $\BCGSIROT \circ \MGSSVL$ and $\BMGST \circ \MGSSVL$ perform as well as $\BMGSSVL \circ \MGSSVL$; in fact, $\BCGSIROT \circ \MGSSVL$ matches the others in residual too.  (Note that \BCGS with the T-fix would reproduce \BMGSSVL.) We note also that, even though \MGSSVL and \MGSCWY are based on the same projector, \BMGSSVL does not work with the muscle \MGSCWY; orthogonality is completely lost.  

\subsection{\monomial} \label{sec:monomial}
\begin{figure}[htbp!]
	\begin{center}
		\caption{Measurements for \monomial \label{fig:monomial}}
		\begin{tabular}{cc}
			\includegraphics[width=.45\textwidth]{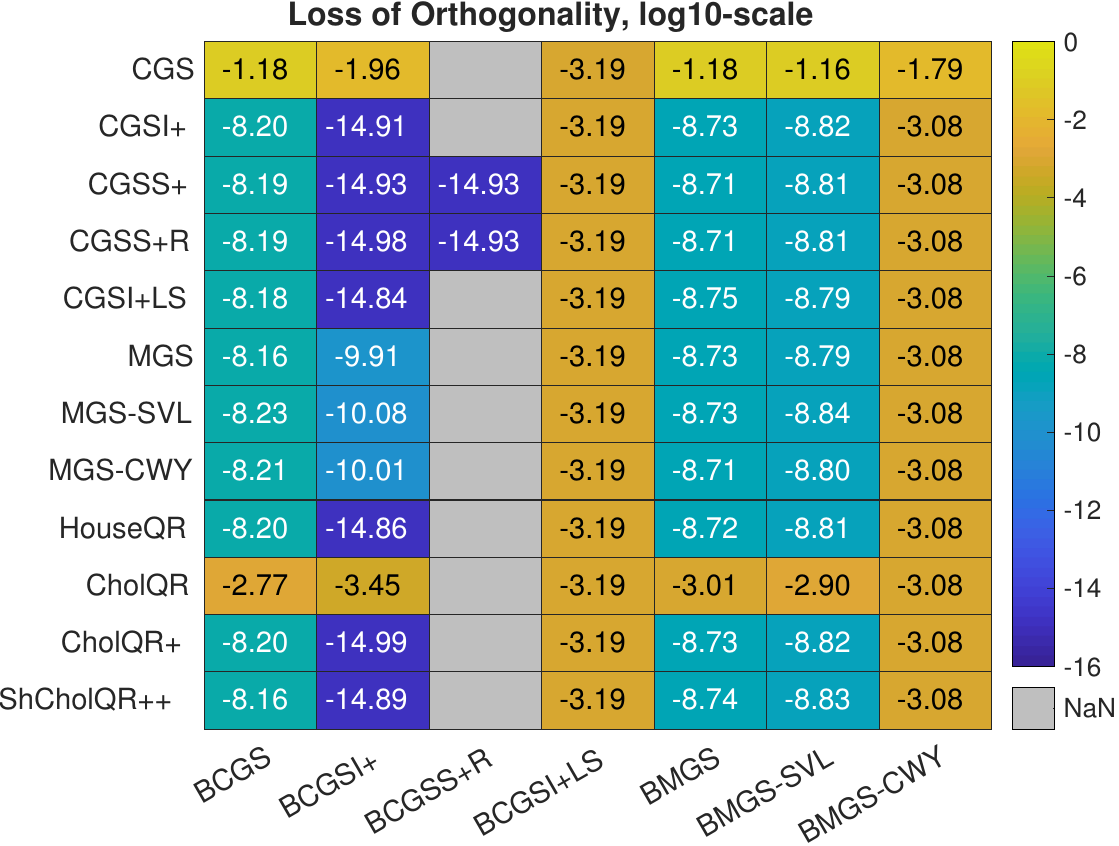} &
			\includegraphics[width=.45\textwidth]{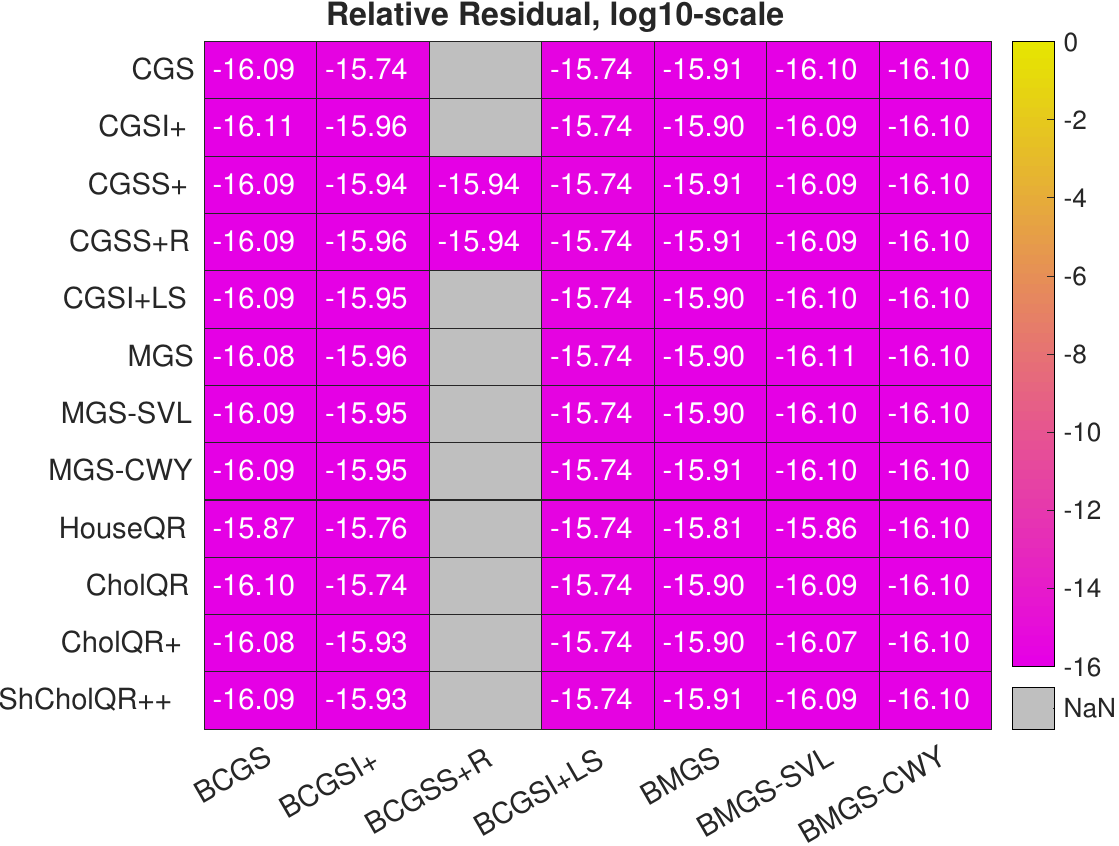}
		\end{tabular}
	\end{center}
\end{figure}

The \monomial matrix also reveals a few interesting properties.  In Figure~\ref{fig:monomial} we now see that all Cholesky-based methods can run, but they all lose quite a bit of orthogonality.  Also, \CGS in particular struggles across the board.  We again see that $\BCGSIRO \circ \MGSSVL$ (and also $\BCGSIRO \circ \MGSCWY$) does not perform as well as other variants of \BCGSIRO; however, the T-fix will not provide any benefit here, because $\BCGSIRO \circ \MGSSVL$ is already as stable as $\BMGSSVL \circ \MGSSVL$.


\subsection{\sstep and \newton} \label{sec:s-step}
The \sstep and \newton matrices simulate behavior one would encounter when implementing s-step Krylov subspace methods; see, e.g., \cite{BallardCarsonDemmel2014, Carson2015, Hoemmen2010}.  The \sstep matrix is poorly conditioned, as expected, because a monomial is used to build the basis; for these particular parameters and test matrix it is even numerically singular.  The \newton matrix simulates what would happen when Newton polynomials are used instead, which leads to much better conditioning.  It is interesting that only \BCGSSROR is capable of factoring the \sstep matrix; not even $\BCGSIRO \circ \HouseQR$ is capable.  One the other hand, both one-sync methods (\BCGSIROLS and \BMGSCWY) factor \newton satisfactorily, almost to full precision.


\begin{figure}[htbp!]
	\begin{center}
		\caption{Measurements for \sstep \label{fig:sstep}}
		\begin{tabular}{cc}
			\includegraphics[width=.45\textwidth]{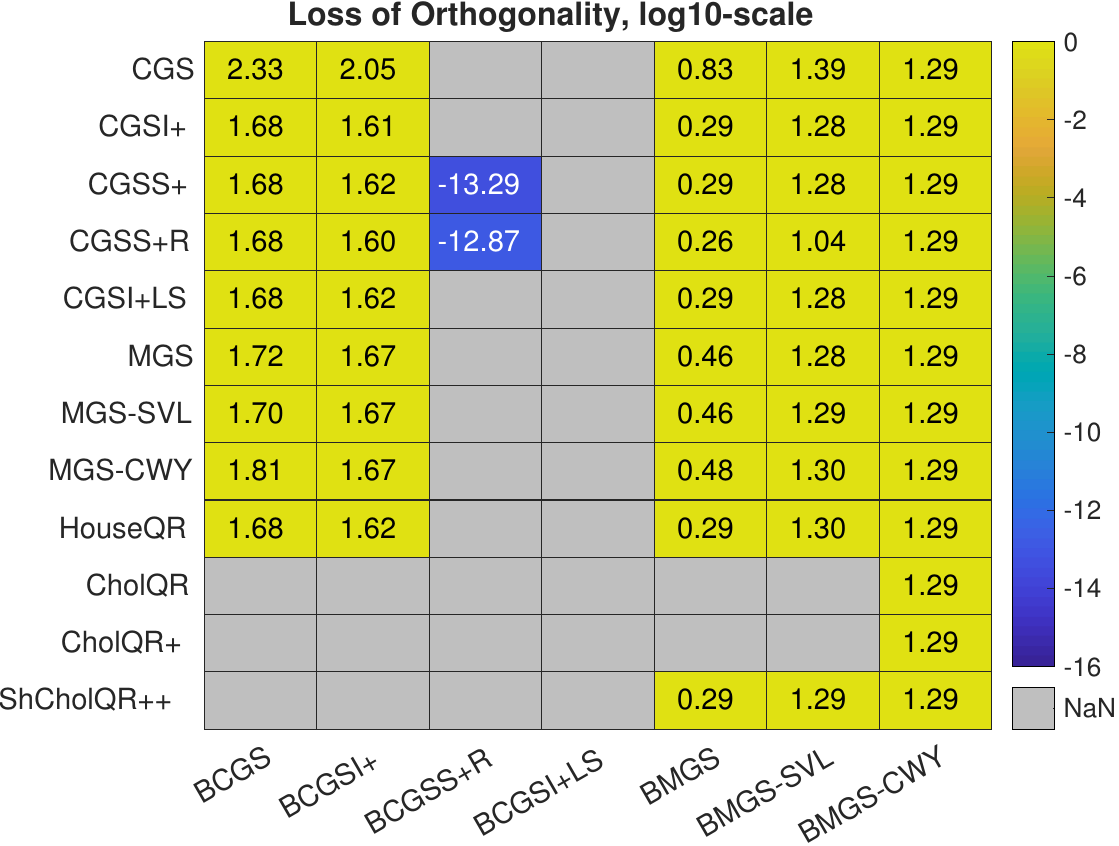} &
			\includegraphics[width=.45\textwidth]{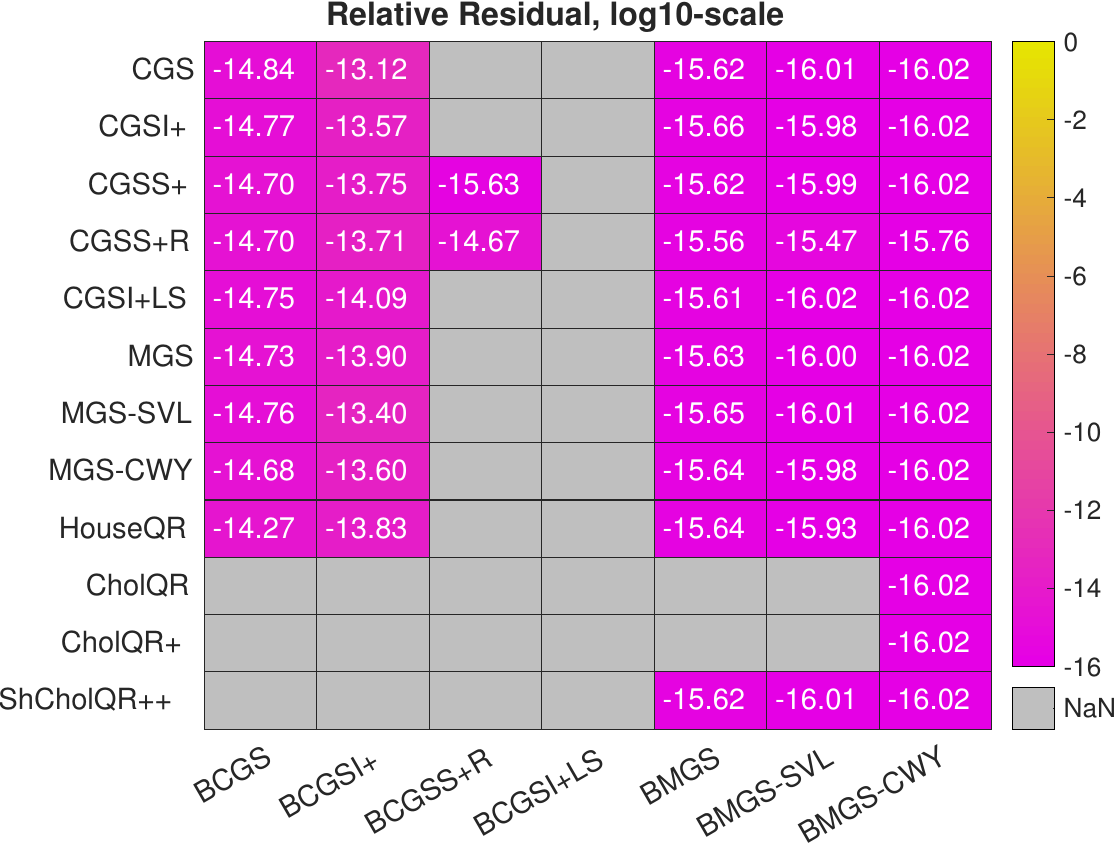}
		\end{tabular}
	\end{center}
\end{figure}

\begin{figure}[htbp!]
	\begin{center}
		\caption{Measurements for \newton \label{fig:newton}}
		\begin{tabular}{cc}
			\includegraphics[width=.45\textwidth]{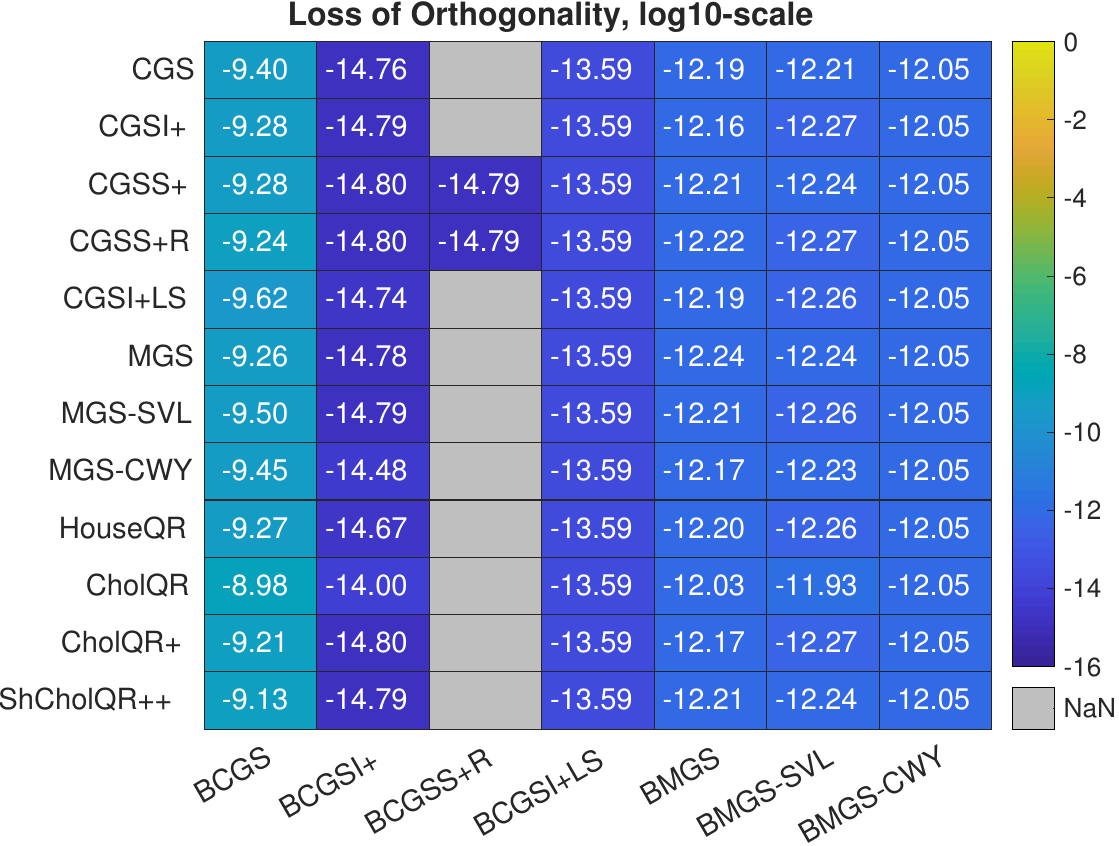} &
			\includegraphics[width=.45\textwidth]{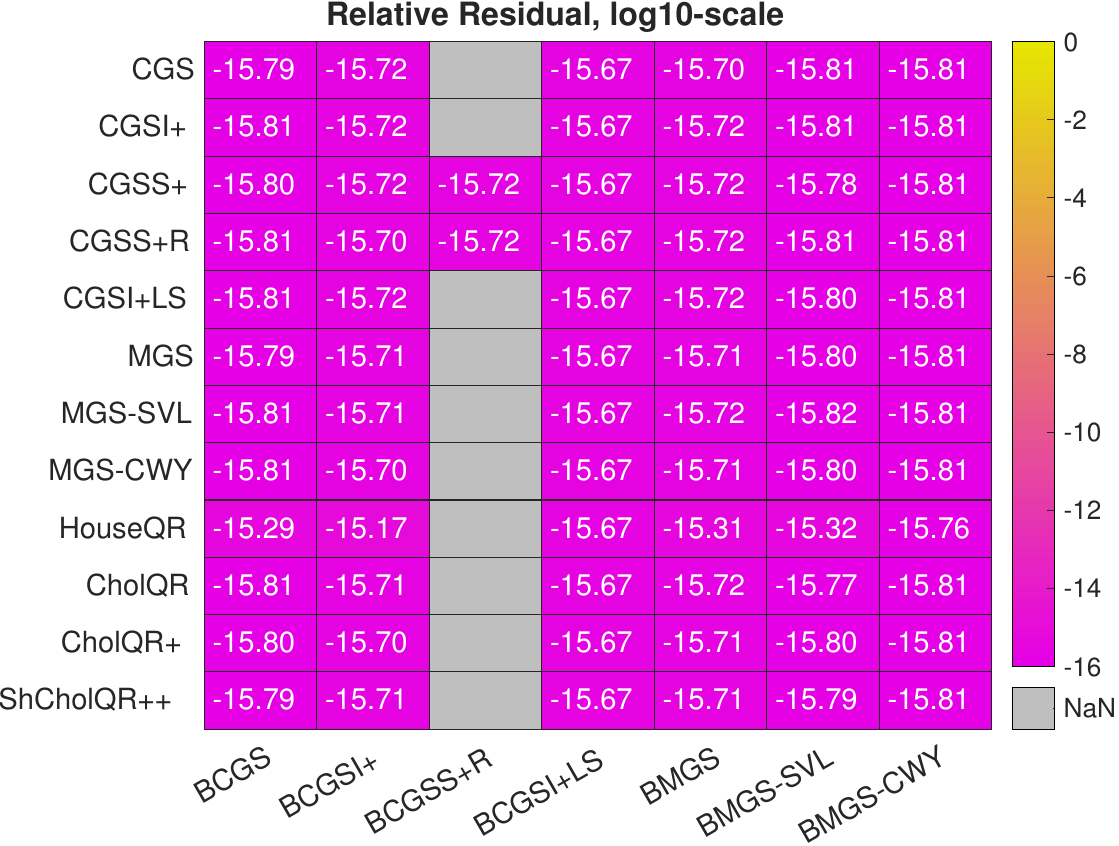}
		\end{tabular}
	\end{center}
\end{figure}

\subsection{\stewart} \label{sec:stewart}
\begin{figure}[htbp!]
	\begin{center}
		\caption{Measurements for \stewart \label{fig:stewart}}
		\begin{tabular}{cc}
			\includegraphics[width=.45\textwidth]{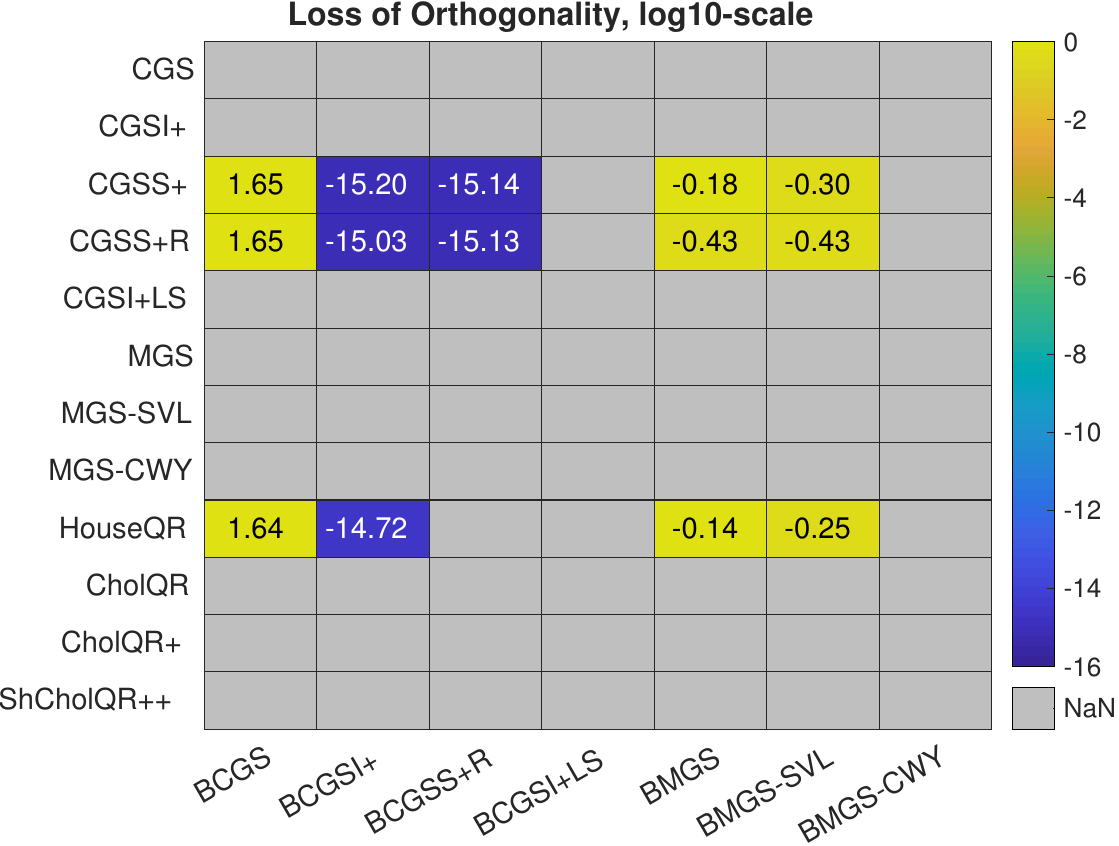} &
			\includegraphics[width=.45\textwidth]{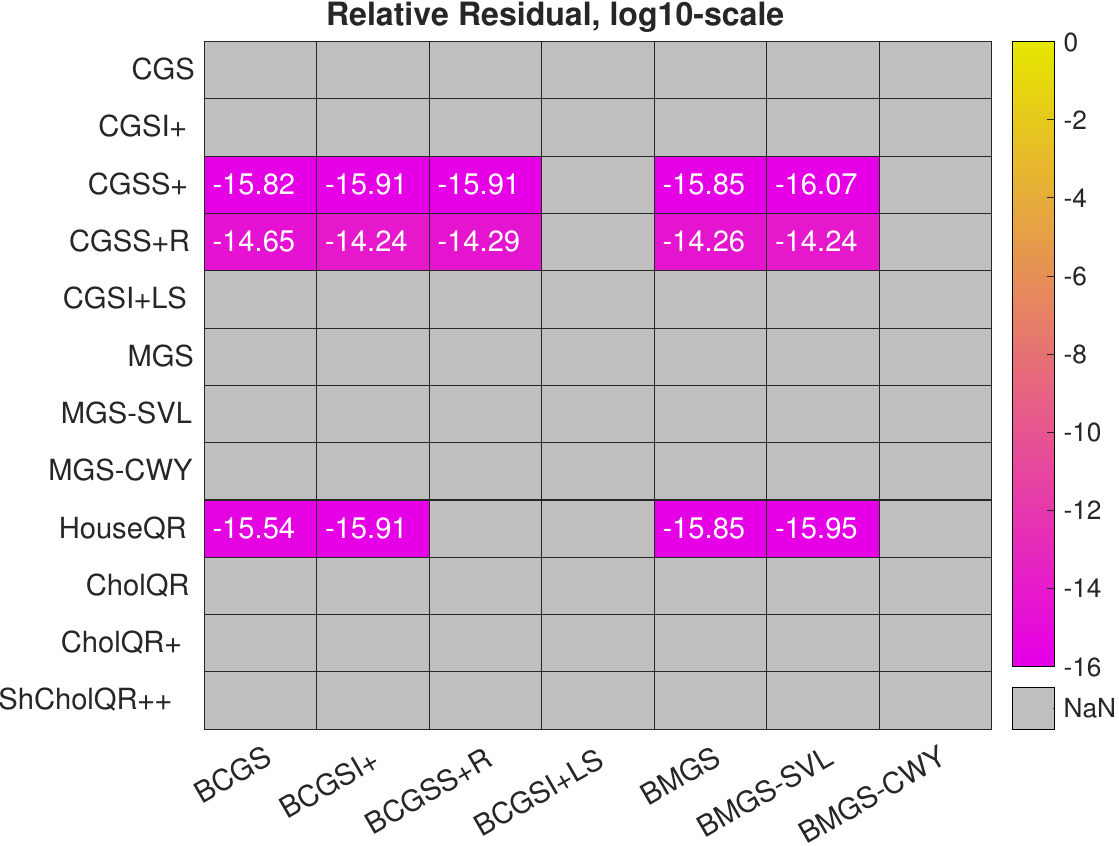}
		\end{tabular}
	\end{center}
\end{figure}

Stewart \cite{Stewart2008} purposefully designed his algorithms to be robust for pathologically bad matrices, and the results in Figure~\ref{fig:stewart} demonstrate the failure of most algorithms for such matrices.  Because of the identically zero column, nearly all methods encounter a \texttt{NaN} at some point after division of $0$ by $0$, which then spoils the rest of the algorithm.  The only methods robust enough for such situations is Stewart's own algorithm, $\BCGSSROR \circ \CGSSROR$, its variant without replacement, and \BCGSIRO with Stewart's \IOnoargs and \HouseQR.

\subsection{\stewartextreme} \label{sec:stewart_extreme}
\begin{figure}[htbp!]
	\begin{center}
		\caption{Measurements for \stewartextreme \label{fig:stewart_extreme}}
		\begin{tabular}{cc}
			\includegraphics[width=.45\textwidth]{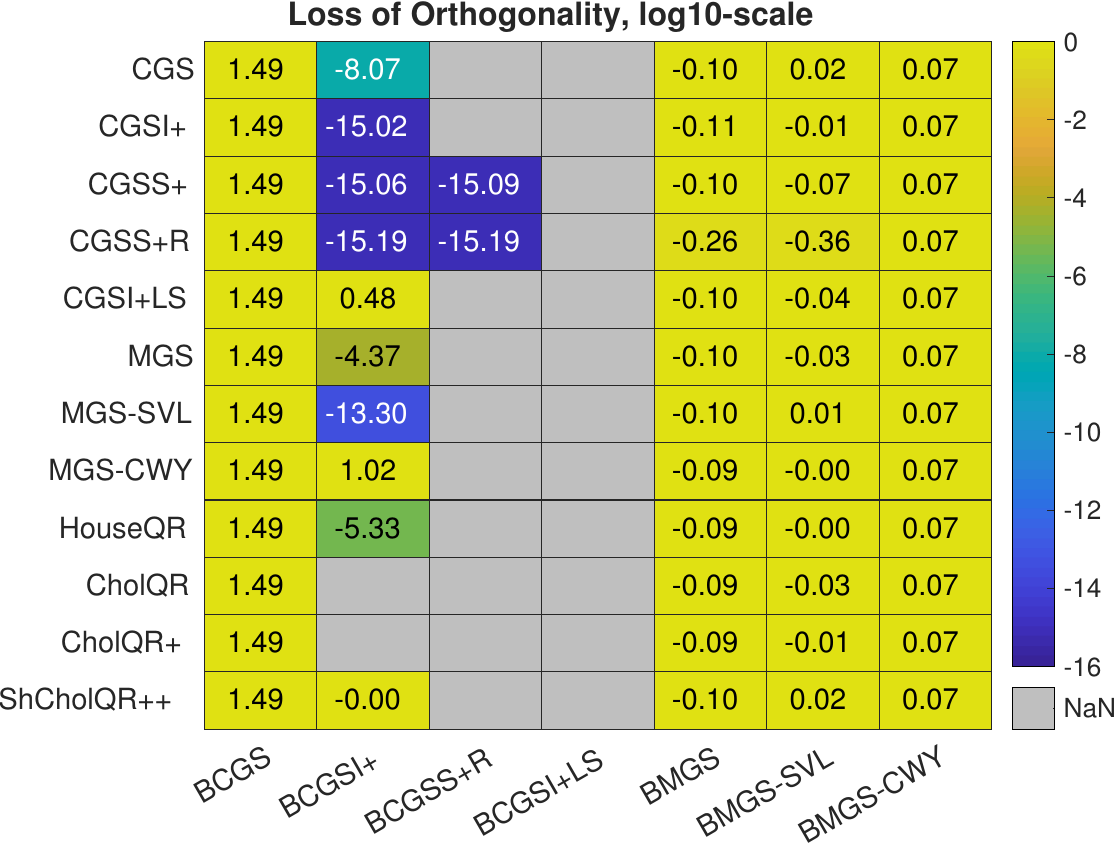} &
			\includegraphics[width=.45\textwidth]{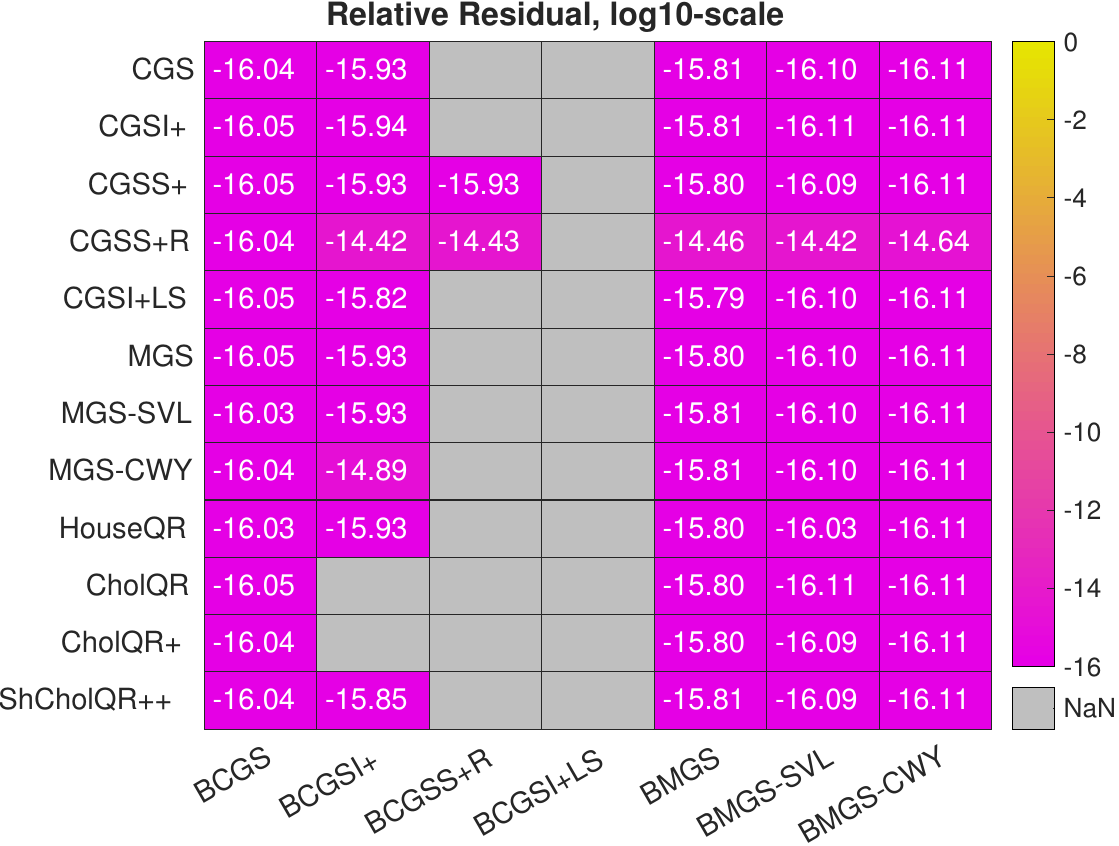}
		\end{tabular}
		\\
		\includegraphics[width=.2\textwidth]{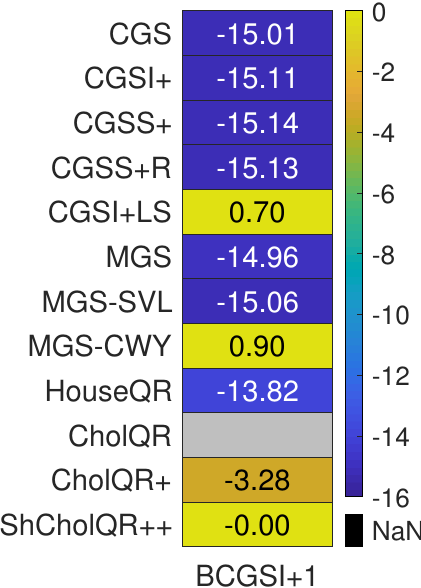}
	\end{center}
\end{figure}

This matrix in our collection demonstrates the utility of aggressive replacement.  Figure~\ref{fig:stewart_extreme} collects the results for \stewartextreme.

For both \BCGSIRO and \BCGSSROR we can see how aggressive replacement can make a significant difference. The only algorithms capable of achieving $\bigO{\eps}$ loss of orthogonality in this case are \BCGSSROR with \CGSSRO and \CGSSROR, and \BCGSIRO with \CGSSRO, \CGSSROR, and \CGSIRO. It is unclear why for this particular case $\BCGSIRO \circ \CGSIRO$ exhibits better orthogonality than $\BCGSIRO\circ\HouseQR$. 

Looking closely at Algorithm~\ref{alg:BCGSIRO}, we see that \BCGSIRO does not reorthogonalize the first block vector in line 2, likely because this is not done for the column-wise algorithm, where it is unnecessary.  By simply rerunning \IOnoarg for the first block vector (denoted ``BCGSI+1'' in Figure~\ref{fig:stewart_extreme}), we can reduce the loss of orthogonality to $\bigO{\eps}$ for some muscles. Notably both the one-sync methods-- \CGSIROLS and \MGSCWY-- still struggle, and \CholQR still cannot run to completion.

This example highlights a problem that remains unaddressed in \cite{FukayaKannanNakatsukasa2020} for \ShCholQRRORO: when $\vX$ contains a column very close to the zero vector.  In this situation, the bounds for the shift $\sigma$ may be too stringent to allow $\vX^T\vX + \sigma I_s$ to be numerically positive definite.  In the \stewartextreme example, setting $\sigma = \norm{\vX}_2^2$ allows $\BCGSIRO \circ \ShCholQRRORO$ to run to completion with $\bigO{10^{-12}}$ loss of orthogonality and $\bigO{10^{-17}}$ residual.  It seems that the upper bound $\sigma \leq \frac{1}{100} \norm{\vX}_2^2$ is more likely an artifact of the analysis and that a larger bound may in fact be tolerable.  To be fully robust, $\ShCholQRRORO$ should be written to allow for an automatically tuned bound based on whether the shifted Gramian is found to be positive definite enough for \chol to run.  This remains a topic of future work.

We note that this trick of reorthogonalizing the first block vector has even better results for other test problems we tried. For both \monomial and \laeuchli, if we reorthogonalize the first block vector then \BCGSIRO with any of the tested muscles provides $\bigO{\eps}$ loss of orthogonality.  We do not display these tests here, but see a related discussion in Section~\ref{sec:BCGSIRO_stab}.

Overall, the experiments in this section give indication that different skeletons have very different requirements on the muscle. For \BCGS, the choice of muscle makes little difference. This is because \BCGS is itself unstable as a skeleton, so there is no benefit to using, say, \HouseQR instead of \CGS (as long as the constraints on condition number are satisfied to make the muscle viable). In constrast, the reorthogonalization in \BCGSIRO stablizes the skeleton, and thus the quality of the \IOnoarg will have a greater affect on the resulting stability. This is also true for \BMGS and variants like \BMGSSVL. We elaborate on these details in the following section. 

\section{Skeleton stability in terms of muscle stability: details} \label{sec:bgs_stab_2}
We now consider the stability of each skeleton in more detail.  We utilize a common plot format, wherein stability quantities like loss of orthogonality or relative residual are plotted against varying condition numbers.  We refer to such plots as \emph{$\kappa$-plots}, where $\kappa$ refers to the fact that we study stability with respect to changes in condition number.  The simplest version of these plots takes a series of matrices $\bXX_t = \bUU \Sigma_t \VV^T \in \spR^{m \times ps}$, where $\bUU \in \spR^{m \times ps}$ is orthonormal, $\Sigma_t \in \spR^{ps \times ps}$ is a diagonal matrix whose entries are drawn from the logarithmic interval $10^{[-t, 0]}$, and $\VV \in \spR^{ps \times ps}$ is unitary.  We also consider \glued $\kappa$-plots, where the matrices are instead built as the ``glued" matrices from \cite{SmoktunowiczBarlowLangou2006}; see \ref{sec:scripts} for how they are generated.  The \monomial matrices from Section~\ref{sec:bgs_stab_1} also lend themselves to $\kappa$-plots, because varying condition numbers can be induced by varying block sizes.  All $\kappa$-plots are run in \MATLAB 9.8.0.1451342 (R2020a) Update 5. 
We provide the command-line calls for each figure in \ref{sec:scripts}.

\subsection{\BCGS} \label{sec:BCGS_stab}
\begin{figure}[htbp!]
	\begin{center}
		\caption{\glued $\kappa$-plots for P-variants of \BCGS. \label{fig:bcgs_p_kappa_plot}}
		\begin{tabular}{cc}
			\includegraphics[trim={4cm 8cm 4cm 8cm},clip,width=.45\textwidth]{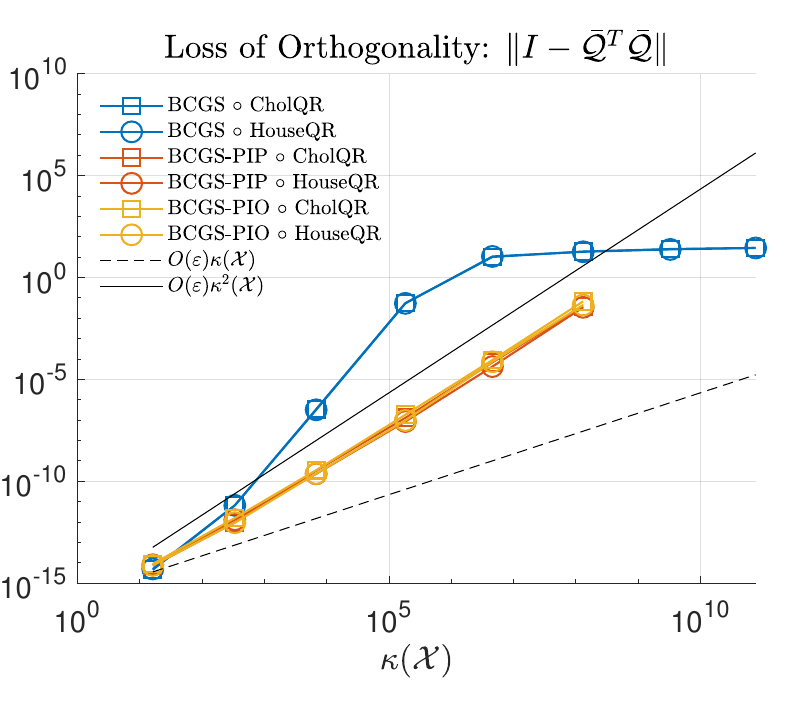} &
			\includegraphics[trim={4cm 8cm 4cm 8cm},clip,width=.45\textwidth]{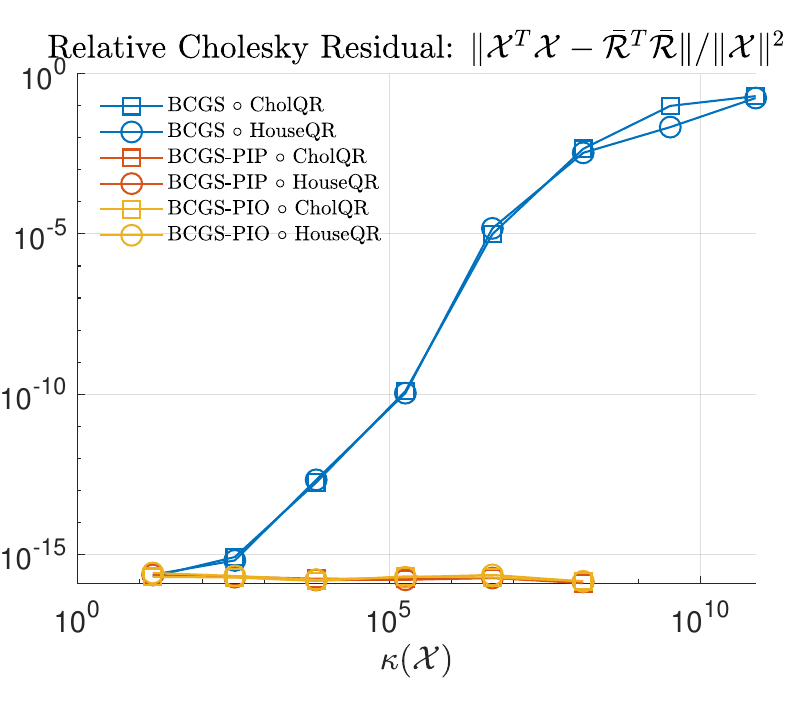}
		\end{tabular}
	\end{center}
\end{figure}

Using similar techniques as for \CGS, we conjecture that the loss of orthogonality in $\BCGS \circ \IOnoarg$ can be bounded as
\[
\norm{I - \bQQbar^T \bQQbar}_2 \leq \bigO{\eps} \kappa(\bXX)^{n-1}
\]
as long as $\bigO{\eps}\kappa(\bXX)<1$. 
A full proof that the \BCGS $P$-variants have $\bigO{\eps} \kappa^2(\bXX)$ loss of orthogonality and $\bigO{\eps}$ relative Cholesky residual as long as $\bigO{\eps}\kappa(\bXX)<1/2$ is given in \cite{CarsonLundRozloznik2021}.  These bounds hold as long as $\IOnoarg(\vX)$ satisfies
\begin{equation} \label{eq:chol_IO}
	\Rbar^T \Rbar = \vX^T \vX + E, \quad \norm{E} \leq \bigO{\eps} \norm{\vX}^2,
\end{equation}
and
\begin{equation} \label{eq:res_IO}
	\vQbar \Rbar = \vX + \vD, \quad \norm{\vD} \leq \bigO{\eps} (\norm{\vX}+ \norm{\vQbar}\norm{\Rbar}).
\end{equation}
We note that these constraints on the \IOnoarg are relatively relaxed, and are satisfied by almost any reasonable choice of \IOnoarg, from \CholQR to \HouseQR. Indeed, in Figure \ref{fig:bcgs_p_kappa_plot}, we see that both \BCGSPIP and \BCGSPIO (as well as \BCGS) behave similarly for both \CholQR and \HouseQR. We also note that it is clear that \BCGS (without the Pythagorean fix) does not satisfy the $\bigO{\eps}\kappa^2{\bXX}$ bound on loss of orthogonality; as the relative Cholesky residual departs from $\bigO{\eps}$, the loss of orthogonality follows suit.   

\subsection{\BCGSIRO} \label{sec:BCGSIRO_stab}
Barlow and Smoktunowicz prove an $\bigO{\eps}$ loss of orthogonality for \BCGSIRO and an unconditionally stable \IOnoarg, as long as $\bigO{\eps} \kappa(\bXX) < 1$ \cite{BarlowSmoktunowicz2013}.  The key to their approach is to first obtain bounds for the following subproblem: given a near left-orthogonal matrix $\bUU \in \spR^{m \times t}$ and a matrix $\vB \in \spR^{m \times s}$, find an $\vS \in \spR^{t \times s}$, $R \in \spR^{s \times s}$ upper triangular, and $\vQ$ left orthogonal such that
\begin{equation} \label{eq:BCGSIRO_step}
	\vB = \bUU \vS_{\vB} + \vQ_{\vB} R_{\vB} \mbox{ and } \bUU^T \vQ_{\vB} \approx 0.
\end{equation}
This subproblem is present at every step of \BCGSIRO, as long as the previous step has produced a near-left orthogonal matrix $\bUU$.  With the unconditionally stable \IOnoarg guaranteeing that $\vQ_1$ (see line~2 in Algorithm~\ref{alg:BCGSIRO}) is near left-orthogonal, induction over bounds on the subproblem leads to the desired result.

One potential drawback to Barlow and Smoktunowicz's approach is the \emph{a posteriori} assumption that $\bigO{\eps} \norm{\vB} \norm{\Rbar_{\vB}^\inv} < 1$, where $\Rbar_{\vB}$ is the computed version of $R_{\vB}$ from the subproblem \eqref{eq:BCGSIRO_step}; see \cite[Equation~(3.27)]{BarlowSmoktunowicz2013}. In addition to simplifying the proof, the authors argue that such an assumption is useful in practice, because the computed quantities can be measured on the fly.

\begin{figure}[htbp!]
	\begin{center}
		\caption{standard $\kappa$-plots for \BCGSIRO, in comparison to \BCGS. \label{fig:bcgs_iro_kappa_plot}}
		\begin{tabular}{cc}
			\includegraphics[trim={4cm 8cm 4cm 8cm},clip,width=.45\textwidth]{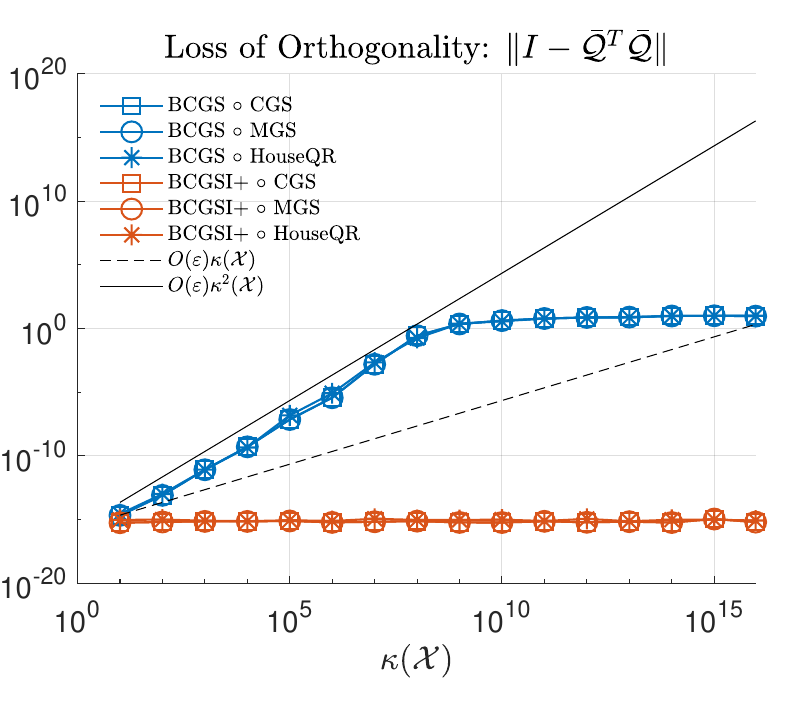} &
			\includegraphics[trim={4cm 8cm 4cm 8cm},clip,width=.45\textwidth]{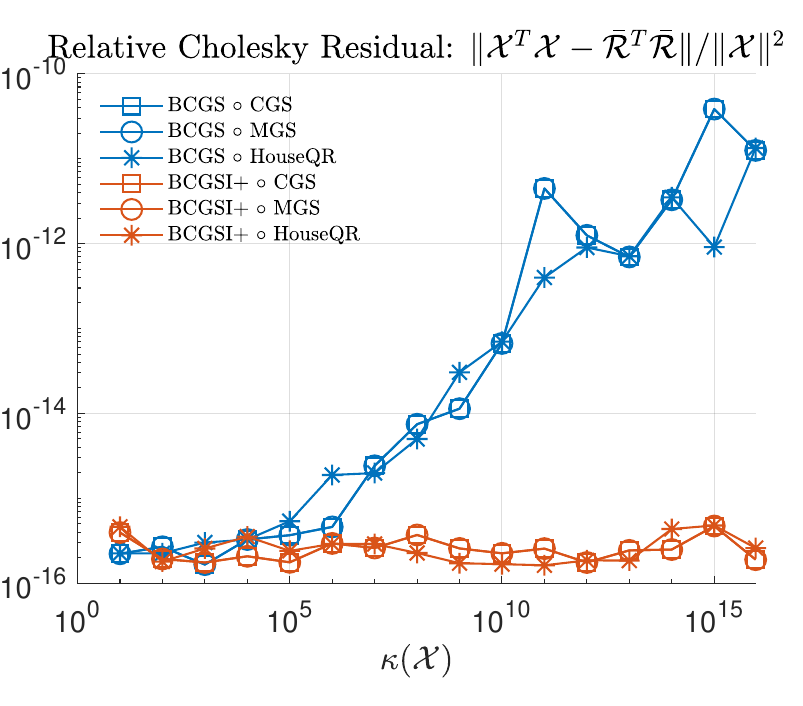}
		\end{tabular}
	\end{center}
\end{figure}
\begin{figure}[htbp!]
	\begin{center}
		\caption{Laeuchli $\kappa$-plots for \BCGSIRO, in comparison to \BCGS. \label{fig:bcgs_iro_laeuchli}}
		\begin{tabular}{cc}
			\includegraphics[trim={4cm 8cm 4cm 8cm},clip,width=.45\textwidth]{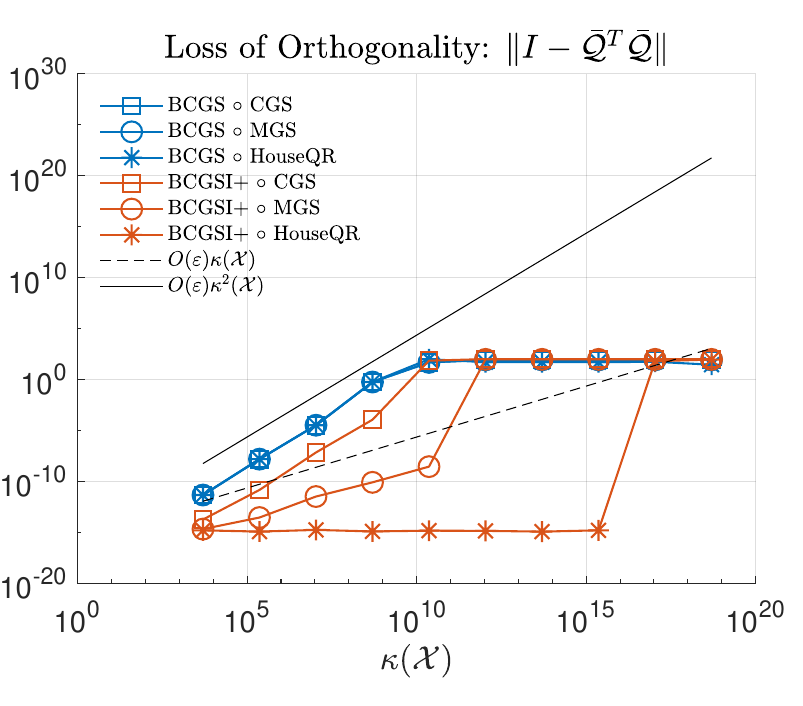} &
			\includegraphics[trim={4cm 8cm 4cm 8cm},clip,width=.45\textwidth]{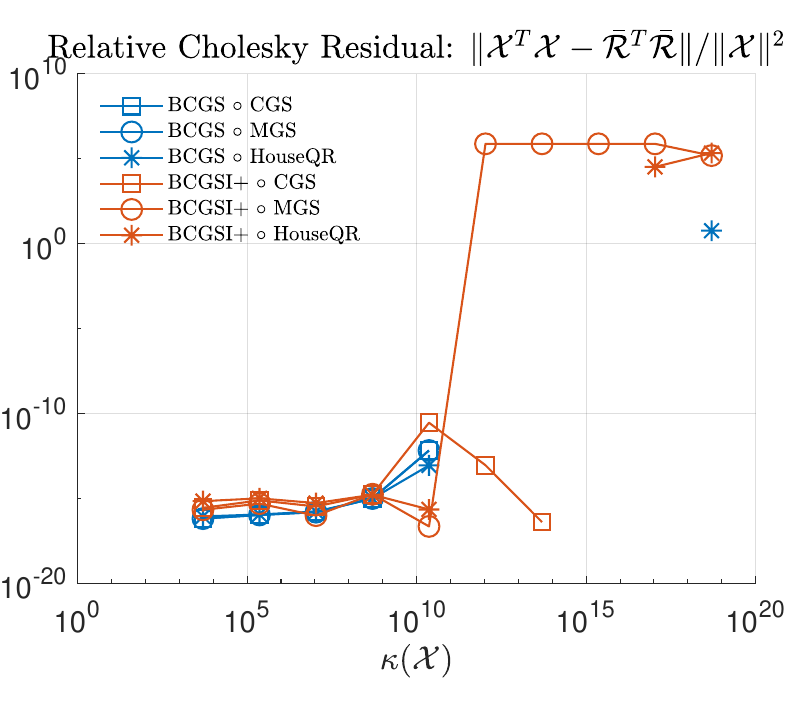}
		\end{tabular}
	\end{center}
\end{figure}

Figure~\ref{fig:bcgs_iro_kappa_plot} displays standard $\kappa$-plots for \BCGSIRO in comparison to \BCGS for muscles \CGS, \MGS, and \HouseQR.  It is tempting to conclude that the skeletons' loss of orthogonality is independent of the muscle, but recalling the results for \laeuchli and \monomial matrices in Sections~\ref{sec:laeuchli} and \ref{sec:monomial} (both of which satisfy $\bigO{\eps} \kappa(\bXX) <1$; cf.\ Table~\ref{tab:matrix_props}), it is clear that we may need an unconditionally stable \IOnoarg for some tough matrices. Figure \ref{fig:bcgs_iro_laeuchli} shows $\kappa$-plots generated using \laeuchli matrices. Indeed, here we see highlighted the different ways that \BCGS and \BCGSIRO are affected by the \IOnoarg. Whereas \BCGS behaves the same regardless of what \IOnoarg is used, \BCGSIRO is much more sensitive. 

\subsection{\BCGSSROR} \label{sec:BCGSSROR_stab}
Although Stewart's development of \BCGSSROR \cite{Stewart2008} is largely driven by stability considerations, a rigorous proof is lacking.  We conjecture that for any $\bXX$, even those with rank deficiency, we have $\bigO{\eps}$ loss of orthogonality, up to a small constant accounting for the replacement
tolerance $\rpltol$.

This conjecture is based largely on numerical results from both Stewart's paper \cite{Stewart2008} and our own studies in Section~\ref{sec:bgs_stab_1}, wherein a clear trade-off between loss of orthogonality and relative residual seems to be driven largely by $\rpltol$.  It is also clear that the orthogonalization fault step, designed precisely to avoid situations leading to instability, plays an important role. A successful stability analysis for \BCGSSROR of course also depends on one for \CGSSROR.

\subsection{\BCGSIROLS} \label{sec:BCGSIROLS_stab} Rigorous stability
analysis for the low-sync algorithms from \cite{SwirydowiczLangouAnanthan2021}, especially their block versions, remains open.  We can, however, formulate some conjectures.

\begin{figure}[htbp!]
	\begin{center}
		\caption{standard $\kappa$-plots for low-sync \BCGSIROLS, in comparison to those of \BCGS and \BCGSIRO. \label{fig:bcgs_iro_ls_kappa_plot}}
		\begin{tabular}{cc}
			\includegraphics[trim={4cm 8cm 4cm 8cm},clip,width=.45\textwidth]{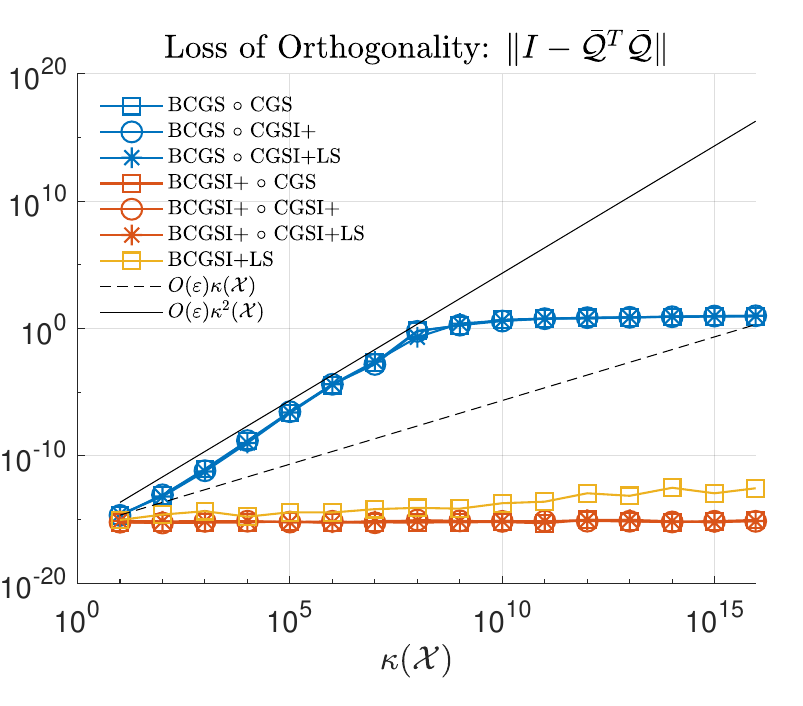} &
			\includegraphics[trim={4cm 8cm 4cm 8cm},clip,width=.45\textwidth]{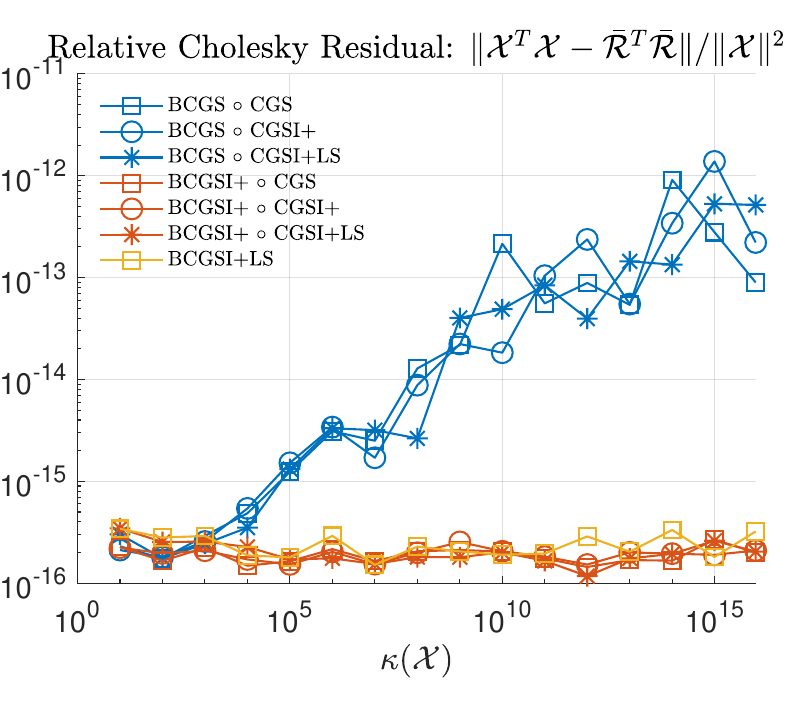}
		\end{tabular}
	\end{center}
\end{figure}

\begin{figure}[htbp!]
	\begin{center}
		\caption{\laeuchli $\kappa$-plots for low-sync \BCGSIROLS, in comparison to those of \BCGS and \BCGSIRO. \label{fig:bcgs_iro_ls_laeuchli_kappa_plot}}
		\begin{tabular}{cc}
			\includegraphics[trim={4cm 8cm 4cm 8cm},clip,width=.45\textwidth]{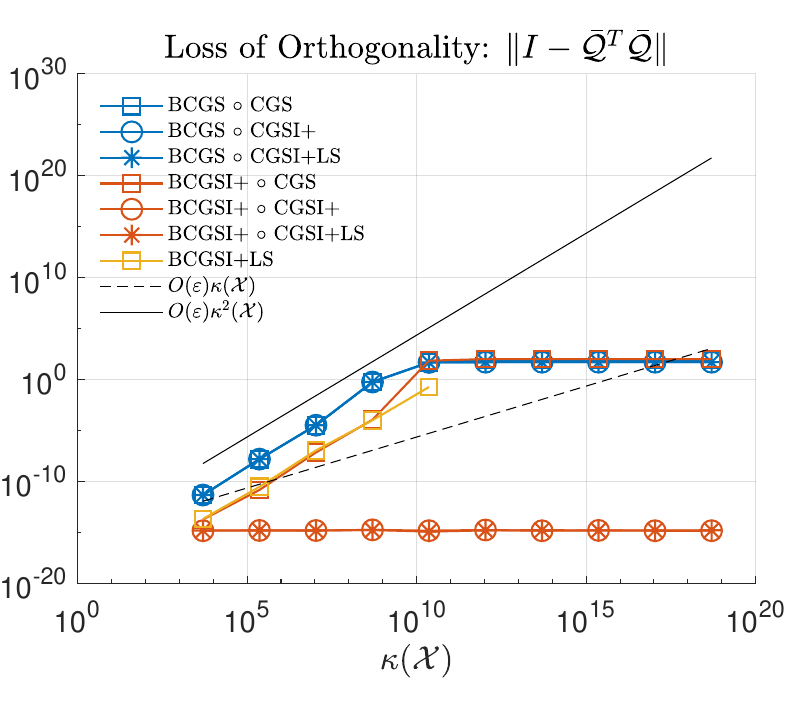} &
			\includegraphics[trim={4cm 8cm 4cm 8cm},clip,width=.45\textwidth]{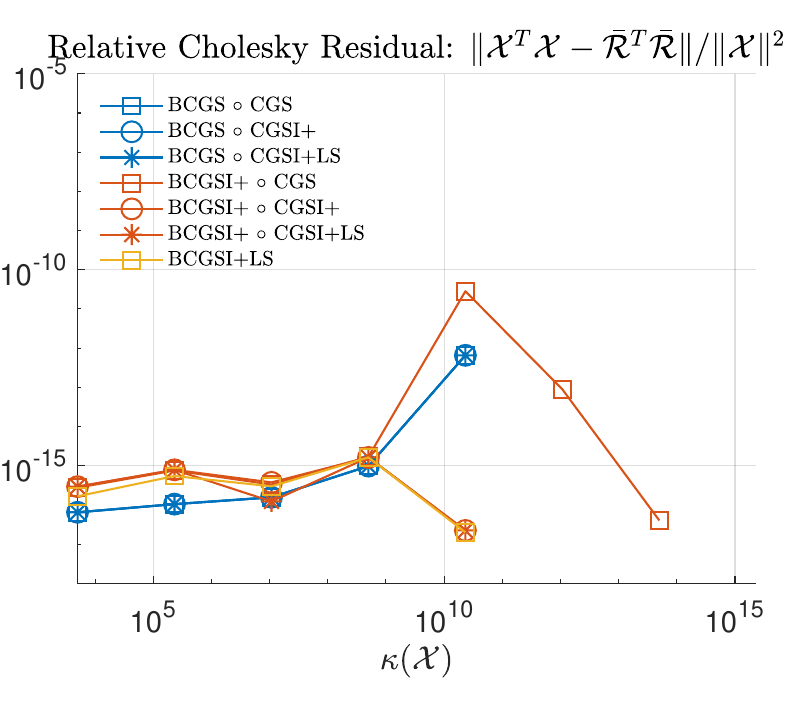}
		\end{tabular}
	\end{center}
\end{figure}

In Figure~\ref{fig:bcgs_iro_ls_kappa_plot}, we see that indicate that \BCGSIROLS remains relatively stable, especially its relative Cholesky residual.  The loss of orthogonality for \BCGSIROLS begins to deviate gradually from $\bigO{\eps}$ in Figure~\ref{fig:bcgs_iro_ls_kappa_plot}. However, looking at the more challenging \laeuchli matrix problems in Figure~\ref{fig:bcgs_iro_ls_laeuchli_kappa_plot}, we see that the loss of orthogonality can be much more significant, exceeding $\bigO{\eps}\kappa(\bXX)$. Based on these results, we conjecture that \BCGSIROLS satisfies the bound $\bigO{\eps}\kappa^2(\bXX)$ as long as $\bigO{\eps}\kappa^2(\bXX) < 1$ (due to the use of \chol within the algorithm). This conjecture is stated in Table \ref{tab:skel_upper_bounds}; a formal proof remains future work. The norm lagging and delayed reorthogonalization seems to have a much more significant effect in \BCGSIROLS versus \CGSIROLS.

\subsection{\BMGS} \label{sec:BMGS_stab}
What appears to be the earliest study on the stability of a block Gram-Schmidt method was conducted by Jalby and Philippe in 1991 \cite{JalbyPhilippe1991}. They focus on $\BMGS \circ \MGS$ and $\BMGS \circ \MGSRO$, and note how the underlying \CGS-like nature of \BMGS makes it prone to instability.  To see this, observe that lines~4-8 of Algorithm~\ref{alg:BMGS} can be written more succinctly (in exact arithmetic) as
\begin{equation} \label{eq:BMGS_like_CGS}
	\vW = (I - \vQ_k \vQ_k^T) \cdots (I - \vQ_1 \vQ_1^T)\vX_{k+1}.
\end{equation}
Because each $\vQ_j$, $j = 1, \ldots, k$, here is a tall-and-skinny matrix, the projector $I - \vQ_j \vQ_j^T$ is equivalent to a step of \CGS.  The authors demonstrate this intuition rigorously with \MGS as the \IOnoarg to arrive at the following bound for the loss of orthogonality:
\begin{equation} \label{eq:BMGS_MGS_ortho_bound}
	\norm{I - \bQQbar^T \bQQbar}_2 \leq \bigO{\eps} C \normF{\bXX} \norm{\RR^\inv}_2 \leq \bigO{\eps} \kappa^2(\bXX),
\end{equation}
wherein the second inequality follows because $C \leq O(\kappa(\bXX))$ and we have $\normF{\bXX} \norm{\RR^\inv}_2 \leq O(\kappa(\bXX))$. The authors further show that the constant $C$ becomes $O(1)$ if \MGS is replaced by \MGSRO, thus giving the MGS-like bound
\begin{equation} \label{eq:BMGS_MGSRO_ortho_bound}
	\norm{I - \bQQbar^T \bQQbar}_2 \leq \bigO{\eps} \kappa(\bXX).
\end{equation}
Barlow \cite{Barlow2019} conjectures that $\BMGS \circ \HouseQR$ is as stable as \MGS but does not prove so explicitly.  A careful look at Jalby and Philippe's proof for their Theorem~4.1 reveals that there is not much to prove; in fact, the very same observation they used to derive \eqref{eq:BMGS_MGSRO_ortho_bound} holds for any \IOnoarg that is unconditionally stable.

We sketch the proof of this here. 
At the top of page 1062 of \cite{JalbyPhilippe1991}, the Jalby and Philippe directly apply \MGS bounds from \cite{Bjorck1967} and define the constant $C = \max_{1 \leq k \leq p-1} \normF{\vW_k} \norm{R_{\vW_k}^\inv}_2$, where $\vW_k$ corresponds to the block vector computed up to line~8 of Algorithm~\ref{alg:BMGS}, and $R_{\vW_k}$ its upper triangular factor from an exact QR decomposition.  The
intermediate constants $\normF{\vW_k} \norm{R_{\vW_k}^\inv}_2$ are absent for any unconditionally stable \IOnoarg.  Tracking $C$ all the way through the proof of \cite[Theorem~4.1]{JalbyPhilippe1991} (by which \eqref{eq:BMGS_MGS_ortho_bound} holds) verifies that, without changing any other lines in their text, we would obtain \eqref{eq:BMGS_MGSRO_ortho_bound} for $\BMGS \circ \IOnoarg$, assuming \IOnoarg is unconditionally stable.

We remark that, to our knowledge, the present work is the first to state bounds for $\BMGS$ with an unconditionally stable \IOnoarg (such as \HouseQR or \TSQR), despite the ubiquity of $\BMGS \circ \HouseQR$ in practice.  More precise bounds for particular \IOnoargs remain open, but can be easily derived on a case-by-case basis from the transparent analysis in \cite{JalbyPhilippe1991}.

Because $\BMGS \circ \HouseQR$ behaves essentially like \MGS, we conjecture that results analogous to those by Paige, Rozlo\v{z}n\'ik, and Strako\v{s} \cite{PaigeRozloznikStrakos2006} hold for block GMRES, which is usually implemented with $\BMGS \circ \HouseQR$.  To be more specific, the loss of orthogonality in $\BMGS \circ \HouseQR$ is tolerable when the end application is the solution of a linear system via the GMRES method.

\begin{figure}[htbp!]
	\begin{center}
		\caption{$\kappa$-plots for \BMGS.  (a) standard $\kappa$-plot.  (b) \laeuchli $\kappa$-plot. \label{fig:bmgs_kappa_plot}}
		\begin{tabular}{cc}
			\includegraphics[trim={4cm 8cm 4cm 8cm},clip,width=.45\textwidth]{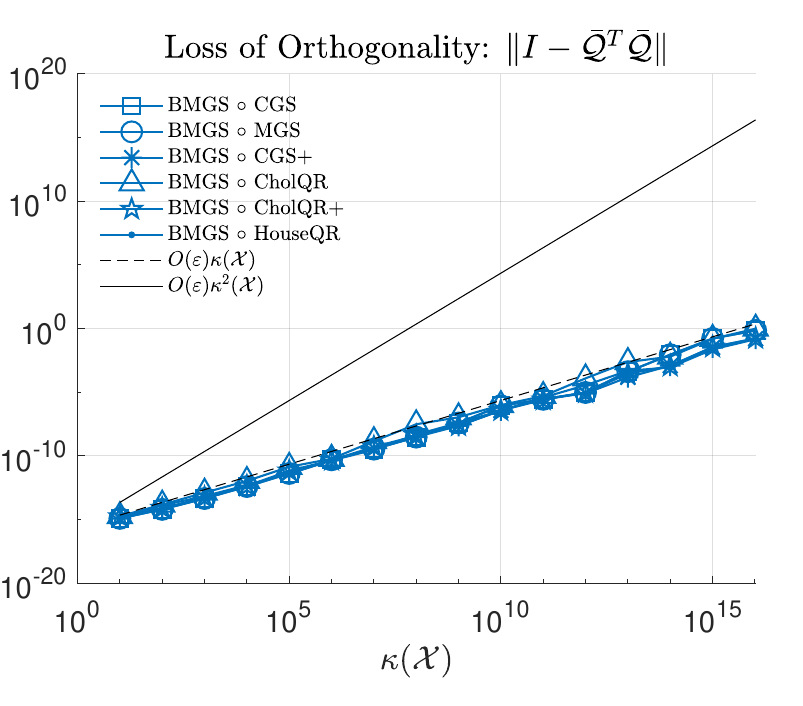} &
			\includegraphics[trim={4cm 8cm 4cm 8cm},clip,width=.45\textwidth]{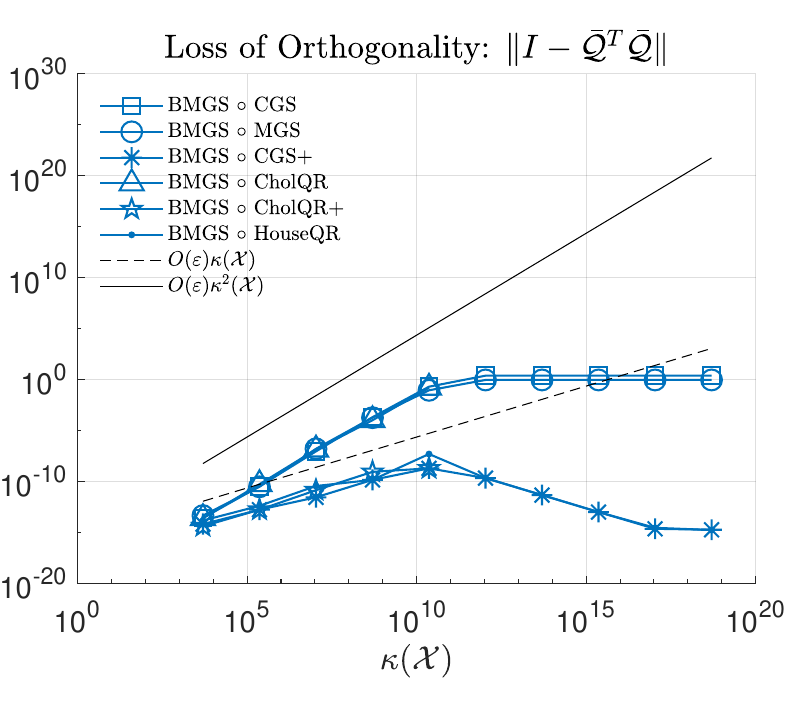} \\
			{\small (a)}	& {\small (b)}
		\end{tabular}
	\end{center}
\end{figure}

Figure~\ref{fig:bmgs_kappa_plot} demonstrates the \BMGS behavior for a variety of muscles and their reorthogonalized variants again for both standard $\kappa$-plots (left) and \laeuchli $\kappa$-plots (right). For the standard $\kappa$-plots, all variants behave the same, giving better orthogonality than indicated by \eqref{eq:BMGS_MGS_ortho_bound}. For the \laeuchli $\kappa$-plots, we do observe behavior like \eqref{eq:BMGS_MGS_ortho_bound} here for the \IOnoargs that do not have $\bigO{\eps}$ loss of orthogonality (\CGS, \MGS, and \CholQR). Indeed, this example demonstrates that an unconditionally stable muscle is really necessary to guarantee $\bigO{\eps} \kappa(\bXX)$ loss of orthogonality.

\subsection{Low-sync \BMGS variants} \label{sec:BMGS_low_sync_stab}
To date, we are only aware of Barlow's stability analysis of \MGSSVL and \BMGSSVL for low-synchronization \BMGS variants \cite{Barlow2019}.  Indeed, it appears that Barlow was even unaware that his algorithms could take on so many other forms, particularly the work by \'{S}wirydowicz et al. \cite{SwirydowiczLangouAnanthan2021}.

Barlow's analysis relies heavily on the Sheffield observation and Schreiber-van-Loan reformulation-- what is sometimes also referred to as the ``augmented" problem-- for the stability analysis of \BMGSSVL; see, in
particular, \cite[Section~4]{Barlow2019}.  The following quantities play key roles:
\begin{align*}
	\resTTSS & := \TTbar \SS - I, \\
	\resQQRR & := \bQQbar \RRbar - \bXX, \mbox{ and }\\
	\resTTRR & := (I - \TTbar) \RRbar.
\end{align*}
The matrix $\SS \in \spR^{n \times n}$ is an exact quantity that does not explicitly arise in the computation.  Although it has several equivalent formulations in exact arithmetic (see especially page 1261 in \cite{Barlow2019}), it is apparently defined as the upper triangular part of the exact multiplication between $\bQQbar^T$ and $\bQQbar$, i.e.,
\[
\SS := \triu(\bQQbar^T \bQQbar).
\]
In exact arithmetic, $\SS = \TT^\inv$, hence why the measure $\resTTSS$ bears significance.

For $\BMGSSVL \circ \IOnoarg$, Barlow shows that
\begin{align}
	\normF{\resTTSS} & \leq \bigO{\eps}, \label{eq:resTTSS} \\
	\normF{\resQQRR} & \leq \bigO{\eps} \normF{\bXX}, \mbox{ and } \label{eq:resQQRR} \\
	\normF{\resTTRR} & \leq \bigO{\eps} \normF{\bXX}, \label{eq:resTTRR}
\end{align}
as long as $[\vQbar, \Rbar, \Tbar] = \IO{\vX}$ satisfies the triad
\begin{align}
	\resTS & := \Tbar S - I_s, \quad \normF{\resTS} \leq \bigO{\eps}; \label{eq:triad1} \\
	\resQR & := \vQbar \Rbar - \vX, \quad \normF{\resQR} \leq \bigO{\eps} \normF{\vX}; \mbox{ and } \label{eq:triad2}\\
	\resTR & := (I_s - \Tbar) \Rbar, \quad \normF{\resTR} \leq \bigO{\eps} \normF{\vX} \label{eq:triad3}
\end{align}
where $S := \triu(\vQbar^T \vQbar)$ and $S = T^\inv$ in exact arithmetic.

Equation~\eqref{eq:resQQRR} already gives the usual residual bound, and Barlow additionally shows that it holds at every step of \BMGS.  To arrive at bounds for loss of orthogonality, he further shows that \eqref{eq:resTTSS}-\eqref{eq:resTTRR} imply
\[
\normF{I - \bQQbar^T \bQQbar} \leq \bigO{\eps} \normF{\bXX} \norm{\RR^\inv},
\]
Recalling (from, e.g., \cite{Bjorck1967}) that $\normF{\bXX} \norm{\RR^\inv} \leq O(\kappa(\bXX))$, this latter condition can be satisfied only if $\bXX$ is sufficiently far from a rank-deficient matrix, i.e., if $\bigO{\eps} \kappa(\bXX) < 1$.  We then have that
\[
\normF{I - \bQQbar^T \bQQbar} \leq \bigO{\eps} \kappa(\bXX),
\]
and thanks to norm equivalence, we obtain the bounds reported in Table~\ref{tab:skel_upper_bounds}.

At first glance through \ref{sec:musc-alg}, it may seem that the \MGS low-sync variants are the only muscles designed to produce a $T$ matrix.  Upon further investigation of how Barlow uses the quantities $\resTS$, $\resQR$, and $\resTR$, it becomes clear that actually all other \IOnoargs implicitly produce
$\Tbar \equiv I$.  This is important for composing \BMGSSVL and \BMGSLTS with other muscles, because they explicitly update the block diagonal entries of $\TT$ with outputs from the \IOnoarg.

It then holds for these muscles that
\[
	\normF{\resTS} \equiv \normF{\triu(I - \vQbar^T \vQbar)} \leq \normF{I - \vQbar^T \vQbar},
\]
and
\[
	\normF{\resTR} = 0.
\]
Is is thus clear that \eqref{eq:triad1}-\eqref{eq:triad3} are satisfied by any \IOnoarg with $\bigO{\eps}$ loss of orthogonality and relative residual (such as \CGSIRO), and thus \eqref{eq:resTTSS}-\eqref{eq:resTTRR} hold. 
Note that although \MGSSVL does not have $\bigO{\eps}$ loss of orthogonality, it still satisfies \eqref{eq:triad1}-\eqref{eq:triad3} as long as $\Rbar$ is assumed to be nonsingular; this is proved directly by Barlow (see \cite[Section~4.2]{Barlow2019}). 
We conjecture that Barlow's framework could be easily repurposed to prove rigorous stability bounds for \MGSLTS and \BMGSLTS.  A full stability analysis for \MGSLTS, \MGSCWY, and \MGSICWY and their block variants is the subject of ongoing work.

\begin{figure}[htbp!]
	\begin{center}
		\caption{$\kappa$-plots for low-sync variants of \BMGS.  (a) standard $\kappa$-plot.  (b) \laeuchli $\kappa$-plot. \label{fig:bmgs_t_kappa_plot}}
		\begin{tabular}{cc}
			\includegraphics[trim={4cm 8cm 4cm 8cm},clip,width=.45\textwidth]{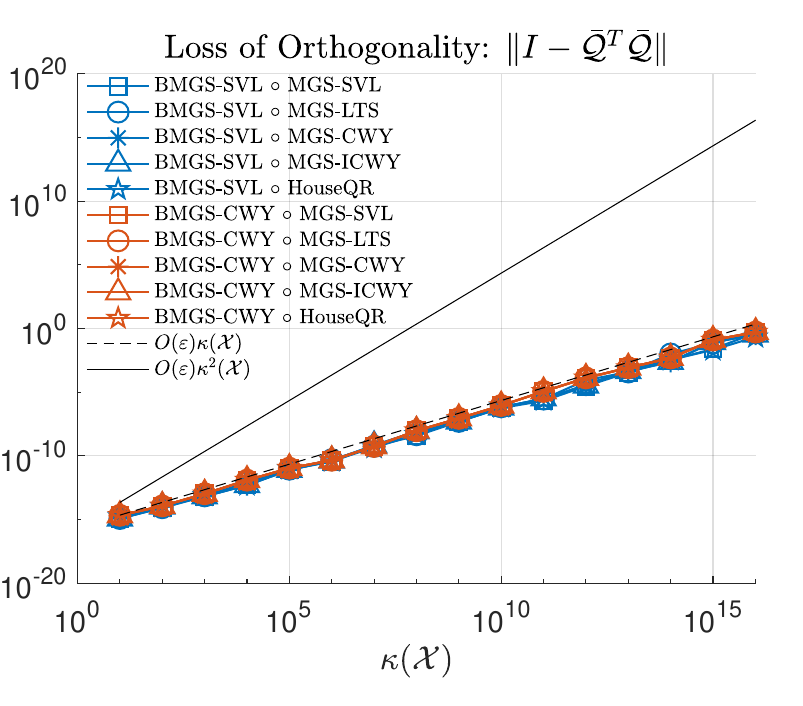} &
			\includegraphics[trim={4cm 8cm 4cm 8cm},clip,width=.45\textwidth]{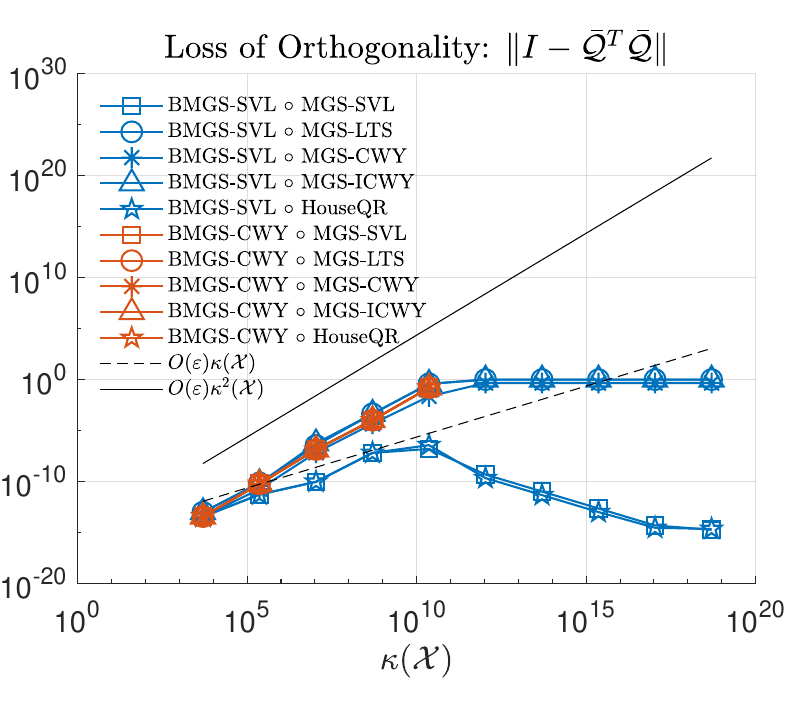} \\
			{\small (a)}	& {\small (b)}
		\end{tabular}
	\end{center}
\end{figure}

In Figure~\ref{fig:bmgs_t_kappa_plot}, we focus on just the \BMGSSVL skeleton and its one-sync version, \BMGSCWY, paired with all four ``T''-variants as muscles as well as \HouseQR. Standard $\kappa$-plots are on the left and \laeuchli $\kappa$-plots are on the right. For the standard $\kappa$-plots, the behavior of all variants is nearly identical, achieving $\bigO{\eps}$ loss of orthogonality. The more challenging \laeuchli matrices are more discerning, allowing us to draw two main conclusions. First, the T-variant skeletons and muscles are not ``mix-and-match''; \BMGSSVL works with \MGSSVL, but in combination with other T-variant muscles, the loss of orthogonality exceeds the $\bigO{\eps}\kappa(\bXX)$ bound. Second, as with \BCGSIROLS (cf.\ Section~\ref{sec:BCGSIROLS_stab}), normalization lagging in \BMGSCWY can lead to complete loss of orthogonality. Here, \BMGSCWY loses orthogonality completely after $\bigO{\eps}\kappa^2(\bXX)$ exceeds 1, even when combined with an unconditionally stable muscle like \HouseQR.



\section{Mixed-precision approaches} \label{sec:mixed_prec}
Due to the trend of commercially available mixed-precision hardware and its potential to provide computational time and energy savings, there is a growing interest in the use of mixed precision within numerical linear algebra routines. We comment briefly on existing efforts towards mixed-precision orthogonalization schemes. 

In~\cite{YamazakiTomovDongarra2015}, a mixed-precision Cholesky QR factorization (\mCholQR) is studied, in which some intermediate computations are performed at double the working precision. The authors prove that in such a mixed-precision setting, the loss of orthogonality depends only linearly (rather than quadratically) on the condition number of the input matrix under the assumption (as in \CholQR) that $\bigO{\eps} \kappa^2(\vX) <1$ (cf.\ Table~\ref{tab:musc_upper_bounds}).

Building on this work, in \cite{YamazakiTomovKurzak2015}, the authors propose using the \mCholQR of \cite{YamazakiTomovDongarra2015} as the \IOnoarg within BGS. The paper does not contain a formal error analysis, but various numerical experiments are performed to examine the loss of orthogonality and residual error when \mCholQR is combined with \BMGS and \BCGS. The experiments suggest that in some cases the numerical stability of the $\BMGS \circ \mCholQR$ variants can be maintained to the same level of \MGS if reorthogonalization is used. 

We note also that in~\cite{YangFoxSanders2019}, the finite-precision behavior of \HouseQR and its blocked variants is analyzed in various mixed-precision settings. The potential performance benefits of mixed precision make the study of the stability properties of mixed-precision BGS variants of great interest going forward. 

\section{Software Frameworks} \label{sec:software}
In the following, we highlight which variants of BGS are implemented in well-known software packages. 

To our knowledge, the Trilinos Library, in particular the Belos package, provides the most opportunity for the user to select among various orthogonalization schemes to be used within Krylov subspace methods. Belos
currently contains implementations of \TSQR, iterated CGS and MGS, and \CGSSRO. For block Krylov subspace methods, Belos implements various BGS methods. For example, in block GMRES, the user may choose between block versions of iterated CGS and MGS, and \CGSSRO, the default being block iterated CGS. A recent presentation \cite{YamazakiHoemmenBoman2019} on the Trilinos Library presents experiments with communication-avoiding (\s-step) GMRES using $\CGS \circ \CholQR$ as well as a ``low-sync CGS'' with \CholQRRO as the \IOnoarg. We also note that \s-step GMRES with low-sync block Gram-Schmidt orthogonalization algorithms including \BCGSPIP, \BMGSICWY, and \BCGSIROLS is described in the recent paper by Yamazaki et~al.  \cite{Yamazakietal2020}. The authors demonstrate significant parallel scaling improvements over the previous Trilinos implementations of these KSM algorithms and lower time-to-solution for the iterative solver. The implementation of communication-avoiding Krylov subspace methods and associated block orthogonalization schemes, in particular low-synch variants, is a goal of the ongoing Exascale Computing Project PEEKS effort.

Block variants of GMRES, CG, break-down free CG and recycling GCRO Krylov solvers have all been implemented recently in PETSc as described in the paper by Jolivet et~al.~\cite{Jolivetetal2020}. The authors also describe blocked eigensolvers employed within the PETSc software framework that are based on SLEPc. These employ the Locally Optimal Block Preconditioned Conjugate Gradient (LOBPCG) algorithm due to Knyazev \cite{Knyazev2001} and the Contour Integral Spectrum Slicing method (CISS) of Sakurai and Suguria \cite{Sakurai2003}.

For prototype algorithms and testing routines for stability conjectures, we direct the reader to the \MATLAB package developed in conjunction with this survey, hosted at \url{https://github.com/katlund/BlockStab}.

\section{Conclusions and outlook} \label{sec:conclusions}
The take-home lesson of stability studies is always the same, but often forgotten: be careful, because the finite precision universe can behave in exceptionally undesirable ways, especially if you have only analyzed your method in exact arithmetic.  At the same time, when an algorithm works well ``most of the time," there are probably good reasons, and we should be realistic about the statistical likelihood that something catastrophic may occur, especially if the algorithm is providing massive benefit in other ways.

Sometimes, however, stability may be a stringent requirement for a software pipeline.  This is especially true if eigenvalue estimates are needed anywhere \cite{Stewart2001}.  Lately, Krylov subspace methods, and particularly block Krylov subspace methods, form an important component in matrix function approximations, which often use eigenvalue estimations to accelerate convergence \cite{EiermannErnstGuettel2011,
FrommerGuettelSchweitzer2014a, FrommerGuettelSchweitzer2014b, FrommerLundSzyld2020, Guettel2013}. Other applications where reliable eigenvalue estimates are needed include Krylov subspace recycling and spectral deflation \cite{ParksDeSturlerMackey2006, ParksSoodhalterSzyld2017}.

It is also worth exploring whether the low-sync BGS variants -- \BCGSIROLS, \BMGSSVL, \BMGSLTS, \BMGSCWY, \BMGSICWY -- coupled with muscles such as \TSQR or \ShCholQRRORO, provide similar communication benefits as $\BCGS \circ \TSQR$ but with better stability properties.  This is especially relevant for \s-step KSMs, which rely on the low communication of \BCGS but can also suffer from its instability.

Although our focus has been on the stability properties of block orthogonalization routines, an overarching goal is to be able to say something about the backward stability of block (and $s$-step) variants of GMRES methods. 
Paige, Rozlo{\v z}n{\'i}k, and Strako{\v s} have proven that (non-block) GMRES with MGS orthogonalization is backward stable, despite loss of orthogonality of the Arnoldi basis vectors due to the finite precision MGS computation~\cite{PaigeRozloznikStrakos2006}. It remains an open problem, even in the non-blocked case, to determine the level to which orthogonality can be lost in the finite precision orthogonalization routine while still obtaining a backward stable GMRES solution. Further, theoretical treatment of the backward stability of block Arnoldi/GMRES methods is entirely lacking from the literature. We hope that the present work provides a path forward in this direction, and a detailed application of the BGS results compiled here to block GMRES and methods like $s$-step GMRES remains as future work.
The natural next step is an examination of the least-squares problem for block upper Hessenberg matrices. Work by Gutknecht and Schmelzer \cite{GutknechtSchmelzer2008} may prove relevant, as well as the low-rank update formulation devised in \cite{FrommerLundSzyld2020} for matrix functions.

\appendix

\section{Pseudocode for muscles} \label{sec:musc-alg}

\begin{algorithm}[htbp!]
	\caption{$[\vQ, R] = \CGS(\vX)$ \label{alg:CGS}}
	\begin{algorithmic}[1]
		\STATE{$\vQ = \vX$}
	\FOR{$k=1, \ldots, n$}
		\STATE{$R_{1:k-1,k} = \vQ_{:,1:k-1}^T \vX_{:,k}$}
		\STATE{$\vQ_{:,k} = X_{:,k} - \vQ_{:,1:k-1}R_{1:k-1,k}$}
		\STATE{$R_{k,k} = \Vert \vQ_{:,k} \Vert$}
		\STATE{$\vQ_{:,k} = \vQ_{:,k}/R_{k,k}$}
	\ENDFOR
		\RETURN{$\vQ$, $R$}
	\end{algorithmic}
\end{algorithm}

\begin{algorithm}[htbp!]
	\caption{$[\vQ, R] = \MGS(\vX)$ \label{alg:MGS}}
	\begin{algorithmic}[1]
		\STATE{$\vQ = \vX$}
	\FOR{$k=1, \ldots, n$}
		\STATE{$R_{k,k} = \Vert \vQ_{:,k} \Vert$}
		\STATE{$\vQ_{:,k} = \vQ_{:,k}/R_{k,k}$}
		\STATE{$R_{k,k+1:n} = Q_{:,k}^T Q_{:,k+1:n}$}
		\STATE{$Q_{:,k+1:n} = Q_{:,k+1:n} - Q_{:,k}R_{k,k+1:n}$}
	\ENDFOR
		\RETURN{$\vQ$, $R$}
	\end{algorithmic}
\end{algorithm}

\begin{algorithm}[htbp!]
	\caption{$[\vQ, \vR] = \CGSSROR(\vX, \rpltol)$ \label{alg:CGSSROR}}
	\begin{algorithmic}[1]
		\STATE{Allocate memory for $\vQ$ and $R$}
		\STATE{$[\vq_1, \sim, r_{11}] = \texttt{cgs\_step\_sror}(\vZero, \vx_1, \rpltol)$}
		\FOR{$k = 1, \ldots, s-1$}
		\STATE{$[\vq_{k+1}, R_{1:k,k+1}, r_{k+1,k+1}] = \texttt{cgs\_step\_sror}(\vQ_{1:k}, \vx_{k+1}, \rpltol)$}
		\ENDFOR
		\RETURN{$\vQ = [\vq_1, \ldots, \vq_s]$, $R = (r_{jk})$}
	\end{algorithmic}
\end{algorithm}

\begin{algorithm}[htbp!]
	\caption{$[\vQ, R, T] = \MGSSVL(\vX)$ \label{alg:MGSSVL}}
	\begin{algorithmic}[1]
		\STATE{Allocate memory for $\vQ$ and $R$}
		\STATE{$T = I_s$}
		\STATE{$r_{11} = \norm{\vx_1}$; $\vq_1 = \vx_1/r_{11}$}
		\FOR{$k = 1, \ldots, s-1$}
			\STATE{$\vw = \vx_{k+1}$}
			\STATE{$R_{1:k,k+1} = \diff{T_{1:k,1:k}^T} \big(\vQ_{1:k}^T \vw \big)$}
			\STATE{$\vw = \vw - \vQ_{1:k} R_{1:k,k+1}$}
			\STATE{$r_{k+1,k+1} = \norm{\vw}$}
			\STATE{$\vq_{k+1} = \vw/r_{k+1,k+1}$}
			\STATE{\diff{$T_{1:k,k+1} = -T_{1:k,1:k} \big(\vQ_{1:k}^T \vq_{k+1}\big)$}}
		\ENDFOR
		\RETURN{$\vQ = [\vq_1, \ldots, \vq_s]$, $R = (r_{jk})$, $T = (t_{jk})$}
	\end{algorithmic}
\end{algorithm}

\begin{algorithm}[htbp!]
	\caption{$[\vQ, R, T] = \MGSLTS(\vX)$ \label{alg:MGSLTS}}
	\begin{algorithmic}[1]
		\STATE{Allocate memory for $\vQ$ and $R$}
		\STATE{$T = I_s$}
		\STATE{$r_{11} = \norm{\vx_1}$; $\vq_1 = \vx_1/r_{11}$}
		\FOR{$k = 1, \ldots, s-1$}
			\STATE{$\vw = \vx_{k+1}$}
			\STATE{$R_{1:k,k+1} = \diff{T_{1:k,1:k}^\tinv} \big(\vQ_{1:k}^T \vw\big)$}
			\STATE{$\vw = \vw - \vQ_{1:k} R_{1:k,k+1}$}
			\STATE{$r_{k+1,k+1} = \norm{\vw}$}
			\STATE{$\vq_{k+1} = \vw/r_{k+1,k+1}$}
			\STATE{\diff{$T_{1:k,k+1} = \vQ_{1:k}^T \vq_{k+1}$}}
		\ENDFOR
		\RETURN{$\vQ = [\vq_1, \ldots, \vq_s]$, $R = (r_{jk})$, $T = (t_{jk})$}
	\end{algorithmic}
\end{algorithm}

\begin{algorithm}[htbp!]
	\caption{$[\vQ, R, T] = \MGSCWY(\vX)$ \label{alg:MGSCWY}}
	\begin{algorithmic}[1]
		\STATE{Allocate memory for $\vQ$ and $R$}
		\STATE{$T = I_s$}
		\STATE{$\vu = \vx_1$}
		\FOR{$k = 1, \ldots, s-1$}
			\STATE{$\vw = \vx_{k+1}$}
			\IF{$k = 1$}
				\STATE{$\begin{bmatrix} r_{k,k}^2 & \rho \end{bmatrix} = \vu^T [\vu \,\, \vw]$}
			\ELSIF{$k > 1$}
				\STATE{$\begin{bmatrix} \vt & \vr \\ r_{k,k}^2 & \rho \end{bmatrix} = [\vQ_{1:k-1} \,\, \vu]^T [\vu \,\, \vw]$}
				\STATE{$T_{1:k-1,k} = \diff{-T_{1:k-1,1:k-1} (\vt / r_{k,k})}$}
			\ENDIF
			\STATE{$R_{1:k,k+1} = \diff{T_{1:k,1:k}^T} \begin{bmatrix} \vr \\ \rho/ r_{k,k} \end{bmatrix}$}
			\STATE{$\vq_k = \vu / r_{k,k}$}
			\STATE{$\vu = \vw - \vQ_{:,1:k} R_{1:k,k+1}$}
		\ENDFOR
		\STATE{$r_{s,s} = \norm{\vu}$}
		\STATE{$\vq_s = \vu / r_{s,s}$}
		\RETURN{$\vQ = [\vq_1, \ldots, \vq_s]$, $R = (r_{jk})$, $T = (t_{jk})$}
	\end{algorithmic}
\end{algorithm}

\begin{algorithm}[htbp!]
	\caption{$[\vQ, R, T] = \MGSICWY(\vX)$ \label{alg:MGSICWY}}
	\begin{algorithmic}[1]
		\STATE{Allocate memory for $\vQ$ and $R$}
		\STATE{$T = I_s$}
		\STATE{$\vu = \vx_1$}
		\FOR{$k = 1, \ldots, s-1$}
			\STATE{$\vw = \vx_{k+1}$}
			\IF{$k = 1$}
				\STATE{$\begin{bmatrix} r_{k,k}^2 & \rho \end{bmatrix} = \vu^T [\vu \,\, \vw]$}
			\ELSIF{$k > 1$}
				\STATE{$\begin{bmatrix} \vt & \vr \\ r_{k,k}^2 & \rho \end{bmatrix} = [\vQ_{1:k-1} \,\, \vu]^T [\vu \,\, \vw]$}
				\STATE{$T_{1:k-1,k} = \diff{\vt / r_{k,k}}$}
			\ENDIF
			\STATE{$R_{1:k,k+1} = \diff{T_{1:k,1:k}^\tinv} \begin{bmatrix} \vr \\ \rho/ r_{k,k} \end{bmatrix}$}
			\STATE{$\vq_s = \vu / r_{s,s}$}
			\STATE{$\vu = \vw - \vQ_{:,1:k} R_{1:k,k+1}$}
		\ENDFOR
		\STATE{$r_{s,s} = \norm{\vu}$}
		\STATE{$\vq_s = \vu / r_{s,s}$}
		\RETURN{$\vQ = [\vq_1, \ldots, \vq_s]$, $R = (r_{jk})$, $T = (t_{jk})$}
	\end{algorithmic}
\end{algorithm}

\begin{algorithm}[htbp!]
	\caption{$[\vQ, R] = \CholQR(\vX)$ \label{alg:CholQR}}
	\begin{algorithmic}[1]
		\STATE{$X = \vX^T\vX$}
		\IF{$X$ is positive definite}
		\STATE{$R = \chol(X)$}
		\STATE{$\vQ = \vX R^\inv$}
		\ELSE
		\RETURN{$\vQ = \nan$, $R = \nan$}
		\ENDIF
		\RETURN{$\vQ$, $R$}
	\end{algorithmic}
\end{algorithm}

\begin{algorithm}[htbp!]
	\caption{$[\vQ, R] = \CholQRRO(\vX)$ \label{alg:CholQRRO}}
	\begin{algorithmic}[1]
		\STATE{$[\vQ, R^{(1)}] = \CholQR(\vX)$}
		\STATE{$[\vQ, R^{(2)}] = \CholQR(\vQ)$}
		\STATE{$R = R^{(2)} R^{(1)}$}
		\RETURN{$\vQ$, $R$}
	\end{algorithmic}
\end{algorithm}

\begin{algorithm}[htbp!]
	\caption{$[\vQ, R] = \ShCholQRRORO(\vX)$ \label{alg:ShCholQRRORO}}
	\begin{algorithmic}[1]
		\STATE{Define shift $\sigma = 11(ms + s(s+1))\eps\norm{\vX}_2^2$}
		\STATE{$R^{(1)} = \chol(\vX^T \vX + \sigma I_s)$}
		\STATE{$\vQ = \vX (R^{(1)})^\inv$}
		\STATE{$[\vQ, R^{(2)}] = \CholQRRO(\vQ)$}
		\STATE{$R = R^{(2)} R^{(1)}$}
	\end{algorithmic}
\end{algorithm}

\FloatBarrier

\section{\MATLAB code for routines in \BCGSSROR and \CGSSROR} \label{sec:sror}
\FloatBarrier
\begin{small} \label{alg:cgs_step_sror}
	\begin{verbatim}
		function [y, r, rho, northog] = cgs_step_sror(Q, x, nu, rpltol)
		% [y, r, rho, northog] = CGS_STEP_SROR(Q, x, nu, rpltol) orthogonalizes x
		% against the columns of Q using the the Classical Gram-Schmidt method.
		% Reorthogonalization is performed as necessary to ensure orthogonality.
		% Specifically, given an orthonormal matrix Q and a vector x, the function
		% returns a vectors y and r, and a scalar rho satisfying
		% 
		% (*)   x = Q*r + rho*y
		% 
		% where y is a normalized vector orthogonal to the column space of Q.
		% The matrix Q may have zero columns.
		% 
		% The optional argument nu is the norm of the original value of x, to be
		% used when x has been subject to previous orthogonalizaton steps.  If
		% nu is absent or is less than or equal to norm(x), is is set to norm(x).
		% 
		% The optional argument rlptol controls when the current vector y is
		% replaced by a random vector.  Its default value is 1.  If it is set
		% to a value greater than one, the relation (*) will be compromised
		% somewhat, but the number of orthogonalizations may be decreased.
		% 
		% The optional output argument northog, if present, contains the number
		% of orthogonalizations.
		%
		% Originally written by G. W. Stewart, 2008. Modified by Kathryn Lund,
		% 2020.
		
		%%
		% Initialize
		[n, nq] = size(Q);
		r = zeros(nq, 1);
		nux = norm(x);
		if nargout >= 4
		northog = 0;
		end
		
		% If Q has no columns, return normalized x if x~=0, otherwise return a
		% normalized random vector.
		if nq == 0
		if nux == 0
		y = rand(n,1)-0.5;
		y = y/norm(y);
		rho = 0;
		else
		y = x/nux;
		rho = nux;
		end
		return
		end
		
		%  If required, intitalize the optional arguments nu and rpltol.
		if nargin < 3
		nu = nux;
		else
		if nu < nux
		nu = nux;
		end
		end
		
		if nargin < 4
		rpltol = 1;
		end
		
		% If norm(x)==0, set it to a nonzero vector and proceed, noting the fact in
		% the flag zeronorm.
		if nux ~= 0
		zeronorm = false;
		y = x/nux;
		nu = nu/nux;
		else
		zeronorm = true;
		y = rand(n,1) - 0.5;
		y = y/norm(y);
		nu = 1;
		end
		
		% Main orthogonalization loop.
		nu1 = nu;
		while true
		if nargout == 4
		northog = northog + 1;
		end
		s = Q'*y;
		r = r + s;
		y = y - Q*s;
		nu2 = norm(y);
		
		% Return if y is orthogonal.
		if nu2 > 0.5*nu1
		break
		end
		
		% Continue orthogonalizing if nu2 is not too small.
		if nu2 > rpltol*nu*eps
		nu1 = nu2;
		else % Replace y by a random vector.
		nu = nu*eps;
		nu1 = nu;
		y = rand(n,1) - 0.5;
		y = nu*y/norm(y);
		end
		end
		
		% Calculate rho and normalize y.
		if ~zeronorm
		rho = norm(y);
		y = y/rho;
		rho = rho*nux;
		r = r*nux;
		else
		y = y/norm(y);
		r = zeros(nq, 1);
		rho = 0;
		end
	\end{verbatim}
\end{small}

\FloatBarrier
\begin{small} \label{alg:bcgs_step_sror}
	\begin{verbatim}
		function [Y, R12, R22] = bcgs_step_sror(QQ, X, rpltol)
		% [Y, R12, R22] = BCGS_STEP_SROR(QQ, X, rpltol) returns an orthonormal
		% matrix Y whose columns lie in the orthogonal complement of the column
		% space of QQ, a matrix R12, and an upper triangular matrix R22 satisfying
		% 
		% (*)   X = QQ*R12 + Y*R22.
		% 
		% The optional argument rpltol has a default value of 1.  Increasing it,
		% say to 100, may improve performance, but will degrade the accuracy of the
		% relation (*).
		%
		% Originally written by G. W. Stewart, 2008.  Modified by Kathryn Lund,
		% 2020.
		
		%%
		%  Initialization.
		reorth = false;
		if nargin < 3
			rpltol = 1;
		else
				if rpltol < 1
						rpltol = 1;
				end
		end

		s = size(X,2);
		sk = size(QQ,2);

		nu = zeros(1,s);
		for k = 1:s
			nu(k) = norm(X(:,k));
		end

		%  Beginning of the first orthogonalization step.  Project Y
		%  onto the orthogonal complement of Q.
		R12 = QQ'*X;
		Y = X - QQ*R12;
		R22 = zeros(s);

		%  Orthogonalize the columns of Y.
		for k = 1:s
				[yk, r, rho] = cgs_step_sror(Y(:,1:k-1), Y(:,k), nu(k), rpltol);
				Y(:,k) = yk;
				R22(1:k-1,k) = r;
				R22(k,k) = rho;
				if rho <= 0.5*nu(k)
						reorth = true;
				end
		end

		%  Return if Q has zero columns.
		if sk == 0 || ~reorth
			 return
		end

		% Beginning of the reorthogonalization.  Project Y onto the orthogonal
		% complement of Q.
		S12 = QQ'*Y;
		Y = Y - QQ*S12;

		% Orthogonalize the columns of Y.
		S22 = zeros(s);
		for k = 1:s
				[yk, r, sig] = cgs_step_sror(Y(:,1:k-1), Y(:,k), 1, rpltol);

				if sig < 0.5
				% The result, yk, is not satisfactorily orthogonal.  Perform an
				% unblocked orthogonalization against Q and the previous columns of Y.
						[yk, r, sig] = cgs_step_sror([QQ, Y(:,1:k-1)], Y(:,k), rpltol);
						S12(:,k) = S12(:,k) + r(1:sk);
						r = r(sk+1:sk+k-1);
				end
				Y(:,k) = yk;
				S22(1:k-1,k) = r;
				S22(k,k) = sig;
		end

		R12 = R12 + S12*R22;
		R22 = S22*R22;
		end
	\end{verbatim}
\end{small}

\section{Glued matrices and command-line calls} \label{sec:scripts}
Glued matrices are generated by the following \MATLAB code, modified from \cite{SmoktunowiczBarlowLangou2006}:
\vspace{10pt}
\begin{small}
	\begin{verbatim}
	function A = CreateGluedMatrix(m, p, s, r, t)
	% Example 2 matrix from [Smoktunowicz et al., 2006]. Generates a glued
	% matrix of size m x ps, with parameters r and t specifying the powers of
	% the largest condition number of the first stage and second stages of the
	% matrix, respectively:
	%
	% Stage 1: A = U * Sigma * V'
	%
	% Stage 2: A = A * kron(I, Sigma_block) * kron(I, V_block)
	
	%%
	n = p * s;
	U = orth(randn(m,n));
	V = orth(randn(n,m));
	Sigma = diag(logspace(0, r, n));
	A = U * Sigma * V';
	
	Sigma_block = diag(logspace(0, t, s));
	V_block = orth(randn(s,s));
	
	ind = 1:s;
	for i = 1:p
		A(:,ind) = A(:,ind) * Sigma_block * V_block';
		ind = ind + s;
	end
	end
	\end{verbatim}
\end{small}
\vspace{10pt}

All the heat-maps in Section~\ref{sec:bgs_stab_1} can be reproduced with the code hosted at \linebreak \url{https://github.com/katlund/BlockStab} and the following commands:
\vspace{10pt}
\begin{verbatim}
XXdim = [10000 50 10];
mat = {'rand_uniform', 'rand_normal', 'rank_def', 'laeuchli',...
    'monomial', 'stewart', 'stewart_extreme', 'hilbert',...
    's-step', 'newton'};
skel = {'BCGS', 'BCGS_IRO', 'BCGS_SROR', 'BCGS_IRO_LS',...
    'BMGS', 'BMGS_SVL', 'BMGS_CWY'};
musc = {'CGS', 'CGS_IRO', 'CGS_SRO', 'CGS_SROR',...
    'CGS_IRO_LS', 'MGS', 'MGS_SVL', 'MGS_CWY', 'HouseQR',...
    'CholQR', 'CholQR_RO', 'Sh_CholQR_RORO'};
rpltol = 100;
verbose = 1;
MakeHeatmap(XXdim, mat, skel, musc, rpltol, verbose)
\end{verbatim}
\vspace{10pt}

All the $\kappa$-plots in Section~\ref{sec:bgs_stab_2} can also be easily reproduced; see Table~\ref{tab:code}.
\begin{table}[H]
	\caption{Command-line calls for $\kappa$-plots \label{tab:code} }
	\begin{tabular}{c|l}
		Figure										& Code \\ \hline
		\ref{fig:bcgs_p_kappa_plot}					& \verb|GluedBlockKappaPlot([1000 50 4], 1:8, ... | \\
													& \verb|{'BCGS', 'BCGS_PIP', 'BCGS_PIO'}, ... | \\
													& \verb|{'CholQR', 'HouseQR'})| \\ \hline
		\ref{fig:bcgs_iro_kappa_plot} 				& \verb|BlockKappaPlot([100 20 2], -(1:16), ... | \\
													& \verb|    {'BCGS', 'BCGS_IRO'}, ...| \\
													& \verb|    {'CGS', 'MGS', 'HouseQR'})| \\ \hline
		\ref{fig:bcgs_iro_laeuchli} 				& \verb|LaeuchliBlockKappaPlot([1000 100 5], ... |\\
													& \verb|    logspace(-1,-16,10), ...| \\
													& \verb|    {'BCGS', 'BCGS_IRO'}, ...| \\
													& \verb|    {'CGS', 'MGS', 'HouseQR'})| \\ \hline											
		\ref{fig:bcgs_iro_ls_kappa_plot} 			& \verb|BlockKappaPlot([100 20 2], -(1:16), ...| \\
													& \verb|    {'BCGS', 'BCGS_IRO', 'BCGS_IRO_LS'}, ...| \\
													& \verb|    {'CGS', 'CGS_IRO', 'CGS_IRO_LS'})| \\ \hline
		\ref{fig:bcgs_iro_ls_laeuchli_kappa_plot}	& \verb|LaeuchliBlockKappaPlot([1000 120 2], ... |\\
													& \verb|		logspace(-1,-16,10), ...| \\
													& \verb|    {'BCGS', 'BCGS_IRO', 'BCGS_IRO_LS'}, ...| \\
													& \verb|    {'CGS', 'CGS_IRO', 'CGS_IRO_LS'})| \\ \hline
		\ref{fig:bmgs_kappa_plot} 					& \verb|BlockKappaPlot([100 20 2], -(1:16), 'BMGS', ...| \\
													& \verb|    {'CGS', 'MGS', 'CGS_RO', ...| \\
													& \verb|    'CholQR', 'Cholqr_RO', 'HouseQR'})| \\
													& \verb|LaeuchliBlockKappaPlot([1000 120 2], ...| \\
													& \verb| 			logspace(-1,-16,10), 'BMGS', ...| \\
													& \verb|    {'CGS', 'MGS', 'CGS_RO', ...| \\
													& \verb|    'CholQR', 'Cholqr_RO', 'HouseQR'})| \\ \hline
		\ref{fig:bmgs_t_kappa_plot}					& \verb|BlockKappaPlot([100 20 2], -(1:16), ...| \\
													& \verb|    {'BMGS_SVL', 'BMGS_CWY'}, ...| \\
													& \verb|    {'MGS_SVL', 'MGS_LTS', ...|\\
													& \verb|    'MGS_CWY', 'MGS_ICWY', 'HouseQR'})| \\
													& \verb|LaeuchliBlockKappaPlot([1000 120 2], ...| \\
													& \verb|    logspace(-1,-16,10), ...| \\
													& \verb|    {'BMGS_SVL', 'BMGS_CWY'}, ...| \\
													& \verb|    {'MGS_SVL', 'MGS_LTS', ...|\\
													& \verb|    'MGS_CWY', 'MGS_ICWY', 'HouseQR'})|
	\end{tabular}
\end{table}
\FloatBarrier


\bibliographystyle{elsarticle-num.bst}
\bibliography{back_stab}
\end{document}